\documentclass[10pt, reqno]{amsart}

\newcommand{\adsg}[1]{\textcolor{red}{Antoine: #1}}

\newcommand{\red}[1]{\textcolor{red}{#1}}

\usepackage{extarrows}
\usepackage{stmaryrd}
\usepackage{MnSymbol}
\numberwithin{equation}{section}

\usepackage{amsxtra}
\usepackage{hyperref}
\usepackage{tikz-cd}
\usepackage[all]{xy}
\usetikzlibrary{patterns}
\usetikzlibrary{backgrounds}

\let\oldtocsection=\tocsection

\let\oldtocsubsection=\tocsubsection

\let\oldtocsubsubsection=\tocsubsubsection

\renewcommand{\tocsection}[2]{\hspace{0em}\oldtocsection{#1}{#2}}
\renewcommand{\tocsubsection}[2]{\hspace{2em}\oldtocsubsection{#1}{#2}}
\renewcommand{\tocsubsubsection}[2]{\hspace{4.5em}\oldtocsubsubsection{#1}{#2}}

\makeatletter
\def\subsection{\@startsection{subsection}{2}
 \z@{.5\linespacing\@plus.7\linespacing}{-.5em}
 {\normalfont\bfseries}}
\def\subsubsection{\@startsection{subsubsection}{3}
 \z@{.5\linespacing\@plus.7\linespacing}{-.5em}
 {\normalfont\bfseries}}
\makeatother

\theoremstyle{plain}

\newtheorem{theorem}{Theorem}[section]
\newtheorem{lemma}[theorem]{Lemma}
\newtheorem{definition-theorem}[theorem]{Definition-Theorem}
\newtheorem{proposition}[theorem]{Proposition}
\newtheorem{corollary}[theorem]{Corollary}
\newtheorem{definition}[theorem]{Definition}
\newtheorem{example}[theorem]{Example}
\newtheorem{remark}[theorem]{Remark}
\newtheorem{notation}[theorem]{Notation}
\newtheorem{assumption}[theorem]{Assumption}
\newtheorem{lemma-definition}[theorem]{Lemma-Definition}
\newtheorem{lemma-notation}[theorem]{Lemma-Notation}
\newtheorem{question}[theorem]{Question}
\newtheorem{remark-definition}[theorem]{Remark-Definition}
\newtheorem{notation-remark}[theorem]{Notation-Remarks}
\newtheorem{conjecture}[theorem]{Conjecture}
\newtheorem{notation-Lemma}[theorem]{Notation-Lemma}
\newcommand \bth[1] { \begin{theorem}\label{t#1} }
\newcommand \ble[1] { \begin{lemma}\label{l#1} }

\newcommand \bpr[1] { \begin{proposition}\label{p#1} }
\newcommand \bco[1] { \begin{corollary}\label{c#1} }
\newcommand \bconj[1] { \begin{conjecture}\label{co#1} }
\newcommand \bde[1] { \begin{definition}\label{d#1}\rm }
\newcommand \bex[1] { \begin{example}\label{e#1}\rm }
\newcommand \bre[1] { \begin{remark}\label{r#1}\rm }
\newcommand \bnota[1] {\begin{notation}\label{n#1}\rm }
\newcommand \bas[1] { \begin{assumption}\label{a#1}\rm }
\newcommand \bld[1] { \begin{lemma-definition}\label{ld#1} }

\newcommand \bqu[1] { \begin{question}\label{q#1}\rm }

\newcommand {\eth} { \end{theorem} }
\newcommand {\ele} { \end{lemma} }

\newcommand {\epr} { \end{proposition} }
\newcommand {\econj} { \end{conjecture} }
\newcommand {\eco} { \end{corollary} }
\newcommand {\ede} { \end{definition} }
\newcommand {\eex} { \end{example} }
\newcommand {\ere} { \end{remark} }
\newcommand {\enota} { \end{notation} }
\newcommand {\eas} {\end{assumption}}
\newcommand {\eld}{ \end{lemma-definition} }

\newcommand {\equ} {\end{question}}
\newcommand \thref[1]{Theorem \ref{t#1}}
\newcommand \leref[1]{Lemma \ref{l#1}}
\newcommand \prref[1]{Proposition \ref{p#1}}
\newcommand \coref[1]{Corollary \ref{c#1}}

\newcommand \deref[1]{Definition \ref{d#1}}
\newcommand \exref[1]{Example \ref{e#1}}
\newcommand \reref[1]{Remark \ref{r#1}}

\newcommand \ldref[1]{Lemma-Definition \ref{ld#1}}

\def \RR {{\mathbb R}}         

\def \ZZ {{\mathbb Z}}

\def \TT {{\mathbb T}}

\def \QQ {{\mathbb Q}}
\def \PP {{\mathbb P}}

\def \calP {{\mathcal{P}}}

\def \wt {\widetilde}
\def \wh {\widehat}

\def \h  {\mathfrak{h}}

\def\mod{\opname{mod}\nolimits}

\renewcommand \max { {\mathrm{max}} }

\newcommand{\opname}[1]{\operatorname{\mathsf{#1}}}

\newcommand\preceqdot{\mathrel{\ooalign{$\preceq$\cr
  \hidewidth\raise0.225ex\hbox{$\cdot\mkern0.5mu$}\cr}}}

\newcommand{\beqa}{\begin{eqnarray*}}                     
\newcommand{\eeqa}{\end{eqnarray*}}
\def \hs {\hspace{.2in}}

\def \wh {\widehat}

\def \sA {{\scriptscriptstyle A}}

\def \FF {\mathbb{F}}

\def \PP {\mathbb{P}}
\def \YY {\mathbb{Y}}
\def \AA {\mathbb{A}}

\def \calC {{\mathcal{C}}}
\def \calP {{\mathcal{P}}}
\def \calU {{\mathcal{U}}}

\def \bb {\backslash\!\!\backslash}

\def \bfx {{\bf x}}
\def \bfy {{\bf y}}

\def \calD {{\mathcal{D}}}

\def \sv {{\scriptscriptstyle \vee}}
\def \tkt{{\begin{xy}(0,1)*+{t}="A",(10,1)*+{t'}="B",\ar@{-}^k"A";"B" \end{xy}}}

\def \tit {{\begin{xy}(0,1)*+{t}="A",(10,1)*+{t'}="B",\ar@{-}^i"A";"B" \end{xy}}}
\def \vkv {\begin{xy}(0,1)*+{v}="A",(10,1)*+{v'}="B",\ar@{-}^k"A";"B" \end{xy}}

\def \val {{\rm val}}
\def \XX {\mathbb{X}}
\def \Zmax {\ZZ^{\rm max}}
\def \bfA {{\bf A}}
\def \bfY {{\bf Y}}

\def \pB {p_{\!{}_{B}}}

\def \fX {\mathfrak{X}}
\def \OO {\mathbb{O}}
\def \pbreveB {p_{\!{}_{\breve{B}}}}
\def \pbarBT {p_{\!{}_{\overline{B}^T}}}
\def \Aprin {\AA_B^{\rm prin}}
\def \Yprin {\YY_{B^T}^{\rm prin}}
\def \ABT {\AA_{B^T}}

\def \ABTZ {\AA_{B^T}(\Zmax)}
\def \YB {\YY_B}
\def \YBZ {\YY_B(\Zmax)}

\def \sdeltasv {{\scriptstyle{\delta^\sv}}}

\def \CfriezeAT {C_{\rm frieze}(A^T)}
\def \CaddA {C_{\rm add}(A)}
\def \sw {{\scriptscriptstyle \wedge}}

\begin{document}

\setlength{\baselineskip}{1.1\baselineskip}
\vspace{-.35in}
\title[Tropical friezes and cluster-additive functions via Fock-Goncharov duality]
{Tropical friezes and cluster-additive functions via Fock-Goncharov duality and a conjecture of Ringel}
\author{Peigen Cao}
\address{
School of Mathematical Sciences, University of Science and Technology of China,
Hefei 230026, Anhui, P. R. China}
\email{peigencao@126.com}
\author{Antoine de St. Germain}
\address{
Department of Mathematics   \\
The University of Hong Kong \\
Pokfulam Road               \\
Hong Kong}
\email{adsg96@hku.hk}
\author{Jiang-Hua Lu}
\address{
Department of Mathematics   \\
The University of Hong Kong \\
Pokfulam Road               \\
Hong Kong}
\email{jhlu@maths.hku.hk}
\date{}

\begin{abstract}
We study tropical friezes and cluster-additive functions
associated to symmetrizable generalized Cartan matrices in the framework of Fock-Goncharov duality in cluster algebras.
In particular, we generalize and prove a conjecture of C. M. Ringel on cluster-additive functions associated to arbitrary
Cartan matrices of finite type. For a Fock-Goncharov dual pair of positive spaces of finite type,
we use tropical friezes and cluster-additive functions to explicitly express  
the Fock-Goncharov pairing between their tropical points and the bijections between  global monomials on one 
 and tropical points of the other.

\end{abstract}

\maketitle
\tableofcontents
\addtocontents{toc}{\protect\setcounter{tocdepth}{1}}

\section{Introduction}
\subsection{Introduction}\label{ss:intro}
Additive functions for a translation quiver
$\Gamma = (\Gamma_0, \Gamma_1, \tau)$, defined as $\ZZ$-valued functions on the
vertex set $\Gamma_0$ and satisfying certain mesh-type conditions,
play an important role in the
representation theory of finite dimensional algebras \cite{HPR_1980,bobinski:additive-fcts},  \cite[p.52]{dlab-gabriel:rep1}. When $\Gamma$ is the repetition quiver of an acyclic valued quiver, the definition of additive functions can be formulated using the underlying Cartan matrix.

\begin{definition}\label{:cluster-additive}
{\rm 
Let $A = (a_{i,j})$ be any $r \times r$ symmetrizable generalized Cartan matrix
(see $\S$\ref{ss:thmA-intro}). An {\it additive function associated to $A$} is a map $d: [1,r]\times \ZZ \to \ZZ$ satisfying
\begin{equation}\label{eq:add}
 	d(i,m) + d(i,m+1) = \sum_{j=i+1}^r (-a_{j,i}) d(j,m) + \sum_{j=1}^{i-1}(-a_{j,i}) d(j,m+1), \;\;\;
  (i,m) \in [1, r] \times \ZZ.
\end{equation}
Here and for the rest of the paper, we write $[1, r] = \{1, 2, \ldots, r\}$ for any integer $r \geq 1$.
}
\end{definition}

Motivated by the combinatorial features of cluster categories \cite{bmrrt_2006} and cluster-tilted algebras \cite{bmr-2007},
C. M. Ringel introduced in \cite[$\S$11]{Ringel:additive} the notion of {\it cluster-additive functions}. At roughly the same time, L. Guo introduced in \cite{Guo:tropical-frieze} (see also \cite{Propp:integers}) {\it tropical friezes} in the context of 
2-Calabi-Yau Hom-finite triangulated categories. We reformulate their definitions again using Cartan matrices as follows. For $a \in \RR$, set $[a]_+ = \max(0,a)$.

\bde{:intro}
{\rm Let  $A = (a_{i,j})$ be any $r \times r$ symmetrizable generalized  Cartan matrix.

1) A {\it cluster-additive function associated to  $A$} is a map
$k: {[1,r]} \times \ZZ \to \ZZ$ satisfying
\begin{equation}\label{eq:cl-additive}
k(i,m) + k(i,m+1) = \sum_{j=i+1}^r (-a_{j,i}) [k(j,m)]_+ + \sum_{j=1}^{i-1}(-a_{j,i}) [k(j,m+1)]_+, \;\;\;
(i, m) \in [1, r]\times \ZZ.
\end{equation}

2) A {\it tropical frieze associated to  $A$} is a map
$f: {[1,r]} \times \ZZ \to \ZZ$ satisfying 
\begin{equation}\label{eq:trop-frieze}
f(i,m) + f(i,m+1) = \left[\sum_{j=i+1}^r (-a_{j,i}) f(j,m) + \sum_{j=1}^{i-1}(-a_{j,i}) f(j,m+1) \right]_+, \;\;\;
(i, m) \in [1, r]\times \ZZ.
\end{equation}
}
\ede
Let $C_{\rm add}(A)$, resp. $C_{\rm frieze}(A)$, be the set of all 
cluster-additive functions, resp. tropical friezes, associated to $A$.
As \eqref{eq:cl-additive} and  \eqref{eq:trop-frieze} are recursive with respect to the total order 
on $[1, r] \times \ZZ$ given by
\begin{equation}\label{eq:order-0}
(i, m) < (i', m') \hs \mbox{iff}\hs m <  m', \;\; \mbox{or} \; \;m = m' \;\; \mbox{and} \; \;i< i',
\end{equation}
one has bijections  $C_{\rm add}(A) \cong \ZZ^r$ and $C_{\rm frieze}(A) \cong \ZZ^r$ 
via evaluations on $[1, r] \times \{0\}$.

Suppose that $A$ is of finite and simply-laced type, and identify $[1,r] \times \ZZ$ with the 
vertex set of the Auslander-Reiten quiver of the
derived category $D^b(\mod \Lambda)$, where $\Lambda$ is the path algebra of the Dynkin quiver associated to $A$.
Ringel \cite[p. 477]{Ringel:additive} conjectured a {\it periodicity property} of every
$k \in C_{\rm add}(A)$, which amounts to saying that $k$ can be regarded as a cluster-additive function on 
the Auslander-Reiten quiver of the {\it cluster category} of $\mod \Lambda$, 
i.e., the orbital quotient of $D^b(\mod \Lambda)$ by its Frobenius isomorphism \cite{bmrrt_2006}.
For every  (not necessarily basic) tilting module $T$ of $\Lambda$, Ringel also constructed in \cite[$\S$10]{Ringel:additive}
a function $d_T \in C_{\rm add}(A)$ expressed via dimensions of indecomposable modules of the 
cluster-tilted algebra of $T$, and gave a decomposition of $d_T$ into 
certain elementary cluster-additive functions called {\it cluster-hammock functions} 
(see \S \ref{ss:hammock}). The main conjecture of Ringel in \cite{Ringel:additive}, stated in 
\cite[$\S$6]{Ringel:additive}, says that {\it every $k \in C_{\rm add}(A)$ is a finite sum of 
cluster-hammock functions}. Ringel's conjecture is thus the cluster categorical analog of 
Butler's theorem \cite{Butler-1981} that additive functions on the Auslander-Reiten quiver of a representation-finite 
algebra are integral linear combinations of hammock functions. 
Ringel's conjecture on cluster-additive functions was proved by Ringel \cite{Ringel:additive} for type $A_n$
and by Guo \cite{Guo:tropical-frieze} for all (finite and simply-laced) type. See \reref{:Guo} for
Guo's work \cite{Guo:tropical-frieze} on tropical friezes.  

In this paper, we study tropical friezes and cluster-additive functions associated to 
{\it arbitrary symmetrizable generalized Cartan matrices}, and, more importantly, in contrast to the categorical approaches of Ringel \cite{Ringel:additive} and Guo 
\cite{Guo:tropical-frieze}, 
we  do so via tropical evaluations of functions on certain Fock-Goncharov dual pairs of positive spaces \cite{FG:ensembles}, 
thereby placing tropical friezes and cluster-additive functions in the framework of Fock-Goncharov duality in the theory of cluster algebras \cite{FG:ensembles, GHKK}. 

To achieve the said goal,  for a Fock-Goncharov dual pair $(\AA, \YY)$ as defined in \cite{FG:ensembles} (see \deref{:ensemble-pair}),
instead of looking for parametrizations of vector space bases of the algebra of (global) functions on $\AA$ (resp. on $\YY)$ by 
the set of tropical points of $\YY$ (resp. of $\AA$) as in the original conjectures stated in \cite{FG:ensembles}, 
we introduce a {\it correspondence} between (global) functions on $\AA$ (resp.
on $\YY$) and tropical points of $\YY$ (resp. of $\AA)$ via the notion of {\it admissible functions}. For 
an arbitrary symmetrizable generalized Cartan matrix $A$, we give two
realizations, of both tropical friezes and cluster-additive functions, via tropical points of one space and admissible functions on the other in a suitable Fock-Goncharov dual pair, and we show that
the two realizations coincide if and only if the tropical points and the admissible functions are in correspondence (see Theorem A in $\S$\ref{ss:thmA-intro} and Theorem C in $\S$\ref{ss:thmC-intro}). In particular, the values of 
 tropical friezes and cluster-additive functions are interpreted via denominator vectors of admissible functions.

The connections between tropical friezes, cluster-additive functions and Fock-Goncharov duality are in full 
display when the Cartan matrix $A$ is of {\it finite type}: 
as a consequence of the validity of the Fock-Goncharov duality
for cluster algebras of finite type,  we 
formulate and prove Ringel's conjecture on cluster-additive functions associated to a Cartan matrix $A$ {\it of any finite type}
 (i.e., not necessarily simply-laced as in the original conjecture of Ringel in \cite{Ringel:additive}),
and we determine, in Lie theoretical terms, the periodicity of all cluster-additive functions and tropical friezes 
associated to such an $A$. 
In the other direction, for any Fock-Goncharov dual pair $(\AA, \YY)$ of finite type, we show that admissible functions 
on $\AA$ and $\YY$ are precisely all the global monomials \cite[Definition 0.1]{GHKK},
and we use tropical friezes and cluster-additive functions to 
express, in formulas that can be
readily implemented by a computer program,
the bijection between tropical points of $\YY$ (resp. of $\AA$) and 
global monomials on $\AA$ (resp. on $\YY$) as well as 
the Fock-Goncharov pairing between tropical points of $\AA$ and $\YY$.

We give more precise statements of the main results of the paper in $\S$\ref{s:main}.

\subsection{Acknowledgments}
The three authors  have been partially supported by grants from  the Research Grants Council of the Hong Kong SAR, 
China (GRF 17307718, GRF 17306621, GRF 17303420). P. Cao has also been supported 
 the National Natural Science Foundation of China Grant No. 12071422, the Guangdong Basic and Applied Basic Research Foundation Grant No. 2021A1515012035, and  the 
JSPS Postdoctoral Fellowship for Research in Japan (Standard) between September 2023 and April 2024. 

\section{Main results}\label{s:main}

\subsection{Tropical friezes and cluster-additive functions via tropical points}\label{ss:thmA-intro}
The starting point of our paper is a simple observation that both 
\eqref{eq:cl-additive} and \eqref{eq:trop-frieze} are tropicalizations of certain relations in 
a field of rational functions. To see this, we first fix some notation to be used throughout the paper.

\bnota{:trop-intro}
{\rm
By $\ZZ^{\max}$ we mean the semi-field $(\ZZ, \cdot, \oplus)$ with 
$a \cdot b = a+b$ and $a \oplus b = {\rm max}(a, b)$ for $a, b \in \ZZ$.
For any semi-field $\PP$, let ${\rm Hom}_{\rm sf}(\PP, \Zmax)$ be the set of all semi-field homomorphisms from $\PP$ to $\Zmax$.
For $r$ independent variables ${\bf x} = (x_1, \ldots, x_r)$,  let $\QQ({\bf x})_{>0}$ be the 
semi-field of rational functions in ${\bf x}$ over $\QQ$ 
that have subtraction-free expressions in ${\bf x}$.  We also write
\[
\val_\delta(x) = \delta(x) \hs \;\;\; \mbox{for}\;\;\; \delta \in {\rm Hom}_{\rm sf}(\QQ({\bf x})_{>0}, \ZZ^{\max})
\;\; \mbox{and} \;\; x \in \QQ({\bf x})_{>0}.
\]
}
\enota

Let  $A = (a_{i, j})$ be any $r \times r$ symmetrizable generalized Cartan matrix, i.e, an integral $r \times r$ matrix such that
$a_{i, i} = 2$ for all $i \in [1, r]$, $a_{i, j} \leq 0$ for all $i \neq j$,
and $DA$ is symmetric for some diagonal matrix $D$ with positive integers on the diagonal. Consider two families 
$\{x(i,m)\}_{(i, m) \in [1,r] \times \ZZ}$ and $\{y(i,m)\}_{(i, m) \in [1,r] \times \ZZ}$
of commuting variables, respectively satisfying, for all $(i, m) \in [1, r] \times \ZZ$,
\begin{align}\label{eq:x(i,m)s}
x(i,m) \, x(i,m+1) &= 1 + \prod_{j=i+1}^r x(j,m)^{-a_{j,i}} \; \prod_{j=1}^{i-1} x(j,m+1)^{-a_{j,i}},\\
\label{eq:y(i,m)s}
y(i,m) \, y(i,m+1) &= \prod_{j=i+1}^r (1+y(j,m))^{-a_{j,i}} \; \prod_{j=1}^{i-1} (1+y(j,m+1))^{-a_{j,i}}.
\end{align}
As 
\eqref{eq:x(i,m)s} and \eqref{eq:y(i,m)s} are again recursive with respect to the total order in \eqref{eq:order-0},
setting
\[
{\bf x} = (x(1, 0), \ldots, x(r, 0)) \hs \mbox{and} \hs  {\bf y} = (y(1, 0), \ldots, y(r, 0))
\]
as independent variables, one can respectively regard 
$\{x(i,m)\}_{(i, m) \in [1,r] \times \ZZ}$ and $\{y(i,m)\}_{(i, m) \in [1,r] \times \ZZ}$
as systems of elements in $\QQ({\bf x})_{>0}$ and
$\QQ({\bf y})_{>0}$.

For $\delta \in {\rm Hom}_{\rm sf}(\QQ({\bf x})_{>0}, \ZZ^{\max})$ and 
$\rho \in {\rm Hom}_{\rm sf}(\QQ({\bf y})_{>0}, \ZZ^\max)$, define
\begin{align*}
&f_\delta:\; [1, r] \times \ZZ \longrightarrow \ZZ, \;\; f_\delta(i, m) = \val_\delta (x(i, m)),\\
&k_\rho:\; [1, r] \times \ZZ \longrightarrow \ZZ, \;\; k_\rho(i, m) = \val_\rho (y(i, m)).
\end{align*}

\ble{:trop-intro}
For any symmetrizable generalized Cartan matrix $A$, one has bijections
\begin{align}
\label{eq:frieze-bijection}
&{\rm Hom}_{\rm sf}(\QQ({\bf x})_{>0},\; \ZZ^\max) \longrightarrow C_{\rm frieze}(A), \;\; \delta \longmapsto  f_\delta,\\
\label{eq:additive-bijection}
&{\rm Hom}_{\rm sf}(\QQ({\bf y})_{>0},\; \ZZ^\max) \longrightarrow C_{\rm add}(A), \;\; \rho \longmapsto  k_\rho.
\end{align}
\ele

\begin{proof}
It follows from \eqref{eq:x(i,m)s} that $f_\delta\in C_{\rm frieze}(A)$ for every 
$\delta \in {\rm Hom}_{\rm sf}(\QQ({\bf x})_{>0}, \ZZ^\max)$.
As an element in ${\rm Hom}_{\rm sf}(\QQ({\bf x})_{>0}, \ZZ^\max)$ is uniquely determined by its values on ${\bf x}$,
\eqref{eq:frieze-bijection} is a bijection. Similarly, \eqref{eq:additive-bijection} is well-defined and is a bijection.
\end{proof}

Tropicalization is best formulated via the notion of positive spaces, as introduced in \cite[\S 1.1]{FG:ensembles}. 

\bde{:pos-space}
A {\it positive space} of dimension $r$ 
is a pair 
\[
\XX = (\FF, \,\{{\bf x}_t\}_{t \in \TT}),
\]
where $\FF$ is a field isomorphic to the field 
of  rational functions over $\QQ$ in $r$ independent variables, $\TT$ is any non-empty index set,
and $\{{\bf x}_t\}_{t \in \TT}$ is an assignment of an ordered
set of free generators of $\FF$ over $\QQ$ to each $t \in \TT$, such that $\QQ({\bf x}_t)_{>0} = 
\QQ({\bf x}_{t'})_{>0} \subset \FF$ for all $t, t' \in \TT$. We also write $\FF = \FF(\XX)$.
\hfill $\diamond$
\ede

For a positive space $\XX = (\FF(\XX), \{{\bf x}_t\}_{t \in \TT})$, 
we call each ${\bf x}_t$, $t \in \TT$, a {\it coordinate chart}, or a {\it chart}, of $\XX$.  Define
{\it the ambient semi-field} of $\XX$ as $\FF_{>0}(\XX) = \QQ({\bf x}_t)_{>0}$ using any $t \in \TT$, and set
\[
\XX(\Zmax) = {\rm Hom}_{\rm sf}(\FF_{>0}(\XX), \, \Zmax).
\]
Elements in $\XX(\Zmax)$ are called {\it $\Zmax$-tropical points}, or simply {\it tropical points}, of $\XX$. 
Note that one can define a tropical point $\delta$ of $\XX$ by specifying its
{\it coordinate vector} $(\val_\delta(x_{t; 1}), \ldots, \val_\delta(x_{t; r}))$ using {\it any} chart ${\bf x}_t = 
(x_{t; 1}, \ldots, x_{t; r})$ of $\XX$. 

The theory of cluster algebras famously gives rise to two classes of positive spaces.
Indeed, recall that a skew-symmetrizable integral matrix is also called a {\it mutation matrix}. Associated to an $r \times r$ mutation matrix $B$, 
one has two $r$-dimensional positive spaces
\[
\AA_B = (\FF(\AA_B), \; \{{\bf x}_t\}_{t \in \TT_r}) \hs \mbox{and} \hs \YY_B = (\FF(\YY_B), \; \{{\bf y}_t\}_{t \in \TT_r}),
\]
where $\TT_r$ is a rooted $r$-regular tree,  $\{{\bf x}_t\}_{t \in \TT_r}$ is the collection of 
{\it $\bfA$-clusters}, pairwise related by sequences of $\bfA$-variable (i.e., cluster variable)  mutations 
(see \eqref{eq:X-seed-mutation-cluster}), 
and $\{{\bf y}_t\}_{t \in \TT_r}$ is the collection of 
{\it $\bfY$-clusters}, pairwise related by sequences of $\bfY$-variable mutations (see \eqref{eq:Y-mut}).
The sub-algebras
\[
\calU(\AA_B) = \bigcap_{t \in \TT_r} \ZZ[{\bf x}_t^{\pm 1}] \subset \FF(\AA_B) \hs \mbox{and} \hs 
\calU(\YY_B) = \bigcap_{t \in \TT_r} \ZZ[{\bf y}_t^{\pm 1}] \subset \FF(\YY_B)
\]
are respectively called the {\it upper cluster algebras} (with trivial coefficients) of $\AA_B$ and of $\YY_B$, and they are the central objects of study in the 
theory of cluster algebras \cite{BFZ:ClusterIII, FZ:ClusterIV}.
By respectively tropicalizing the $\bfA$-variable and $\bfY$-variable mutation rules,  one obtains the
respective mutation rules for the coordinate vectors of tropical points of $\AA_B$ and  tropical points of $\YB$.  See 
\leref{:construct-trop-point-A} 
and \leref{:construct-trop-point}.

\bde{:ensemble-pair} \cite{FG:ensembles}
For any $r \times r$ mutation matrix $B$, 

1) the pair $(\AA_B, \YY_B)$ of positive spaces is called a {\it Fock-Goncharov ensemble}; and

2) the pair $(\AA_{B^T}, \YY_B)$ of positive spaces is called a {\it Fock-Goncharov dual pair}. 

\noindent
Here and for the rest of the paper, $M^T$ for a matrix $M$ denotes the transpose of $M$.
\ede

For an $r \times r$ symmetrizable generalized Cartan matrix $A = (a_{i, j})$, we define
the mutation matrix
\begin{equation}\label{eq:BA-intro}
	 B = B_A = \begin{pmatrix}
		0 & a_{1,2} & a_{1,3} & \dots & a_{1,r} \\
		-a_{2,1} & 0 & a_{2,3} & \dots & a_{2,r} \\
		-a_{3,1} & -a_{3,2} & 0 & \dots & a_{3,r} \\
		\dots& \dots & \dots & \dots & \dots \\
		-a_{r,1} & -a_{r,2} & -a_{r,3}& \dots &0
	\end{pmatrix}. 
\end{equation}
Note that $B_A$ is {\em acyclic} in the sense of \cite[Definition 1.14]{BFZ:ClusterIII}, and every acyclic mutation matrix is, 
up to simultaneous re-labelling of rows and columns, of the form 
$B_A$ for a symmetrizable generalized Cartan matrix $A$ \cite[p.36]{BFZ:ClusterIII}.

The starting point of our paper is the fact that relations in \eqref{eq:x(i,m)s} and \eqref{eq:y(i,m)s} are respectively
satisfied by 
certain $\bfA$-variables and certain $\bfY$-variables in the Fock-Goncharov dual pair
\begin{equation}\label{eq:OA-1}
\OO_A \, \stackrel{\rm def}{=}\, (\ABT, \; \YB).
\end{equation}

To give more detail, 
let $t_0$ be the root vertex of $\TT_r$ and label by $t(i, m)$, $(i, m) \in [1, r] \times \ZZ_{\geq 0}$,
the vertices on the infinite path in $\TT_r$ that starts from $t_0$ and follows
the sequence of edges $(1, 2, \ldots, r, 1, 2, \ldots)$ via
\begin{equation}\label{eq:right}
\begin{xy}
(-6,1)*+{t_0 = t(1, 0)}="A2",(14,1)*+{t(2, 0)}="A3", (28,1)*+{\cdots}="A4",(42,1)*+{t(r, 0)}="A5",
(58,1)*+{t(1,1)}="A6",(74,1)*+{t(2, 1)}="A7", (88,1)*+{\cdots}="A8",
\ar@{-}^{\;\;\;\;1}"A2";"A3", \ar@{-}^{\;\;\;\;2}"A3";"A4",
\ar@{-}^{r-1\;\;}"A4";"A5",\ar@{-}^{r}"A5";"A6",\ar@{-}^{1}"A6";"A7",\ar@{-}^{\;\;\;2}"A7";"A8",
\end{xy}.
\end{equation}
Similarly, label by $t(i, m)$, $(i, m) \in [1, r] \times \ZZ_{\leq -1}$,
the vertices on the infinite path in $\TT_r$ that starts from $t_0$ and
follows the sequence of edges $(r, r-1, \ldots,2,  1, r, r-1, \ldots)$ 
via
\begin{equation}\label{eq:left}
\begin{xy}
(27,1)*+{\cdots}="A4",(42.5,1)*+{t(r, -2)}="A5",
(60.5,1)*+{t(1,-1)}="A6",(78.5,1)*+{t(2, -1)}="A7", (93.5,1)*+{\cdots}="A8",(108.5,1)*+{t(r,-1)}="A9",
(129,1)*+{t(1, 0) = t_0}="A10",
\ar@{-}^{r-1\;\;\;\;}"A4";"A5",\ar@{-}^{r}"A5";"A6",\ar@{-}^{1}"A6";"A7",
\ar@{-}^{\;\;\;\;2}"A7";"A8",\ar@{-}^{r-1\;\;\;\;}"A8";"A9",\ar@{-}^{r\;\;}"A9";"A10"
\end{xy}.
\end{equation}

\bnota{:xim-yim-intro}
{\rm
Denote the $\bfA$-clusters of $\ABT$ by $\{\bfx_t^\sv\}_{t \in \TT_r}$ and the $\bfY$-clusters of
$\YB$ by $\{\bfy_t\}_{t \in \TT_r}$. 
Let $\fX(\ABT)$ be the set of all cluster variables of $\ABT$ 
and $\fX(\YY_B)$ that of all ${\bf Y}$-variables of $\YY_B$. For $(i, m) \in [1, r] \times \ZZ$, 
let $x^\sv(i, m) \in \fX(\ABT)$ be the $i^{\rm th}$ variable in 
${\bf x}^\sv_{t(i, m)}$, and  let $y(i, m) \in \fX(\YB)$
be the  $i^{\rm th}$ variable in ${\bf y}_{t(i, m)}$. Consider the maps
\begin{align}\label{eq:xsv-im-intro}
&[1, r] \times \ZZ \longrightarrow \fX(\ABT), \;\; (i, m) \longmapsto x^\sv(i, m),\\
\label{eq:yim-intro}
&[1, r] \times \ZZ \longrightarrow \fX(\YY_{B}), \;\; (i, m) \longmapsto y(i, m).
\end{align}
For $\delta^\sv \in \ABT(\Zmax)$ and $\rho \in \YY_B(\Zmax)$, let $(\delta^\sv_{t; 1}, \ldots, \delta^\sv_{t; r})$ 
and $(\rho_{t; 1}, \ldots, \rho_{t; r})$ be the respective coordinate vectors of 
$\delta^\sv$ and $\rho$ at $t \in \TT_r$, and define $f_{\sdeltasv}, k_\rho: [1, r] \times \ZZ \to \ZZ$ by
\begin{align}\label{eq:f-delta-intro}
&f_{\sdeltasv}(i, m) = \val_{\delta^\sv}(x^\sv(i, m)) = \delta^\sv_{t(i, m); i}, \hs \hs (i, m) \in [1, r] \times \ZZ,\\
\label{eq:k-rho-intro}
&k_\rho(i,m) = \val_\rho(y(i, m)) = \rho_{t(i, m); i}, \hs \hs (i, m) \in [1, r] \times \ZZ.
\end{align}
}
\enota

The following is a combination of 
\prref{:xim-yim}, \prref{:trop-frieze-dual-1}, and \prref{:additive-dual-1}.

\medskip
\noindent
{\bf Theorem A.}
{\it For any symmetrizable generalized Cartan matrix $A$ and with $B = B_A$ as in \eqref{eq:BA-intro}, the 
map in \eqref{eq:xsv-im-intro} satisfies \eqref{eq:x(i,m)s} with $A$ replaced by $A^T$, i.e., 
\begin{equation}\label{eq:relations-xsv-intro-00}
x^{{\sv}}(i,m) \, x^{{\sv}}(i,m+1) = 1 + \prod_{j=i+1}^r x^{\sv}(j,m)^{-a_{i, j}} \; 
 \prod_{j=1}^{i-1} x^{\sv}(j,m+1)^{-a_{i, j}}, \hs (i, m) \in [1, r] \times \ZZ,
\end{equation}
while the map in 
\eqref{eq:yim-intro} satisfies \eqref{eq:y(i,m)s}.
Consequently, one has the well-defined bijections
\[
\ABT(\Zmax) \longrightarrow C_{\rm frieze}(A^T), \;\; \delta^\sv \longmapsto  f_{\sdeltasv}, \hs \mbox{and} \hs 
\YY_B(\Zmax) \longrightarrow C_{\rm add}(A), \;\; \rho \longmapsto  k_\rho.
\]
}

\medskip
We refer to the bijections in Theorem A as {\it realizations} of tropical friezes
and cluster-additive functions {\it via tropical points}
of the Fock-Goncharov dual pair $(\ABT, \YB)$. 

The map in \eqref{eq:xsv-im-intro}, called the 
{\it generic $\bfA$-frieze pattern}\footnote{The  
recursive procedure in \eqref{eq:relations-xsv-intro-00}
of assigning cluster variables of $\ABT$ is also known as the {\it knitting algorithm} 
\cite{keller-knitting-algo}, and the seeds of $\ABT$
(see \deref{:cluster-pattern})  at $t(i, m)\in \TT_r$ for $(i, m) \in [1, r] \times \ZZ$ belong to the {\it acyclic belt} of $\ABT$ (see \cite{RSW1}).}, has been studied in 
\cite{ARS:friezes, GHL:friezes}.  Correspondingly, we call the map in 
\eqref{eq:yim-intro} the
{\it generic $\bfY$-frieze pattern}.
In view of Theorem A, tropical friezes are tropicalizations of generic $\bfA$-frieze patterns, 
while cluster-additive functions are 
tropicalizations of generic $\bfY$-frieze patterns. 
We can thus also call cluster-additive functions {\it tropical $\bfY$-friezes}.

Theorem A sets the stage for the main investigation in  our paper, namely that of the interplay between tropical friezes,
cluster-additive functions, and Fock-Goncharov duality. For this purpose, we devote $\S$\ref{s:admissible-elements} 
and $\S$\ref{s:FG-pairing} to a general discussion, summarized in the next $\S$\ref{ss:F-G-intro}, 
on certain aspects of Fock-Goncharov duality via the notion of {\it admissible functions} on Fock-Goncharov dual pairs.

\subsection{Fock-Goncharov duality via admissible functions}\label{ss:F-G-intro}


Let $B$ be any $r \times r$ mutation matrix, not necessarily of the form $B_A$ in \eqref{eq:BA-intro}, 
and consider the Fock-Goncharov
dual pair $(\ABT, \YB)$. 

Following \cite[Definition 0.1]{GHKK}, for $\XX = \ABT$ or $\YB$, elements in $\calU(\XX)$ that are
Laurent monomials in some cluster of $\XX$  are called
{\it global monomials} on $\XX$. Global monomials on $\ABT$ are precisely the cluster monomials on $\ABT$ defined 
in \cite{fz_2002} (see \leref{:A-global-mono}). By \cite{GHKK} (see \leref{:Delta-plus} and \ldref{:Y-g-vector}), global monomials on
$\ABT$ are in bijection with their $g$-vectors which lie in $\YBZ$, and global monomials on
$\YB$ are in bijection with their $g$-vectors which lie in $\ABTZ$.

V. V. Fock and A. B. Goncharov conjectured in \cite[\S 4]{FG:ensembles} that 
there is a $\ZZ$-basis $\{x^\sv_\rho\}$ of $\calU(\ABT)$, parameterized by $\rho \in \YBZ$ and extending the
parametrization of cluster monomials on $\ABT$ by their $g$-vectors, and a $\ZZ$-basis $\{y_{\scriptstyle{\delta^\sv}}\}$
of $\calU(\YB)$, parameterized by $\delta^\sv \in \ABT(\Zmax)$ and extending the parametrization of 
global monomials on $\YB$ by their $g$-vectors, such that 
\begin{equation}\label{eq:FG-val-1}
\val_{\delta^\sv}(x^\sv_\rho) = \val_\rho(y_{\sdeltasv}), \hs \;\;\forall \; 
\delta^\sv \in \ABT(\Zmax), \;  \; \rho \in \YBZ.
\end{equation}
The common integer value in \eqref{eq:FG-val-1} is the conjectured {\it Fock-Goncharov pairing} between the tropical points
$\delta^\sv \in \ABT(\Zmax)$ and $\rho \in \YBZ$.
The conjecture as stated does not hold in general \cite{GHK}, as $\calU(\YY_B)$ may have ``too few" elements. In \cite{GHKK} the authors modified the 
conjecture by constructing the so-called {\it theta functions}
$\theta_\rho \in \calU(\ABT)$ and $\theta_{\delta^\sv}
\in \calU(\YY_B)$
for $\rho$ in a certain subset of
$\YBZ$ and $\delta^\sv$ in a certain subset of $\ABT(\Zmax)$, and it is conjectured in \cite[Remark 9.11]{GHKK}
that \eqref{eq:FG-val-1} holds 
for $x^\sv_\rho = \theta_\rho$  and $y_{\scriptstyle{\delta^\sv}}= \theta_{\delta^\sv}$. When $B$ is skew-symmetric,
the conjecture in 
\cite[Remark 9.11]{GHKK} has been proved in \cite{Muller_talk} (see  \reref{:theta-pairing}). 

While theta functions are important, their construction in \cite{GHKK} requires the rather involved 
techniques of scattering diagrams and broken lines in tropical spaces.  
In place of theta functions, in this paper we introduce and work with the so-called
{\it admissible functions}. More precisely, for any $\rho \in \YBZ$, 
an element $x^\sv \in \calU(\ABT)$ is said to be {\it $\rho$-admissible}, and we write  
\begin{equation}\label{eq:rho-admi}
x^\sv \in \calU^{\rm admi}_\rho(\ABT),
\end{equation}
if the Laurent expansion of $x^\sv$ in every cluster
${\bf x}_t^\sv$ of $\ABT$ has a certain prescribed form 
involving the coordinate vector of $\rho$ at $t$  (see \deref{:admissible-A} for detail). 
An element $x^\sv \in \calU(\ABT)$ is called an {\it admissible function} on $\ABT$ if it is $\rho$-admissible for some $\rho \in \YBZ$.
Similarly, for $\delta^\sv \in \ABTZ$ we can define
$\delta^\sv$-admissible functions on $\YB$ (see \deref{:admissible-Y}), and we write
\begin{equation}\label{eq:delta-admi}
y \in \calU^{\rm admi}_{\delta^\sv}(\YB)
\end{equation}
if $y \in \calU(\YB)$ is $\delta^\sv$-admissible. An element in $\calU(\YB)$ is called an {\it admissible function} 
on $\YB$ if
it is $\delta^\sv$-admissible for some  $\delta^\sv \in \ABTZ$. Admissible functions, 
while still capturing the essential properties of theta functions, are more general (see $\S$\ref{ss:theta-admi} and \exref{:unique}) 
and easier to work with 
due to the descriptive nature of
their definition.

We regard relations expressed 
in \eqref{eq:rho-admi} and \eqref{eq:delta-admi} as {\it correspondences} between tropical points in one positive space 
and global functions on the other in the dual pair $(\ABT, \YB)$. The correspondences are in general not bijective
(see \exref{:unique}).
Inspired by \cite[Conjectures 4.3]{FG:ensembles} of Fock and Goncharov, we ask the following natural question 
(a similar question has been studied by J. Fei \cite{Fei19b}).

\medskip
\noindent
{\bf Question.} (Question \ref{q:FG-pairing-conj}) {\rm Let $B$ be any $r \times r$ mutation matrix. When does one have
\begin{equation}\label{eq:val-val-intro}
\val_{\delta^\sv}(x^\sv) = \val_\rho(y),
\end{equation}
where $\delta^\sv \in \AA_{B^T}(\Zmax)$,  $\rho \in \YY_B(\Zmax)$, 
 $x^\sv \in \calU^{\rm admi}_\rho(\AA_{B^T})$, and $y \in \calU^{\rm admi}_{\delta^\sv}(\YY_{B})$?
}

\medskip
For the purpose of this paper, it is enough to consider admissible functions on positive 
$\bfA$-spaces and positive $\bfY$-spaces with trivial coefficients (except in 
$\S$\ref{ss:admissible-Aprin} where we also consider principal coefficients). 
When the mutation matrices are of full rank, admissible functions are precisely the {\it universally positive and compatibly pointed
elements} in the sense of F. Qin \cite{qin-2022} (see \reref{:admissible-1}). Admissible functions
on positive spaces with frozen variables (and full rank extended mutation matrices) are used by the first author in \cite{cao:F-invariant} to study {\it $F$-invariants}. 

In $\S$\ref{s:admissible-elements}, we review the bijections between global monomials on $\ABT$ and on $\YB$ with 
their $g$-vectors. For a cluster monomial $x^\sv$  on $\ABT$, let  $\rho(x^\sv) \in \YBZ$ be its $g$-vector,
and for a global monomial $y$ on $\YB$ let $\delta^\sv(y) \in \ABTZ$ be its $g$-vector.
The following is a summary of our results on admissible functions.

\medskip
\noindent
{\bf Theorem B.} {\it Let $B$ be an $r \times r$ mutation matrix and 
consider the Fock-Goncharov dual pair $(\AA_{B^T}, \YY_{B})$.

1) The following statements hold.

\hspace{3mm} 1a) \ldref{:globals-A-are-admissible}: a cluster monomial $x^\sv$  on $\ABT$ is, and only is, $\rho(x^\sv)$-admissible.

\hspace{3mm} 1b) \ldref{:Y-g-vector}: a global monomial $y$ on $\YB$ is, and only is, $\delta^\sv(y)$-admissible.

\hspace{3mm} 1c) $\S$\ref{ss:theta-admi}: theta functions on $\ABT$ and on $\YB$ are admissible.

2) \coref{:admissible-exist} and \coref{:acyclic-rho-unique}: if $B$ is mutation equivalent to an acyclic mutation matrix, 
then $\rho$-admissible functions on $\ABT$ exist for every 
$\rho \in \YBZ$;
an admissible function on $\ABT$ is $\rho$-admissible for a unique $\rho \in \YBZ$; and 
$\delta^\sv$-admissible functions on $\YB$ exist for every 
$\delta^\sv \in \ABTZ$ (by definition an admissible function on $\YB$ is $\delta^\sv$-admissible for a unique 
$\delta^\sv \!\in\! \ABTZ$, see \reref{:delta-unique}).

3) \thref{:finite-admi-mono}: When $B$ is of finite type, admissible functions on $\ABT$ 
 are precisely all the cluster monomials,   
and admissible functions on $\YB$ are precisely all the global monomials.
}

\smallskip
In \prref{:duality-pairing} we show that identity \eqref{eq:val-val-intro} holds 
if either $x^\sv$ or $y$ is a global monomial.

\subsection{Tropical friezes and cluster-additive functions via admissible functions}\label{ss:thmC-intro}
Let $A$ be any $r \times r$ symmetrizable generalized Cartan matrix, let $B = B_A$ as in \eqref{eq:BA-intro}, and
consider again the
Fock-Goncharov dual pair
$\OO_A \, = \,  (\AA_{B^T}, \;\YY_{B})$.
Recall from \eqref{eq:xsv-im-intro} and \eqref{eq:yim-intro}
the generic frieze patterns 
 \[
 (i, m) \longmapsto x^\sv(i, m) \in \fX(\ABT) \hs \mbox{and} \hs
 (i, m) \longmapsto y(i, m) \in \fX(\YB).
 \]
We show in \leref{:B-im-column} 
that each $y(i, m)$ is a {\it global variable} of $\YB$, i.e., a
$\bfY$-variable of $\YB$ that is also a global monomial. 
For $(i, m) \in [1, r] \times \ZZ$, denote the $g$-vectors of $x^\sv(i, m)$ and $y(i, m)$ respectively by
\[
\rho(i,m) \in \YBZ \hs \mbox{and} \hs 
\delta^\sv(i, m) \in \ABTZ.
\]
For
$x^\sv \in \calU^{\rm admi}(\ABT)$ and $y \in \calU^{\rm admi}(\YB)$, define
$f_y, k_{x^\sv}: [1, r]\times \ZZ \to \ZZ$ by 
\begin{equation}\label{eq:f-y-k-x-intro}
f_y(i, m) = \val_{\rho(i,m)}(y) \hs \mbox{and} \hs 
k_{x^\sv}(i,m) = \val_{\delta^\sv(i,m)}(x^\sv), \hs (i, m) \in [1, r]\times \ZZ.
\end{equation}
Recall from Theorem A  the bijections 
\[
\ABTZ \longrightarrow \CfriezeAT, \;\; \delta^\sv \longmapsto f_{\sdeltasv}, \hs \mbox{and} \hs
\YBZ \longrightarrow \CaddA, \;\; \rho \longmapsto k_\rho.
\]
The following Theorem C is a combination of 
\prref{:trop-frieze-dual-2} and \prref{:additive-dual-2}.

\medskip
\noindent
{\bf Theorem C.} {\it Let $A$ be any symmetrizable generalized Cartan matrix and let $B = B_A$ as in \eqref{eq:BA-intro}.

 1) One has $f_{y} \in C_{\rm frieze}(A^T)$ for every $y\in \calU^{\rm admi}(\YB)$, and one has a surjective map
\[
\calU^{\rm admi}(\YB) \longrightarrow C_{\rm frieze}(A^T), \;\; y \longmapsto f_{y}.
\]
Moreover, for $y \in \calU^{\rm admi}(\YB)$ and $\delta^\sv \in \ABTZ$, one has $f_{y} = f_{\sdeltasv}$
 if and only if $y$ is $\delta^\sv$-admissible. 

2) One has $k_{x^\sv} \in C_{\rm add}(A)$ for every $x^\sv \in \calU^{\rm admi}(\ABT)$,  and one has a surjective map
\begin{equation}\label{eq:k-xsv-1}
\calU^{\rm admi}(\ABT) \longrightarrow C_{\rm add}(A), \;\; x^\sv \longmapsto k_{x^\sv}.
\end{equation}
Moreover, for $x^\sv \in \calU^{\rm admi}(\ABT)$ and $\rho \in \YBZ$, one has $k_{x^\sv} = k_\rho$
 if and only if $x^\sv$ is $\rho$-admissible. 
}

\smallskip
Combining Theorem A and Theorem C, we have thus two realizations of each $f \in C_{\rm frieze}(A^T)$ and $k \in \CaddA$
using the Fock-Goncharov dual pair $\OO_A = (\ABT, \YB)$, namely 
\begin{align}\label{eq:f-dual}
&f = f_{\sdeltasv} = f_{y} \hs \mbox{for a unique}\;\, \delta^\sv \in \ABTZ \; \,\mbox{and any} \;\,
y \in \calU^{\rm admi}_{\delta^\sv}(\YB),\\
\label{eq:k-dual}
&k = k_\rho = k_{x^\sv}  \hs \mbox{for a unique}\;\, \rho \in \YBZ \; \,\mbox{and any} \;\,
x^\sv \in \calU^{\rm admi}_\rho(\ABT).
\end{align}
\bde{:slice}
{\rm
For a map $g: [1, r]\times \ZZ \to \ZZ$ and for any $m \in \ZZ$, define the {\it $m^{\rm th}$ slice} of $g$ to be 
\[
g|_{\{m\}} = (g(1, m), \, \ldots, \, g(r, m))^T \in \ZZ^r.
\]
}
\ede

We have the following  consequences of Theorem C.

\medskip
\noindent
{\bf Theorem C'}. {\it Let $A$ be any symmetrizable generalized Cartan matrix and let $B = B_A$ as in \eqref{eq:BA-intro}.

1) \leref{:f-y-alter}: For $y \in \calU^{\rm admi}(\YB), f = f_y \in C_{\rm frieze}(A^T)$ and $m \in \ZZ$, 
$f|_{\{m\}}$ is the denominator vector of $y$ at $t(1, m)$;

2) \coref{:kxsv-slice}: For $x^\sv \in \calU^{\rm admi}(\ABT), k = k_{x^\sv} \in \CaddA$ and $m \in \ZZ$,
 $k|_{\{m\}}$ is the 
denominator vector of $x^\sv$ at $t(1, m)$;

3) \leref{:trop-frieze-decomposition}: If $y_1, y_2 \in \calU^{\rm admi}(\YB)$ are such that $y_1 y_2 \in \calU^{\rm admi}(\YB)$, then 
$f_{y_1 y_2} = f_{y_1} + f_{y_2}$;

4) \leref{:ca-decomposition}: If $x_1^\sv, x_2^\sv \in \calU^{\rm admi}(\AA_{B^T})$ are such that $x_1^\sv x_2^\sv \in \calU^{\rm admi}(\AA_{B^T})$, then 
$k_{x_1^\sv x_2^\sv} = k_{x_1^\sv} + k_{x_2^\sv}$.
}

\medskip
When the Cartan matrix $A$ is of finite and simply-laced type, for each $(i, m) \in [1, r] \times \ZZ$, Ringel defined in 
\cite{Ringel:additive} the {\it cluster-hammock function} $h_{(i, m)}$ to be the unique element in $C_{\rm add}(A)$ satisfying
\begin{equation}\label{eq:him-intro}
h_{(i, m)} (j, m) = -\delta_{j, i},  \hs j \in [1, r].
\end{equation}
Defining $h_{(i, m)} \in C_{\rm add}(A)$ as in \eqref{eq:him-intro} 
for {\it arbitrary symmetrizable generalized Cartan matrix $A$}, we show in 
\prref{:hammock-1} and \coref{:hammock-3} that for every $(i, m) \in [1, r] \times \ZZ$, 
\begin{equation}\label{eq:him-intro-1}
h_{(i, m)} = k_{x^\sv(i, m)}  \in C_{\rm add}(A) \cap C_{\rm frieze}(A).
\end{equation}
As a consequence of \eqref{eq:him-intro-1}, we generalize in \prref{:d-dual} 
a result of N. Reading and S. Stella
\cite[Theorem 2.2 and Theorem 2.3]{Reading-Stella_2018} on a duality of $d$-compatibility degrees in the finite type and rank $2$ cases. 
 

\subsection{The case of finite type: Ringel's conjecture and applications}\label{ss:finite-intro}
In $\S$\ref{s:finite} we apply the results summarized in $\S$\ref{ss:thmA-intro}-$\S$\ref{ss:thmC-intro} to 
the case when the Cartan matrices are of finite type. 

Assume thus that $A$ is an $r \times r$ Cartan matrix of finite type and 
let $B = B_A$ as in \eqref{eq:BA-intro}.
We show in \prref{:F-invariant-frieze} that both generic frieze patterns
\[
(i, m) \longmapsto x^\sv(i, m)\in \fX(\ABT) \hs \mbox{and} \hs (i, m) \longmapsto y(i, m) \in \fX(\YB)
\]
are invariant under a bijection 
${\mathfrak{F}}_A: [1, r] \times \ZZ \rightarrow [1, r] \times \ZZ$, which is described using the root system of $A$
and generalizes the {\it gliding symmetry} of the classical Coxeter frieze patterns 
(see $\S$\ref{ss:Lie} and \reref{:FA}).
Combining  with results in 
Theorem A, Theorem C, and 4) of Theorem C', we prove the following generalization of 
Ringel's conjecture which is stated in \cite[$\S$6]{Ringel:additive} for the simply-laced case.

\medskip
\noindent
{\bf Theorem D.} [Ringel's conjecture] {\it For every Cartan matrix $A$ of finite type, the following hold.

1) Periodicity property: 
all $f \in \CfriezeAT$ and all $k \in \CaddA$ are ${\mathfrak{F}}_A$-invariant.

2) Decomposition property: every $k \in C_{\rm add}(A)$ is a (finite) sum of cluster-hammock functions.
}

\medskip
Decomposition of tropical friezes will depend on a better understanding of the fan structure of the tropical space $\ABTZ$ 
(but see \leref{:trop-frieze-decomposition}) and will be 
investigated in a future paper.

Turning in the other direction, since $A$ is of finite type, 
the full Fock-Goncharov conjectures hold
for the Fock-Goncharov pair $\OO_A = (\ABT, \YB)$. Namely, one has the bijective Fock-Goncharov correspondences
\begin{align}\label{eq:deltasv-to-y}
&\ABTZ \longrightarrow \calU^{\rm mono}(\YY_B): \;\;\; 
\delta^\sv\; \longmapsto \;y_{\sdeltasv},\\
\label{eq:rho-to-xsv}
&\YBZ \longrightarrow  \calU^{\rm mono}(\AA_{B^T}): \;\;\;\rho\; \longmapsto \;x_{\rho}^\sv,
\end{align}
where $y_{\sdeltasv}$ has $g$-vector $\delta^\sv$ and $x_{\rho}^\sv$ has $g$-vector $\rho$, as well as the 
Fock-Goncharov pairing
\begin{equation}\label{eq:pairing-intro}
\langle \delta^\sv, \; \rho\rangle = \val_{\delta^\sv} (x^\sv_\rho) = \val_\rho(y_{\sdeltasv}), \hs \delta^\sv \in \ABTZ, \; \,
\rho \in \YBZ.
\end{equation}
Let $\calD_A \subset [1, r]\times \ZZ$ be the fundamental domain of ${\mathfrak{F}}_A$ given in \eqref{eq:calD}, so one has the bijection
\[
\calD_A \longrightarrow \fX(\ABT), \;\; (i, m) \longmapsto x^\sv(i,m).
\]
In particular, $|\calD_A| = r + \frac{1}{2} |\Phi_A|$, where $\Phi_A$ is the 
root system of $A$. For $\delta^\sv \in \ABTZ$ and $\rho \in \YBZ$, let $f_{\sdeltasv} \in \CfriezeAT$ and  
$k_\rho \in \CaddA$ be as in Theorem A. 
As applications of tropical friezes and cluster-additive functions to Fock-Goncharov duality, we prove in 
$\S$\ref{ss:thmE}  the following Theorem E.

\medskip
\noindent
{\bf Theorem E.} {\it Let $A$ be a Cartan matrix of finite type, and let $B = B_A$ as in \eqref{eq:BA-intro}.

1) For every $\rho \in \YBZ$, one has 
\[
x^\sv_\rho= \prod_{(i,m) \in \calD_A} x^\sv (i,m)^{[-k_\rho(i,m)]_+};
\]

2) For $\delta^\sv \in \ABTZ$ and $\rho \in \YBZ$, one has
\[
    \langle \delta^\sv , \rho \rangle = \sum_{(i,m) \in \calD_A} f_{\sdeltasv}(i,m) [-k_{\rho}(i,m)]_+;
\]

3) For $\delta^\sv \in \ABTZ$ with (row) coordinate vector 
$\delta^\sv_{t_0}$ at $t_0$, letting $-p(\delta^\sv)$ be the unique point in $\YY_{-B^T}(\Zmax)$ with coordinate vector  
$-\delta^\sv_{t_0} B^T$ at $t_0$, one has
$k_{-p(\delta^\sv)} \in C_{\rm add}(A^T)$ and
\begin{equation}\label{eq:y-deltasv}
y_{\sdeltasv} = 
{\bf y}_{t_0}^{-(\delta_{t_0}^\sv)^T} \prod_{(i, m) \in \calD_A} \left(F(i, m)({\bf y}_{t_0})\right)^{[-k_{-p(\delta^\sv)}(i, m-1)]_+},
\end{equation}
where ${\bf y}_{t_0}$ is the $\bfY$-cluster of $\YB$ at $t_0$, and 
$F(i, m)$ for 
$(i, m) \in \calD_A$ is the $i^{\rm th}$ F-polynomial at $t(i, m)$ 
associated to the mutation matrix $B$ with principal coefficients at $t_0$.
}

\medskip

We remark that a geometric model of $\ABT$ is given in \cite{YZ:Lcc} using certain Coxeter reduced 
double Bruhat cells in the complex semi-simple Lie group associated to $A^T$, 
and the cluster variables $x^\sv(i, m)$ and the $F$-polynomials $F(i, m)$,
for $(i, m) \in \calD_A$, can be explicitly expressed using generalized minors.

We also remark that 
C. Melo and A. N\'ajera Ch\'avez obtain in \cite[Lemma 4.1]{MeloChavez:Afan} a formula expressing 
each global monomial $y_{\sdeltasv}$ on $\YB$ in terms of $\delta_{t_0}^\sv$ and 
$F$-polynomials from a vertex $w \in \TT_r$ at which  $\delta^\sv$ is optimized
(see \S \ref{ss:a-question} for the definition). 

The formulas in 
Theorem E can be easily implemented using a computer program.  
Some of it can be done using Keller's Java applets \cite{Keller-Java}. A more direct program is available at the
{\it Visual Cluster Algebra} website \cite{visualca} maintained by the second author.

\subsection{Notation}\label{ss:nota-intro} For an integer $r \geq 1$, we denote by $\ZZ^r$ 
the set of all {\it column} integral vectors of size $r$ and by $\ZZ^r_{\rm row}$ that of {\it row} vectors. The standard basis of $\ZZ^r_{\rm row}$ is denoted as $(e^1, \ldots, e^r)$.
For an integral matrix $M$,  we write 
$M\geq 0$ if all of its entries are non-negative. 

For $r$ independent variables ${\bf x} = (x_1, \ldots, x_r)$, $\ZZ_{\geq 0}[{\bf x}^{\pm 1}]$ denotes the set of all Laurent 
polynomials in ${\bf x}$ with non-negative integer coefficients.
If $R$ is commutative ring and $R^\times$ its group of units, for 
${\bf a} = (a_1, \ldots, a_n) \in (R^\times)^n$ we write
${\bf a}^L = a_1^{l_1} \cdots a_n^{l_n} \in R^\times$ for 
$L = (l_1, \ldots, l_n)^T \in \ZZ^n$, and if
$L$ is an $n \times m$ integral matrix with columns $L_1, \ldots, L_m \in \ZZ^n$, we write
${\bf a}^L = ({\bf a}^{L_1}, \ldots, {\bf a}^{L_m}) \in (R^\times)^m$.

For two sets $T$ and $S$, 
we will sometimes use $\{s_t\}_{t \in T}$ to denote (not a subset of $S$ but) a map $T \to S, t \mapsto s_t$, and
we call it {\it an assignment}. For
an $r$-tuple of elements ${\bf s} = (s_1, \ldots, s_r)$ in $S$, 
we also write $s \in {\bf s}$ to mean that $s = s_i$ for some $i \in [1, r]$.

\section{Positive spaces and tropical points from cluster algebras}\label{s:dual-pairing}

\subsection{Tropical points and global functions on positive spaces}
Let $\XX = (\FF(\XX), \,\{{\bf x}_t\}_{t \in \TT})$ be an $r$-dimensional positive space as in \deref{:pos-space}, and recall that 
\[
\XX(\Zmax) = {\rm Hom}_{\rm sf}(\FF_{>0}(\XX), \, \Zmax).
\]
For
$\rho \in \XX(\Zmax)$ and  $t \in \TT$, the vector 
$\val_\rho({\bf x}_t) \in \ZZ^r_{\rm row}$ is called the {\it coordinate vector} of $\rho$ at $t \in \TT$.
 We also say that an assignment $\{\rho_t \in \ZZ^r_{\rm row}\}_{t \in \TT}$ is a tropical point of $\XX$ if 
there exists $\rho \in \XX(\Zmax)$ such that
$\rho_t = \val_\rho({\bf x}_t)$ for every $t \in \TT$, and in this case we write 
$\rho = \{\rho_t \in \ZZ^r_{\rm row}\}_{t \in \TT}\in \XX(\Zmax)$. 
We also define
\[
\calU(\XX) = \bigcap_{t \in \TT} \ZZ[{\bf x}_t^{\pm 1}] \subset \FF(\XX) \hs \mbox{and} \hs
\calU^+(\XX) = \bigcap_{t \in \TT} \ZZ_{\geq 0}[{\bf x}_t^{\pm 1}] \subset \calU(\XX).
\]
 An element of $\FF(\XX)$ that lies in 
$\calU(\XX)$ (resp. in $\calU^+(\XX)$) is also said to be {\it universally Laurent} (resp. 
{\it universally positive Laurent}). Elements in $\calU(\XX)$ are also called {\it global functions} on $\XX$.

\bre{:XX-XX-prime}
{\rm
While any two positive spaces $\XX$ and $\XX^\prime$ with the same dimension have 
(non-natural) isomorphic ambient semi-fields and thus one can identify $\XX(\Zmax)$ with $\XX^\prime(\Zmax)$, the algebras
$\calU(\XX)$ and $\calU(\XX')$ depend on the coordinate charts of $\XX$ and $\XX'$ and are in general not isomorphic. In particular, 
given a positive space $\XX = (\FF, \{{\bf x}_t\}_{t \in \TT})$,
any non-empty subset $\TT'$ of $\TT$ gives rise to a positive space
$\XX' = (\FF, \{{\bf x}_t\}_{t \in \TT'})$ with $\calU(\XX) \subset \calU(\XX')$.
\hfill $\diamond$
}
\ere

\subsection{Matrix patterns}\label{ss:matrix-pattern} Recall again that an $r \times r$ integral skew-symmetrizable matrix is called 
a {\it mutation matrix of size $r$}.
Given a mutation matrix $B  = (b_{i, j})$ of size $r$ and  $k \in [1,r]$, 
the {\it mutation} of $B$ in direction $k$ is the $r \times r$  matrix $B' = (b_{i,j}')$, where
\begin{equation}\label{eq:matrix-mut}
		b_{i,j}' = \begin{cases}
			-b_{i,j}, & \text{ if } i = k \text{ or } j = k, \\
			b_{i,j} + [b_{i,k}]_+[b_{k,j}]_+ - [-b_{i,k}]_+[-b_{k,j}]_+, &\text{ otherwise},
		\end{cases}
	\end{equation}
and we write $B' =\mu_k(B)$. Then $\mu_k(B)$ is again a mutation matrix and  $\mu_k^2(B)= B$. 

Throughout the paper, $\TT_r$ denotes an
 $r$-regular tree with the $r$ edges from each vertex labeled bijectively by the set $[1, r]$,
 and we write $t \in \TT_r$ for a vertex $t$ of $\TT_r$. 
 If $t, t' \in \TT_r$ are
joined by an edge labeled by $k \in [1, r]$, we write \tkt.  We also fix $t_0 \in \TT_r$ and call it 
 {\it  the root} of $\TT_r$.
 
\bde{:b-pattern}
A {\it matrix pattern of rank $r$} is an assignment $\{B_t\}_{t \in \TT_r}$ of mutation matrix $B_t$ of size $r$ to each
$t \in \TT_r$, such that  $\mu_k(B_t) = B_{t'}$ whenever \tkt. 
Note that   
$\{B_t\}_{t \in \TT_r}$ is determined by  $B_{t_0}$. \hfill $\diamond$
\ede


\subsection{Positive $\bfA$-spaces, tropical mutation rule,  and $d$-compatibility degrees with cluster variables}\label{ss:X-seed-mut}
Let $\FF$ be a field isomorphic to the field of rational functions over $\QQ$ in $r$ independent variables. 

\bde{:seed}
{\rm
A {\it labelled seed (with trivial coefficients)}, or simply a {\it seed}, in $\FF$ is
a pair $({\bf z}, B)$, where ${\bf z}$ is an ordered set of free generators of $\FF$ over $\QQ$ and $B$ is 
a mutation matrix of size $r$.
}
\ede

Given a seed $({\bf x}, B)$ in $\FF$ and $k \in [1,r]$, the {\it $\bfA$-mutation} of $({\bf x}, B)$ in  direction $k$ is the seed 
$({\bf x}', \;\mu_k(B))$,
where 
 $\mu_k(B)$ is given 
 in \eqref{eq:matrix-mut} and  ${\bf x}' = (x_1', \dots, x_r')$ is given by

 \begin{equation}\label{eq:X-seed-mutation-cluster}
 x_i' = \begin{cases}
  x_k^{-1}(\prod_{j=1}^r x_j^{[b_{j,k}]_+} +\prod_{j=1}^r x_j^{[-b_{j,k}]_+}) ,& i = k; \\
  x_i,& \text{otherwise}.
  \end{cases}
\end{equation}

\bde{:cluster-pattern}
1) An {\it $\bfA$-pattern} in $\FF$ is an assignment $\{({\bf x}_t, B_t)\}_{t \in \TT_r}$ 
of a seed $({\bf x}_t, B_t)$ in $\FF$ to each $t \in \TT_r$ such that 
$({\bf x}_{t'}, B_{t'})$ is the $\bfA$-mutation of $({\bf x}_t, B_t)$ in direction $k$ whenever\tkt.

2) The {\it positive space associated to an $\bfA$-pattern $\{({\bf x}_t, B_t)\}_{t \in \TT_r}$} is the pair
$\AA =(\FF, \,\{{\bf x}_t\}_{t \in \TT_r})$.
We call each ${\bf x}_t$ an {\it $\bfA$-cluster}, or simply a {\it cluster}, of $\AA$,
each $x \in {\bf x}_t$ an {\it $\bfA$-variable}, or a {\it cluster variable}, of $\AA$, and $\{B_t\}_{t \in \TT_r}$ the {\it matrix pattern} of $\AA$.

3) We will also refer to positive spaces associated to $\bfA$-patterns as {\it positive $\bfA$-spaces}.
\hfill $\diamond$
\ede

\bnota{:AB}
For a mutation matrix $B$, 
we denote 
the positive space associated to any $\bfA$-pattern $\{({\bf x}_t, B_t)\}_{t \in \TT_r}$ with $B_{t_0} = B$ 
by $\AA_B = (\FF(\AA_B), \{{\bf x}_t\}_{t \in \TT_r})$ 
and call it
{\it the 
positive $\bfA$-space associated to $B$}.
\hfill $\diamond$
\enota

For the positive space 
$\AA_B$ with clusters $\{{\bf x}_t = (x_{t; 1}, \ldots, x_{t; r})\}_{t \in \TT_r}$ and
matrix pattern
$\{B_t = (b^t_{i, j}) \}_{t \in \TT_r}$, 
the following is obtained by tropicalizing \eqref{eq:X-seed-mutation-cluster}. 

\ble{:construct-trop-point-A}
An assignment $\{\delta_t \!= (\delta_{t;1}, \dots, \delta_{t;r}) \!\in \!\ZZ^r_{\rm row}\}_{t\in\mathbb T_r}$ 
is a tropical point of $\AA_B$ iff
	\begin{equation}\label{eq:ZT-point-comb-A}
		\delta_{t';i} = \begin{cases}
			 -\delta_{t;k} + \max(\sum_{j=1}^r [b_{j,k}^t]_+\delta_{t;j}\,, \; \sum_{j=1}^r [-b_{j,k}^t]_+\delta_{t;j}),& i = k; \\
			 \delta_{t;i},& \text{otherwise},
		\end{cases}
	\end{equation}
for all $i \in [1, r]$ whenever\tkt in $\TT_r$. 
\ele

Let $\mathfrak{X}(\AA_B)$ denote the set of cluster variables of $\AA_B$. Two cluster variables $x$ and $x'$ are said to be
{\it compatible} if there exists a cluster of $\AA_B$ containing both $x$ and $x'$. 

\bld{:x-degree}
For each $x \in \mathfrak{X}(\AA_B)$, there exists a unique $d_x \in \AA_B(\ZZ^{\max})$, called the {\it $d$-tropical point
of $x$}, such that 

    1) ${\rm val}_{d_x}(x) = -1$; and

    2) ${\rm val}_{d_x}(x') = 0$ for every $x' \in \mathfrak{X}(\AA_B)$ that is compatible with $x$ and $x'\neq x$.
\eld

\begin{proof}
Recall that $(e^1, \ldots, e^r)$ is the standard basis of $\ZZ^r_{\rm row}$.
Choose any $t \in \TT_r$ and $i \in [1,r]$ such that $x = x_{t;i}$, and define $d_x \in \AA_B(\Zmax)$ by setting  
${\rm val}_{d_x}({\bf x}_t) = -e^i$.
Suppose that 
$x= x_{t'; i'}$ for some $t' \in \TT_r$ and $i' \in [1,r]$. 
By \cite[Theorem 10]{cao-li:enough-g-pair}, there exists a sequence
\begin{xy}(-6,1)*+{t=}="A0",(0,1)*+{t_1}="A1",(10,1)*+{t_2}="A2",(20,1)*+{t_3}="A3",(30,1)*+{\ldots}="A4",
(42,1)*+{t_{l-1}}="A5",(55,1)*+{t_{l}}="A6",
\ar@{-}^{k_1}"A1";"A2", \ar@{-}^{k_2}"A2";"A3",\ar@{-}"A3";"A4"
\ar@{-}"A4";"A5",\ar@{-}^{~~k_{l-1}}"A5";"A6"\end{xy}
in $\TT_r$ such that $k_j \neq i$ for every $j \in [1, l-1]$ and 
 ${\bf x}_{t_l} = {\bf x}_{t^\prime}$ 
up to a re-ordering of their variables. By \eqref{eq:ZT-point-comb-A}, one has $\val_{d_x}({\bf x}_{t_2}) = -e^i$, 
and by induction, 
$\val_{d_x}({\bf x}_{t_l}) = -e^i$. Thus ${\rm val}_{d_x}({\bf x}_{t'}) = -e^{i'}$.
This shows that $d_x$ is independent of the choice of $t$ and $i$ and has the two properties as described.
\end{proof}

Recall that elements in the upper cluster algebra $\calU(\AA_B) =\bigcap_{t \in \TT_r} \ZZ[{\bf x}^{\pm 1}_t]$ are also
called global functions on 
$\AA_B$, and we have set $\calU^+(\AA_B) = \bigcap_{t \in \TT_r} \ZZ_{\geq 0}[{\bf x}^{\pm 1}_t]$.

\bde{:x-degree}
{\rm 
For $x\in \mathfrak{X}(\AA_B)$
and $x' \in \FF_{>0}(\AA_B)$, define the {\it $d$-compatibility degree of $x'$ with  $x$} to be
\begin{equation}\label{eq:valuation-d-comp-A}
    ( x \, \|\,  x'  )_d \stackrel{\text{def}}{=} {\rm val}_{d_x}(x').
\end{equation}

}
\ede

Recall that for any $x' \in \calU^+(\AA_B)\backslash\{0\}$ and any $t \in \TT_r$, one can write
\[
x^\prime = \frac{P({\bf x}_t)}{{\bf x}_t^{{\bf d}}},
\]
where ${\bf d} \in \ZZ^r$ and $P({\bf x}_t) \in \ZZ_{\geq 0}[{\bf x}_t]$
is not divisible by $x_{t; i}$ for any $i \in [1, r]$. The vector ${\bf d}$ is called the {\it denominator vector} of $x^\prime$ 
at $t$. For any $x \in \fX(\AA_B)$ and $x^\prime \in \calU^+(\AA_B)\backslash\{0\}$, 
the $d$-compatibility degree $(x \, \|\,  x')_d$ is thus the $i^{\rm th}$-component of the denominator vector of $x'$ 
at $t$ for any $t \in \TT_r$ such that $x_{t; i} = x$
(see \cite{cao_li_2020b} and \cite[Remark 2.31]{ckq_2022}).
The following  facts on {$(\cdot \, \|\, \cdot)_d$}  will be used
$\S$\ref{ss:hammock}.

\ble{:d}\cite[Remark 2]{cao-li:enough-g-pair} Let $B$ be any mutation matrix. For any $x, x' \in \mathfrak{X}(\AA_B)$, 

 1) $(x \,\|\, x')_d \geq -1$, and $(x \,\|\, x')_d = -1$ if and only if $x = x'$;

 2) $(x \,\|\, x')_d = 0$ if and only if $x$ and $x'$  are compatible and $x \neq x'$.
\ele

\ble{:u-v}
Let $B$ be any mutation matrix. Suppose that $u$ and $v$ are two cluster monomials on  $\AA_B$ such that 
$(x \, \| \, u)_d \leq  (x \, \| \, v)_d$ for every $x \in \fX(\AA_B)$, then $u = v$.
\ele

\begin{proof}
Choose any $t \in \TT_r$ such that  $u = x_{t; j_1}^{m_1} \cdots x_{t; j_l}^{m_l}$, where $\{j_1, \ldots, j_l\} \subset [1, r]$
and  $m_{k} > 0$ for each $k \in [1, l]$.  
Suppose that $v$ is a cluster monomial in the cluster ${\bf x}_{t'}$ for some $t' \in \TT_r$. 
By \leref{:d}, for each $i \in [1, r]$ one has $(x_{t'; i} \, \| \, v)_d \leq 0$. Thus
\[
(x_{t'; i} \, \| \, u)_d \leq (x_{t'; i} \, \| \, v)_d \leq 0, \hs \forall \; i \in [1, r].
\]
Suppose that $i \in [1, r]$ is such that $x_{t'; i} \notin \{x_{t; j_1}, \ldots, x_{t; j_l}\}$. It then follows from 
\leref{:d} and 
\[
(x_{t'; i} \, \| \, u)_d = m_1 (x_{t';i} \, \| \, x_{t; j_1})_d + \cdots + m_l (x_{t';i} \, \| \, x_{t; j_l})_d \leq 0
\]
that $(x_{t';i} \, \| \, x_{t; j_k})_d = 0$ for every $k \in [1, l]$. Thus 
$S=\{x_{t'; 1}, \ldots, x_{t'; r}, x_{t; j_1}, \ldots, x_{t; j_l}\}$
is a compatible subset of  $\fX(\AA_B)$  in the sense
that its elements are pairwise compatible. By \cite[Theorem 13]{cao-li:enough-g-pair}, $S$ is a subset of some cluster of $\AA_B$.
We can thus assume without loss of generality that $t = t'$. For each $i \in [1, r]$, one then has
$(x_{t; i} \, \| \, u)_d \leq (x_{t; i} \, \| \, v)_d$, and, if $x_{t; i}^\prime$ is the unique new
variable in $\mu_i({\bf x}_t)$, then one also has
\[
-(x_{t; i} \, \| \, u)_d  = (x_{t; i}^\prime  \, \| \, u)_d \leq (x_{t; i}^\prime  \, \| \, v)_d = 
-(x_{t; i} \, \| \, v)_d.
\]
So $(x_{t; i} \, \| \, u)_d = (x_{t; i} \, \| \, v)_d$ for each $i\in[1,r]$.
It follows that $u = v$.
\end{proof}


\subsection{Positive $\bfY$-spaces, tropical mutation rule, and $d$-compatibility degrees with $\bfY$-variables}\label{ss:Y-seed-mut}
Let again $\FF$ be a field isomorphic to the field of rational functions over $\QQ$ in $r$ independent variables. 

Given a seed $({\bf y}, B)$ in $\FF$ and $k \in [1,r]$, the {\it $\bfY$-mutation} of $({\bf y}, B)$ in 
 direction $k$ is the seed $({\bf y}^\prime, \mu_k(B))$, where 
 $\mu_k(B)$ is given 
 in \eqref{eq:matrix-mut} and  ${\bf y}' = (y_1', \dots, y_r')$ is given by\footnote{If necessary and clear from the context, 
 we denote by $\mu_k(\Sigma)$ both
 the $\bfA$-mutation and the $\bfY$-mutation of a seed $\Sigma$.}
	\begin{equation}\label{eq:Y-mut}
		y_i' = \begin{cases}
			y_k^{-1}, &  i = k; \\
			y_i y_k^{[b_{k,i}]_+} (1+y_k)^{-b_{k,i}}, &\text{otherwise}.
		\end{cases}
	\end{equation}

\bde{:y-pattern}
1) A {\it $\bfY$-pattern} in $\FF$ is an assignment $\{({\bf y}_t, B_t)\}_{t \in \TT_r}$ 
of a seed $({\bf y}_t, B_t)$ in $\FF$ to each $t \in \TT_r$ such that 
$({\bf y}_{t'}, B_{t'})$ is the $\bfY$-mutation of $({\bf y}_t, B_t)$ in direction $k$ whenever\tkt.

2) The {\it positive space associated to a $\bfY$-pattern} $\{({\bf y}_t, B_t)\}_{t \in \TT_r}$ is the pair $\YY = ( \FF, \,\{{\bf y}_t\}_{t \in \TT_r})$. We call each ${\bf y}_t$ a {\it \bfY-cluster} of $\YY$,
each $y \in {\bf y}_t$ a {\it $\bfY$-variable of $\YY$}, and $\{B_t\}_{t \in \TT_r}$ the {\it matrix pattern} of $\YY$.

3) We also refer to positive spaces associated to $\bfY$-patterns as {\it positive $\bfY$-spaces}.
\hfill $\diamond$
\ede

\bnota{:YB}
For a mutation matrix $B$, denote by 
$\YY_B = (\FF(\YY_B), \,\{{\bf y}_t\}_{t \in \TT_r})$
the positive space associated to any $\bfY$-pattern $\{({\bf y}_t, B_t)\}_{t \in \TT_r}$  
with  $B_{t_0} = B$, and call 
$\YY_B$ {\it the} positive $\bfY$-space associated to $B$.
\hfill $\diamond$
\enota

For the positive space 
$\YY_B$ with $\bfY$-clusters $\{{\bf y}_t = (y_{t; 1}, \ldots, y_{t; r})\}_{t \in \TT_r}$ and
matrix pattern
$\{B_t = (b^t_{i, j}) \}_{t \in \TT_r}$, 
the following is obtained by tropicalizing \eqref{eq:Y-mut}.

\ble{:construct-trop-point}
An assignment $\{\rho_t = (\rho_{t;1}, \dots, \rho_{t;r}) \in \ZZ_{\rm row}^r\}_{t\in\mathbb T_r}$ is a tropical point of $\YY_B$ iff
	\begin{equation}\label{eq:ZT-point-comb}
		\rho_{t';i} = \begin{cases}
			-\rho_{t;k}, & i = k; \\
			\rho_{t;i} + [b_{k,i}^t]_+\, \rho_{t;k} - b^t_{k,i} \, [\rho_{t;k}]_+, & \text{otherwise},
		\end{cases}
	\end{equation}
for all $i \in [1, r]$ whenever\tkt. 
\ele

Denote by $\mathfrak{X}(\YY_B)$ the set of $\bfY$-variables of $\YY_B$. 
Two elements $y, y' \in \mathfrak{X}(\YY_B)$ are said to be
{\it compatible} if there exists a $\bfY$-cluster of $\YY_B$ containing both $y$ and $y'$. 

\bld{:y-degree}
For every $y \in \mathfrak{X}(\YY_B)$, there exists a unique $d_y \in \YY_B(\ZZ^{\max})$, called the 
{\it $d$-tropical point of $y$}, such that 

1) ${\rm val}_{d_y}(y) = -1$; and 

2) ${\rm val}_{d_y}(y') = 0$ for every  $y' \in \mathfrak{X}(\YY_B)$ that is compatible with $y$ and $y'\neq y$.
\eld

\begin{proof}
Choose any $t \in \TT_r$ and $i \in [1,r]$ such that $y = y_{t;i}$, and define $d_y \in \YY_B(\Zmax)$ by setting  
${\rm val}_{d_y}({\bf y}_t) = -e^i$,
where again $(e^1, \ldots, e^r)$ is the standard basis of $\ZZ^r_{\rm row}$. 
Suppose that 
$y= y_{t'; i'}$ for some $t' \in \TT_r$ and $i' \in [1,r]$. 
By \cite[Theorem 7.10]{ckq_2022}, there exists a sequence
\[\begin{xy}(-6,1)*+{t=}="A0",(0,1)*+{t_1}="A1",(10,1)*+{t_2}="A2",(20,1)*+{t_3}="A3",(30,1)*+{\ldots}="A4",
(42,1)*+{t_{l-1}}="A5",(55,1)*+{t_{l}}="A6",
\ar@{-}^{k_1}"A1";"A2", \ar@{-}^{k_2}"A2";"A3",\ar@{-}"A3";"A4",
\ar@{-}"A4";"A5",\ar@{-}^{~~k_{l-1}}"A5";"A6"\end{xy}
\]
in $\TT_r$ such that $y \in {\bf y}_{t_j}$ for every $j \in [1, l]$ and ${\bf y}_{t_l} = {\bf y}_{t^\prime}$ 
up to a re-ordering of their variables. By the mutation formula in \eqref{eq:Y-mut} and by the algebraic independence of 
the variables in any $\bfY$-cluster, one sees that $k_j \neq i$ and $b^{t_j}_{k_j, i}= 0$ for every $j \in [1, l-1]$,
and
$y = y_{t_j; i}$ for every $j \in [1, l]$. By \eqref{eq:ZT-point-comb}, one has $\val_{d_y}({\bf y}_{t_2}) = -e^i$, 
and by induction, 
$\val_{d_y}({\bf y}_{t_l}) = -e^i$. Thus ${\rm val}_{d_y}({\bf y}_{t'}) = -e^{i'}$.
This shows that $d_y$ is well-defined and has the two required properties.
\end{proof}

Recall that elements in $\calU(\YY_B) = \bigcap_{t \in \TT_r} \ZZ[{\bf y}^{\pm 1}_t]$ are called global functions on $\YY_B$, and 
we have set $\calU^+(\YY_B) = \bigcap_{t \in \TT_r} \ZZ_{\geq 0}[{\bf y}^{\pm 1}_t]$.

\bde{:y-degree}
{\rm
For $y\in \mathfrak{X}(\YY_B)$ 
and $y' \in \FF_{>0}(\YY_B)$, define the {\it d-compatibility degree of $y'$ with $y$} to be
\[
    ( y \, \|\,  y' )_d \stackrel{\text{def}}{=} {\rm val}_{d_y}(y').
\]
For $y' \in \calU^+(\YY_B)\backslash\{0\}$, 
the integer $( y \, \|\,  y' )_d$ is the $y$-exponent of the
denominator of $y'$ as a Laurent polynomial of ${\bf y}_t$ for any $t \in \TT_r$ such that ${\bf y}_t$ contains $y$.
}
\ede


\bre{:not-in}
{\rm
Note that $\mathfrak{X}(\YY_B) \subset \FF_{>0}(\YY_B)$ but, unlike  the case of positive $\bfA$-spaces, 
$\mathfrak{X}(\YY_B)$ is {\it not} necessarily a subset of $\calU^+(\YY_B)$.
Moreover, for $y, y' \in \mathfrak{X}(\YY_B)$ and $y \neq y'$, if $y$ and $y'$ are compatible, then 
$(y\, \|\, y')_d = (y'\, \|\, y)_d = 0$ by \ldref{:y-degree}, but the converse is not necessarily true. 
Indeed, for $B = \left(\begin{array}{cc} 0 & -1\\ 1 & 0\end{array}\right)$ and ${\bf y}_{t_0} = (y_1, y_2)$, 
we have
\begin{align*}
(y_1, \, y_2) &\stackrel{\mu_1}{\longrightarrow} \left(\frac{1}{y_1}, \, y_2(1+y_1)\right) 
\stackrel{\mu_2}{\longrightarrow} \left(\frac{1+y_2+y_1y_2}{y_1}, \; \frac{1}{y_2(1+y_1)}\right) \\
&\stackrel{\mu_1}{\longrightarrow} \left(\frac{y_1}{1+y_2+y_1y_2}, \; \frac{1+y_2}{y_1y_2}\right) 
\stackrel{\mu_2}{\longrightarrow} \left(\frac{1}{y_2}, \; \frac{y_1y_2}{1+y_2}\right)
\stackrel{\mu_1}{\longrightarrow} (y_2, \, y_1).
\end{align*}
Thus $\mathfrak{X}(\YY_B)$ has $10$ elements, while we have 
\[
\mathfrak{X}(\YY_B) \cap \calU^+(\YY_B) = \left\{ y_1, \; y_2(1+y_1), \; \frac{1+y_2+y_1y_2}{y_1}, \; \frac{1+y_2}{y_1y_2}, \; 
\frac{1}{y_2}\right\}.
\]
Let $z = y_2(1+y_1)$. One also has $(y_1\, \|\, z)_d = (z\, \|\, y_1)_d = 0$, but $y_1$ and $z$ are not compatible.  
\hfill $\diamond$
}
\ere

\subsection{The ensemble map}\label{ss:p-map}
Let again $B$ be any mutation matrix of size $r$. Recall from \deref{:ensemble-pair} that the pair $(\AA_B, \YY_B)$ 
of positive spaces is called a Fock-Goncharov ensemble. 
Following \cite{FG:ensembles}, the map 
\[
\pB: \;\AA_B(\ZZ^{\max}) \longrightarrow \YY_B(\ZZ^{\max})
\]
induced by the semi-field homomorphism 
\[
\pB^*: \;\;\FF_{>0}(\YY_B) \longrightarrow  \FF_{>0}(\AA_B), \;\;
\pB^*({\bf y}_{t_0}) = {\bf x}_{t_0}^B,
\]
is called 
the {\it $p$-map}, or the {\it ensemble map}, of $(\AA_B, \YY_B)$.
By the well-known fact \cite[Proposition 3.9]{FZ:ClusterIV} that
$p_{\!{}_B}^*({\bf y}_t) = {\bf x}_t^{B_t}$ for every $t \in \TT_r$, 
where $\{B_t\}_{t \in \TT_r}$ is the matrix pattern with $B_{t_0} = B$. Thus
the map $\pB$ is linear in each coordinate chart of
$\AA_B(\Zmax)$, i.e.,
\begin{equation}\label{eq:p-map}
\pB:\; \AA_B(\ZZ^{\max}) \longrightarrow \YY_B(\ZZ^{\max}),\;\; 
\{\delta_t\in\mathbb Z_{\rm row}^r\}_{t\in\TT_r}\longmapsto 
\{\delta_tB_t\in\mathbb Z_{\rm row}^r\}_{t \in \TT_r}.
\end{equation}

\subsection{Positive spaces with principal coefficients}\label{ss:Aprin-trop}
Let $l$ be any non-negative integer, and consider any $(r+l) \times (r+l)$ skew-symmetrizable 
matrix
\begin{equation}\label{eq:breve-B}
\breve{B}  = \left(\begin{array}{cc} B & Q \\ P & R\end{array}\right),
\end{equation}
where  $B$ again has size $r \times r$. We then have  
the Fock-Goncharov ensemble $(\AA_{\breve{B}}, \YY_{\breve{B}})$, written as
\[
\AA_{\breve{B}} = (\FF(\AA_{\breve{B}}), \,
\{{{\bf x}}_t\}_{t \in \TT_{r+l}}) \quad \text{ and } \quad 
\YY_{\breve{B}} = (\FF(\YY_{\breve{B}}), \,\{{{\bf y}}_t\}_{t \in \TT_{r+l}}).
\]
Identify the root $t_0$ of $\TT_r$ with that of $\TT_{r+l}$ and 
 $\TT_r$ with the sub-tree of $\TT_{r+l}$ 
with edges labeled by indices from $[1, r] \subset [1, r+l]$. 
Let
$\{\breve{B}_t\}_{t \in \TT_{r+l}}$ be the matrix pattern with $\breve{B}_{t_0} = \breve{B}$, and let
\[
\breve{B}_t = \left(\begin{array}{cc} B_t & Q_t \\ P_t & R_t\end{array}\right) \hs \mbox{and} \hs
(B_t \bb P_t) = \left(\begin{array}{c} B_t \\ P_t\end{array}\right), 
\hs t \in \TT_{r+l}.
\]
One checks directly from the definitions that the assignment $\{({\bf x}_t, \, (B_t \bb P_t))\}_{t \in \TT_r}$
depends only on the $(r+l) \times r$ matrix $(B \bb P)$,  and the assignment
$\{({\bf y}_t, \, (B_t \; Q_t))\}_{t \in \TT_r}$
depends only on the $r \times (r+l)$ matrix $(B \; Q)$.
Correspondingly one has two positive spaces
\[
\AA_{(B \bb P)} = (\FF(\AA_{(B \bb P)}), \; \{{\bf x}_t\}_{t \in \TT_r}) \hs \mbox{and} \hs 
\YY_{(B \; Q)} = (\FF(\YY_{(B \; Q)}), \; \{{\bf y}_t\}_{t \in \TT_r}). 
\]
The
last $l$ elements in ${\bf x}_t$, being the same for all $t \in \TT_r$, are called the {\it frozen variables}
of $\AA_{(B \bb P)}$. 
Tropical points of $\AA_{(B \bb P)}$ are thus assignments
$\{{\delta}_t \in \ZZ^{r+l}_{\rm row}\}_{t \in \TT_r}$
that obey the mutation rule \eqref{eq:ZT-point-comb-A} with $B$ replaced by $\breve{B}$,
which depends only on $(B \bb P)$ when $t$ is
restricted to $\TT_r$.  Similarly, 
tropical points of $\YY_{(B \; Q)}$ are assignments 
$\{{\rho}_t \in \ZZ^{r+l}_{\rm row}\}_{t \in \TT_r}$
obeying the mutation rule \eqref{eq:ZT-point-comb} 
with $B$ replaced by $\breve{B}$,   which depends only on $(B \; Q)$ when restricted to $t \in \TT_r$.

\ble{:breve-B-p-map}
For any $\breve{B}$  as in \eqref{eq:breve-B},
the $p$-map for the ensemble $(\AA_{\breve{B}}, \YY_{\breve{B}})$ gives the map
\[
\pbreveB:\;\; \AA_{(B \bb P)}(\Zmax) \longrightarrow \YY_{(B \; Q)}(\Zmax), \;\;
\{\widetilde{\delta}_t \in \ZZ^{r+l}_{\rm row}\}_{t \in \TT_r} 
\mapsto \{\widetilde{\delta}_t \breve{B}_t\in \ZZ^{r+l}_{\rm row}\}_{t \in \TT_r},
\]
which is bijective if $\breve{B} \in GL(r+l, \ZZ)$. Furthermore, one has the well-defined maps
\begin{align*}
&\varepsilon_{\!{}_{(B\bb P)}}:\;\; \AA_B(\Zmax) \longrightarrow \AA_{(B \bb P)}(\Zmax), \;\; \{\delta_t \in \ZZ^r_{\rm row}\}_{t \in \TT_r} 
\mapsto \{(\delta_t, 0) \in \ZZ^{r+l}_{\rm row}\}_{t \in \TT_r},\\
&\omega_{\!{}_{(B \; Q)}}:\;\; \YY_{(B \; Q)}(\Zmax) \longrightarrow \YY_B(\Zmax), \;\; 
\{(\rho_t, \, \rho^\prime_t) \in \ZZ^{r+l}_{\rm row}\}_{t \in \TT_r} 
\mapsto \{\rho_t \in \ZZ^{r}_{\rm row}\}_{t \in \TT_r},
\end{align*}
with $\varepsilon_{\!{}_{(B\bb P)}}$ injective, 
$\omega_{\!{}_{(B \; Q)}}$ surjective, and $\omega_{\!{}_{(B \; Q)}} \circ \pbreveB \circ \varepsilon_{\!{}_{(B\bb P)}} = \pB$. 
\ele

\begin{proof}
The first statement follows from \eqref{eq:p-map} and the fact that if $\breve{B} \in GL(r+l, \ZZ)$, so does
$\breve{B}_t$ for all $t \in \TT_r$. The rest of the lemma follows
from the mutation rules for coordinate vectors of tropical points.
\end{proof}

Let $I_r$ be the identity matrix of size $r$. For any $r \times r$ mutation matrix $B$,  define
\[
\Aprin = \AA_{(B \bb I_r)} \hs \mbox{and} \hs \YY_B^{\rm prin} = \YY_{(B \; I_r)}
\]
and call them the {\it positive $\bfA$ and $\bfY$-spaces with principal coefficients} defined by $B$. 
Let
\begin{equation}\label{eq:o-B}
\overline{B} = \begin{pmatrix}
        B & -I_r \\
        I_r & 0
    \end{pmatrix},
\end{equation}
and let 
$\{\overline B_t\}_{t\in\mathbb T_{2r}}$ be the matrix pattern with $\overline B_{t_0}=\overline B$. For
$t \in \TT_{r}$, write
\begin{equation}\label{eq:wt-Bt}
\widetilde{B}_t = \begin{pmatrix}
        B_t \\
        C_t
\end{pmatrix} = (B_t \bb C_t) \hs \mbox{if} \hs \overline{B}_t = \begin{pmatrix}
        B_t & Q_t \\
        C_t & R_t
\end{pmatrix}.
\end{equation}
The following fact will be used in $\S$\ref{ss:admissible-YB}.

\ble{:YB-to-Aprin-trop}
With the notation as above, one has the well-defined injective map
\begin{equation}\label{eq:beta-B}
\beta_{B^T}:\;\; 
\AA_{B^T}(\Zmax) \longrightarrow \YY^{\rm prin}_{B^T}(\Zmax), \;\{\delta_t^\sv \in \ZZ^r_{\rm row}\}_{t \in \TT_r}
\longmapsto \{(\delta^\sv_t B_t^T, \, \delta_t^\sv C_t^T) \in \ZZ^{2r}_{\rm row}\}_{t \in \TT_r}.
\end{equation}
\ele

\begin{proof}
In the notation of  \leref{:breve-B-p-map}, one has 
$\beta_{B^T} =\pbarBT \circ \varepsilon_{\!{}_{(B^T \bb -I_r)}}$, and since $\overline{B} \in GL(2r, \ZZ)$, 
$\beta_{B^T}$ is injective by \leref{:breve-B-p-map}.
\end{proof}

\bre{:YB-to-Aprin-trop}
{\rm
The embedding $\beta_{B^T}$ in \eqref{eq:beta-B} appears in \cite[Page 554]{GHKK}. In  
\cite[Definition 7.12]{GHKK}, $\AA_{B^T}(\Zmax)$ is identified with its image under  $\beta_{B^T}$ as a {\it linear subspace}
of $\Yprin(\Zmax)$. Indeed, if $G_t = G_t^{B, t_0}$ for $t \in \TT_r$ in the notation of
\cite[\S 6]{FZ:ClusterIV}, then $G_tB_t = BC_t$ for every $t \in \TT_r$ by \cite[(6.14)]{FZ:ClusterIV}.
If we regard $\beta_{B^T}$ in \eqref{eq:beta-B} as a map from the set $R^\sv$ of all assignments 
$\delta^\sv =\{\delta_t^\sv \in \ZZ^r_{\rm row}\}_{t \in \TT_r}$ to the set $\wt{R}$ of all assignments
$\wt{\rho} =\{(\xi_t, \eta_t) \in \ZZ^r_{\rm row} \times \ZZ^r_{\rm row}\}_{t \in \TT_r}$, then 
$\wt{\rho} \in \wt{R}$ belongs to the image of $\beta_{B^T}$ if and only if 
\[
G_t \xi_t^T - B \eta_t^T = 0, \hs \;\; \forall \; t \in \TT_r,
\]
and in that case $\wt{\rho} = \beta_{B^T}(\delta^\sv)$ with $\delta^\sv = \{\eta_t(C_t^{-1})^T\}_{t \in \TT_r} \in R^\sv$.
Moreover, for any $\delta^\sv \in R^\sv$, one has 
$\delta^\sv \in \AA_{B^T}(\Zmax)$ if and only if $\beta_{B^T}(\delta^\sv) \in \Yprin(\Zmax)$. 
\hfill $\diamond$
}
\ere

\section{Global monomials and admissible functions}\label{s:admissible-elements}
For a mutation matrix $B$, recall from \deref{:ensemble-pair} that $(\AA_{B^T}, \YY_{B})$ is called a Fock-Goncharov dual pair 
(of positive spaces).
In this section, we introduce admissible functions on $\ABT$, $\YB$ and $\Aprin$. 
We also recall some facts about global monomials on $\ABT$, $\YB$ and $\Aprin$ and their $g$-vectors, which are either well-known or implicit in \cite{GHKK} but rephrased here
in the language of admissible functions.
In particular, 
we show that global monomials are admissible, and, for the convenience of the reader,  we give direct proofs that global monomials 
are in bijection with their $g$-vectors. 

Global monomials on $\YB$ that are also ${\bfY}$-variables are called {\it global variables
of $\YB$} (see \deref{:global-Y-variables}). As a new result,
we show in $\S$\ref{ss:d-g} that there is a well-defined injective map from the set of 
global variables of $\YB$ to the set of cluster variables of $\ABT$.

\subsection{Global monomials}\label{ss:global-mono}
Let $B$ be any $r \times r$ mutation matrix, and write
\begin{equation}\label{eq:AveeB-YB}
	\AA_{B^T} = (\FF(\AA_{B^T}), \{{\bf x}^\sv_t = (x^\sv_{t; 1}, \ldots, x^\sv_{t; r})\}_{t \in \TT_r}) \quad \text{ and } 
 \quad \YY_B = (\FF(\YY_B), \{{\bf y}_t=(y_{t; 1}, \ldots, y_{t; r})\}_{t \in \TT_r}).
\end{equation}
Elements in $\FF(\ABT)$ of the form $({\bf x}_t^\sv)^{\bf m}$ (resp. in $\FF(\YY_B)$ of the form ${\bf y}_t^{\bf m}$) 
for some $t \in \TT_r$ and ${\bf m} \in \ZZ^r$ are called {\it local Laurent monomials} on $\ABT$ (resp. on $\YY_B$).

\bde{:global-monomial} \cite[Definition 0.1]{GHKK} By a {\it global monomial on $\ABT$ (resp. on $\YY_B$)} we mean 
a local Laurent monomial on $\ABT$ (resp. on $\YY_B$) that lies in  $\calU(\ABT)$ (resp. in $\calU(\YY_B)$). 
\ede

Recall that  a {\em cluster monomial} on $\ABT$ is an element of the form $({\bf x}_t^\sv)^{\bf m}$ for some $t \in \TT_r$ and
${\bf m} \in \ZZ_{\geq 0}^r$. 
We prove the following well-known fact (see \cite[Definition 0.1]{GHKK}) for the convenience of the reader.

\ble{:A-global-mono}
For any local Laurent monomial $x^{\sv}$ on $\ABT$, the following are equivalent:

    1) $x^\sv$ is a cluster monomial on $\ABT$;

    2) $x^\sv$ is a global monomial on $\ABT$;
    
    3) for any $t \in \TT_r$ and ${\bf m} \in \ZZ^r$ such that $x^\sv = ({\bf x}_{t}^\sv)^{{\bf m}}$, one has ${\bf m} \geq 0$.
\ele

\begin{proof}
By the Laurent phenomenon, 1) implies 2). By definition, 3) implies 1). 

It remains to show that 2) implies 3). 
Let $x^\sv = ({\bf x}_t^\sv)^{\bf m}$, where $t \in \TT_r$ and ${\bf m} = (m_1, \ldots, m_k)^T \in \ZZ^r$, and assume that $x^\sv$ is a
global monomial. If $m_k < 0$ for some $k \in [1, r]$, then with  $t' \in \TT_r$ such that \tkt,
the mutation relation \eqref{eq:X-seed-mutation-cluster} gives
    \[
        x^\sv = M \, \left(\frac{P_k^+ + P_k^-}{x_{t';k}^\sv}\right)^{m_k}
    \]
for some $M$, $P_k^+$ and $P_k^-$ that are monomials in the variables in ${\bf x}_t^\sv \cap {\bf x}_{t'}^\sv$. So $x^\sv$
can not be 
in $\ZZ[({\bf x}_{t'}^\sv)^{\pm 1}]$, contradicting the assumption that $x^\sv$ is a global monomial.
Thus ${\bf m} \geq 0$.
\end{proof}

Let $\{B_t\}_{t \in \TT_r}$ be the matrix pattern with $B_{t_0} = B$.

\ble{:Y-global-mono}
For any local Laurent monomial $y$ on $\YY_B$, the following are equivalent:

    1) $y$ is a global monomial on $\YY_B$;

    2) for any $t \in \TT_r$ and ${\bf m} \in \ZZ^r$ such that $y = {\bf y}_{t}^{{\bf m}}$, one has $B_{t} {\bf m} \geq 0$;

    3) there exist $t \in \TT_r$ and ${\bf m} \in \ZZ^r$ such that $y = {\bf y}_t^{\bf m}$ and $B_{t} {\bf m} \geq 0$.
\ele

\begin{proof}
We first prove that 1) implies  2). Assume thus that $y$ is a global monomial on $\YY_B$, and let
$t \in \TT_r$ and ${\bf m} \in \ZZ^r$ be such that $y={\bf y}_{t}^{\bf m}$. Write 
$B_t{\bf m} =(v_1, \ldots, v_r)^T$ and suppose that there exists  $k \in [1,r]$ such that $v_k <0$.
Let $t' \in \TT_r$ be such that \tkt. 
Using \eqref{eq:Y-mut}, one sees that there exists a Laurent monomial $M$ in ${\bf y}_{t'}$ such that $y=M (1+y_{t';k})^{v_k}$,
contradicting the fact that $y \in \calU(\YY_B)$. Thus $B_t{\bf m} \geq 0$. This proves that 1) implies  2).
It is clear that 2) implies 3). 

To prove that 3) implies 1),
let $t \in \TT_r$ and ${\bf m}$ be as in 3), and let $w \in \TT_r$ be arbitrary. 
By the Separation Formulas \cite[Proposition 3.13]{FZ:ClusterIV}, 
there exist an $r \times r$ integral matrix $C^w_t$ and ${\bf F}^w_t = (F^w_{t; 1}({\bf y}_w), \ldots, F^w_{t; r}({\bf y}_w))$,
where $F^w_{t; j} \in \ZZ[{\bf y}_w]$ for each $j \in [1, r]$, such that 
${\bf y}_t = {\bf y}_w^{C^w_t} \left({\bf F}^w_t({\bf y}_w)\right)^{B_t}$ (recall notation in $\S$\ref{ss:nota-intro}).
Thus $y = {\bf y}_t^{{\bf m}} = {\bf y}_w^{C^w_t{\bf m}} \left({\bf F}^w_t({\bf y}_w)\right)^{B_t{\bf m}} \in 
\ZZ[{\bf y}_w^{\pm 1}]$ for any $w\in\mathbb T_r$. This implies  $y\in\calU(\YY_B)$ and thus $y$ is a global monomial on $\YY_B$.
\end{proof}

\subsection{Admissible functions on $\AA_{B^T}$}\label{ss:admissible-ABT}
Let again $B$ be any $r \times r$ mutation matrix and let $(\ABT, \YB)$ be the Fock-Goncharov dual pair 
as in \eqref{eq:AveeB-YB}.

\bde{:admissible-A}
1) Given $\rho = \{\rho_t \in \ZZ^r_{\rm row}\}_{t \in \TT_r} \in \YY_{B}(\Zmax)$, 
an element  $x^\sv \in \mathcal{U}(\AA_{B^T})$ is said to be {\it $\rho$-admissible} 
if for each $t \in \TT_r$ one has 
\begin{equation}\label{eq:a-g-x}
x^\sv = ({\bf x}_t^\sv)^{-\rho_t^T} \left(1 + \sum_{0 \neq {\bf u} \in U_t} \lambda_{{\bf u}} ({\bf x}_t^\sv)^{B_t^T{\bf u}}\right) 
\end{equation}
for some finite subset $U_t$ of $\ZZ_{\geq 0}^r$ and $\lambda_{\bf u} \in \ZZ_{\geq 0}$ for every ${\bf u} \in U_t$.  
Denote by 
$\calU^{\rm admi}_{\rho}(\AA_{B^T})$ the set of all $\rho$-admissible elements in $\calU(\AA_{B^T})$.

2) An element in $\calU(\AA_{B^T})$ is said to be {\it admissible}, and called an
{\it admissible function on $\AA_{B^T}$},
if it is $\rho$-admissible for some $\rho \in \YY_{B}(\Zmax)$.
The set of all admissible functions on $\ABT$ 
 is denoted by $\calU^{\rm admi}(\AA_{B^T})$. \hfill $\diamond$
\ede

For a given mutation matrix $B$, a natural question is whether or not $\calU^{\rm admi}_{\rho}(\AA_{B^T}) \neq \emptyset$ 
for every $\rho \in \YBZ$. We will show in \coref{:admissible-exist} that the answer is positive when $B$ is mutation equivalent to an acyclic matrix.  See also \reref{:large-complex}. On the other hand, 
the set $\calU^{\rm admi}_{\rho}(\AA_{B^T})$ may contain more than one element as in the following example.

\bex{:unique}
{\rm
For $B =\left(\begin{array}{cc} 0 & -2 \\ 2 & 0\end{array}\right)$ there are 
two well-known subsets of $\calU^{\rm admi}(\AA_{B^T})$: 

    1) the set of {\it greedy elements} $\{y_\rho : \rho \in \YY_{B}(\ZZ^{\max})\}$ defined in \cite{Lee-Li-Zelevinsky:greedy};

    2) the set of {\it generic variables} $\{z_\rho : \rho \in \YY_{B}(\ZZ^{\max})\}$ defined in \cite{Dupont-2011}.

\noindent
For $n \in \ZZ_{>0}$, let $\rho_n \in \YY_{B}(\ZZ^{\max})$ with coordinate vector $(n,-n) \in \ZZ_{\rm row}^2$
at $t_0$. By \cite[Example 3.11]{rank2theta}, 
$y_{\rho_1} = z_{\rho_1}$, $y_{\rho_2} =  (z_{\rho_1})^2-2 =z_{\rho_2} -2$, 
and both $y_{\rho_2}$ and $z_{\rho_2}$ are $\rho_2$-admissible. 
}
\eex

\bre{:admissible-1} 
Suppose that $B$ has full rank. Then for every $t \in \TT_r$, $B_t^T$ has full rank and
one has the partial order $\preceq_t$, called
the {\it dominance order} in \cite{Qin17},  on $\ZZ^{r}$,  
defined as 
\[
{\bf m} \preceq {\bf m}' \hs \mbox{if and only if} \hs {\bf m}' = {\bf m} + {{B}_t^T {\bf u}} \;\;\; \mbox{for some}\;\;\;
{\bf u} \in \ZZ_{\geq 0}^{r}.
\]
For any $x^\sv \in \calU^{\rm admi}(\AA_{B^T})$ and $t \in \TT_r$,
$-{\rho}_t^T$ is the unique minimal element among all the exponent vectors in the Laurent
expansion of $x^\sv$ in ${\bf x}^\sv_t$. In particular, when $B$ has full rank,
every $x^\sv \in \calU^{\rm admi}(\AA_{B^T})$ is $\rho$-admissible for a unique $\rho \in \YY_B(\Zmax)$, and admissible functions
on $\ABT$ are precisely all 
the {\it universally positive compatibly pointed elements} in $\calU(\AA_{B^T})$ in the sense of \cite[Definition 3.4.2]{qin-2022}.
If  $B$ does not have full rank, an arbitrary element $x^\sv \in \calU(\AA_{B^T})$ may have more than
one expression of the form in \eqref{eq:a-g-x} at some $t \in \TT_r$, so there is no well-defined 
vector $\rho_t \in \ZZ^r_{\rm row}$ associated to $x^\sv$ and $t \in \TT_r$. 
 \hfill $\diamond$
\ere

\bld{:globals-A-are-admissible}
Let $B$ be any $r \times r$ mutation matrix, and let $x^\sv$ be any cluster monomial on $\AA_{B^T}$. Let $t \in \TT_r$
and ${\bf m} \in \ZZ_{\geq 0}^r$ be such that $x^\sv  = ({\bf x}_t^\sv)^{\bf m}$, and let  
$\rho(x^\sv) = \{\rho(x^\sv)_w \in \ZZ^r_{\rm row}\}_{w \in \TT_r}$ be the 
unique tropical point of $\YY_B$ such that $\rho(x^\sv)_t = -{\bf m}^T$. Then 

1) $\rho(x^\sv) \in \YBZ$ depends only on $x^\sv$ and not on the choice of $t$;

2) $x^\sv$ is $\rho(x^\sv)$-admissible, and if $x^\sv$ is $\rho$-admissible for $\rho \in \YY_B(\Zmax)$, then $\rho = \rho(x^\sv)$.

\noindent
The tropical point $\rho(x^\sv) \in \YY_B(\Zmax)$ 
is called the {\it $g$-vector}\footnote{This terminology is slightly different from that in \cite[\S 6]{FZ:ClusterIV}, where it is the column vector $-\rho(x^\sv)_w^T$ which Fomin and Zelevinsky call {\it the $g$-vector} of $x^\sv$ (with respect to the fixed vertex $w\in\mathbb T_r$).} of $x^\sv$ in \cite[Definition 5.8]{GHKK}. 
\eld

\begin{proof} 
We prove 2) first. Let $t$, ${\bf m}$ and  
$\rho(x^\sv) \in \YBZ$ be as stated. For $w \in \TT_r$, let $G_t^w:=G_t^{B_w^T; w}$ be the $G$-matrix and
${\bf F}_t^w := (F_{t; 1}^{B_w^T; w}, \ldots, F_{t; r}^{B_w^T; w})$ 
be the $F$-polynomials in $r$ independent variables 
defined in \cite[\S 6]{FZ:ClusterIV}. Then  $G_t^t = I_r$ and 
${\bf x}_t^\sv = ({\bf x}_w^\sv)^{G_t^w} {\bf F}_t^w(({\bf x}_w^\sv)^{B_w^T})$. 
It follows that 
\[
x^\sv = ({\bf x}_t^\sv)^{\bf m} =  ({\bf x}_w^\sv)^{G_t^w{\bf m}} \left({\bf F}_t^w(({\bf x}_w^\sv)^{B_w^T})\right)^{{\bf m}}, \hs w \in \TT_r.
\]
Let $i \in [1, r]$ and  
$\rho_{(t; i)} =\{-(G_t^we_i)^T  \in \ZZ_{\rm row}^r \}_{w\in \TT_r}$, where
$e_i$  is the $i^{\rm th}$ standard basis vector of $\mathbb Z^r$.
By \cite[Proposition 4.2]{NZ12}, 
the {\it sign-coherence of $c$-vectors} \cite[Corollary 5.5]{GHKK} implies that 
$\rho_{(t; i)} \in \YY_{B}(\Zmax)$ and the $F$-polynomial $F_{t; i}^{B_w^T, w}$ has constant term $1$ for $i\in[1,r]$.
By the
{\it $G$-matrix sign coherence} \cite[Theorem 5.11]{GHKK} and the fact that  ${\bf m} \in \ZZ_{\geq 0}^r$, the assignment
$\left\{-(G_t^w {\bf m})^T  \in \ZZ_{\rm row}^r\right \}_{w\in \TT_r}$
is a tropical point of $\YY_B$, which must equal $\rho(x^\sv)$ as the two have the same coordinate vector 
$-{\bf m}^T$ at $t$. Moreover, for every $w \in \TT_r$, one knows by \cite[Corollary 0.4]{GHKK} that the $F$-polynomial
$F_{t; i}^{B_w^T, w}$ has non-negative integer coefficients for every $i \in [1, r]$, so 
$\left({\bf F}_t^w(({\bf x}_w^\sv)^{B_w^T})\right)^{{\bf m}}$ is a polynomial in 
$({\bf x}_w^\sv)^{B_w^T}$ with non-negative integer coefficients and constant term $1$. It follows that
$x^\sv$ is $\rho(x^\sv)$-admissible.

Suppose now that $\rho = \{\rho_w \in \ZZ^r_{\rm row}\}_{w \in \TT_r}$ is any point in $\YBZ$ 
such that $x^\sv$ is $\rho$-admissible. Then 
\[
({\bf x}_t^\sv)^{{\bf m}} = x^\sv =({\bf x}_t^\sv)^{-\rho_t^T} F_t(({\bf x}_t^\sv)^{B_t^T})
\]
for some $F_t \in \ZZ_{\geq 0}[y_1, \ldots, y_r]$
with constant term $1$.
By the linear independence of 
Laurent monomials in  ${\bf x}_{t}^\sv$, one must have $F_t = 1$ and $\rho_t = -{\bf m}^T$.  Thus $\rho = \rho(x^\sv)$, which concludes the proof of 2). 

1) is \cite[(2) of Lemma 7.10]{GHKK}. However, 1) also follows easily from 2). Indeed, 
assume that $x^\sv=({\bf x}_{t'}^\sv)^{{\bf m}'}$ for some 
$t'\in\mathbb T_r$ and ${\bf m}'\in \mathbb Z_{\geq 0}^r$. Let $\rho'$ be the unique tropical point 
of $\YB$ with coordinate vector $-({\bf m}^\prime)^T$ at $t'$. Then by 2), $x^\sv$ is $\rho'$-admissible, 
and $\rho' = \rho(x^\sv)$.
\end{proof}

Let $\Delta^+(\YY_B, \ZZ) \subset \YY_B(\Zmax)$ be the set\footnote{Our
$\Delta^+(\YY_B, \ZZ)$ is $\Delta^+_{\AA_{B^T}}(\ZZ)$ in the notation of \cite[Definition 7.9]{GHKK}.} of
$g$-vectors
of all cluster monomials on $\AA_{B^T}$, i.e., 
\begin{equation}\label{eq:Delta-plus-YBZ}
\Delta^+(\YY_B, \ZZ)=\{\rho = \{\rho_t\in \ZZ^r_{\rm row}\}_{t \in \TT_r}\in \YY_B(\Zmax): 
-\rho_t \geq 0\, \mbox{for some}\,  t \in \TT_r\}. 
\end{equation}
Denote by $\calU^{\rm mono}(\AA_{B^T})$ the set of all cluster monomials on  $\AA_{B^T}$.  By 
\ldref{:globals-A-are-admissible}, one has 
\[
\calU^{\rm mono}(\AA_{B^T}) \subset \calU^{\rm admi}(\ABT) \subset \calU(\ABT).
\]

\ble{:Delta-plus}
For any mutation matrix $B$, one has the well-defined bijection 
\begin{equation}\label{eq:mono-bijection-A}
\calU^{\rm mono}(\AA_{B^T}) \longrightarrow  \Delta^+(\YY_B, \ZZ),\;\; x^\sv \longmapsto \rho(x^\sv).
\end{equation}
\ele

\begin{proof} By \ldref{:globals-A-are-admissible}, 
the map in \eqref{eq:mono-bijection-A} is well-defined and it is surjective. 
The injectivity of the map in \eqref{eq:mono-bijection-A} seems well-known and 
is implied in \cite[(5) of Theorem 0.3]{GHKK} via the construction of theta functions, but 
we provide an elementary proof for the convenience of the reader. 

Suppose that $x^\sv, z^\sv \in \calU^{\rm mono}(\AA_{B^T})$ are such that 
$\rho(x^\sv) = \rho(z^\sv) = \{\rho_t \in \ZZ^r_{\rm row}\}_{t \in \TT_r}$, and let $t \in \TT_r$ be such that $-\rho_t \geq 0$ and
$x^\sv= ({\bf x}_t^\sv)^{-\rho_t^T}$. As $z^\sv$ is also $\rho$-admissible, one has
\[
z^\sv = ({\bf x}_t^\sv)^{-\rho_t^T} F(({\bf x}_t^\sv)^{B_t^T}) = x^\sv F(({\bf x}_t^\sv)^{B_t^T})
\]
for some $F \in \ZZ_{\geq 0}[y_1, \ldots, y_r]$ with constant term $1$. Let $f =F(({\bf x}_t^\sv)^{B_t^T})-1 \in \FF(\AA_{B^T})$,
so either $f = 0$ or is in $\FF_{>0}(\AA_{B^T})$. 
Similarly, there exists $t^\prime \in \TT_r$ and 
$F^\prime \in \ZZ_{\geq 0}[y_1, \ldots, y_r]$ with constant term $1$ such that 
$x^\sv = z^\sv F^\prime(({\bf x}_{t'}^\sv)^{B_{t'}^T})$. Let $f' = F^\prime(({\bf x}_{t'}^\sv)^{B_{t'}^T}) - 1$, which is again either $0$ or in $\FF_{>0}(\AA_{B^T})$. 
It follows that 
$1 = F(({\bf x}_t^\sv)^{B_t^T})F^\prime(({\bf x}_{t'}^\sv)^{B_{t'}^T})  = (1+f) (1+f^\prime)$,
so 
$f + f^\prime + ff^\prime = 0$ which can happen only if $f = f' = 0$. Thus $x^\sv = z^\sv$.  
Hence the map in \eqref{eq:mono-bijection-A} is injective.
\end{proof}

\bre{:g-fan}
{\rm
Consider the semi-field $\RR^{\rm max} = (\RR, \cdot, \oplus)$ with 
$a \cdot b = a+b$ and $a \oplus b = {\rm max}(a, b)$ for $a, b \in 
\RR$ and the set $\YY_B(\RR^{\rm max}) = {\rm Hom}_{\rm sf}(\FF_{>0}(\YY_B), \RR^{\rm max})$ 
of $\RR^{\rm max}$-tropical points
of $\YY_B$. A point in $\YY_B(\RR^{\rm max})$ is then an assignment $\rho = \{\rho_t \in \RR^r_{\rm row}\}_{t \in \TT_r}$ 
to each $t \in \TT_r$ a (row) vector of size $r$ with entries in $\RR$ which
 obeys the same
mutation rule \eqref{eq:ZT-point-comb} as points in $\YBZ$, so $\YBZ\subset \YY_B(\RR^{\rm max})$ as sets.
For $t \in \TT_r$, let $C_t^+(\RR) \subset \YY_B(\RR^{\rm max})$ 
be the set of all $\rho \in \YY_B(\RR^{\rm max})$ such that 
$-\rho_t \geq 0$, and call $C_t^+(\RR)$ the {\it Fock-Goncharov cluster chamber} associated to $t$
\cite[Definition 2.9]{GHKK}. Define 
\[
\Delta^+(\YY_B, \RR) = \bigcup_{t \in \TT_r} C_t^+(\RR) \subset \YY_B(\RR^{\rm max}). 
\]
By \cite[Theorem 2.13]{GHKK}, under the identification 
\begin{equation}\label{eq:identification}
\YY_B(\RR^{\rm max}) \cong \RR^r, \;\; \rho = \{\rho_t \in \RR^r_{\rm row}\}_{t \in \TT_r} \longmapsto -\rho_{t_0}^T,
\end{equation}
the collection $\Delta^+(\YY_B)$ of the subsets $C_t^+(\RR)$ over $t \in \TT_r$ form the maximal cones of
a simplicial fan in $\RR^r$, called the {\it $g$-vector fan of $\ABT$}. 
Note that under the identification in \eqref{eq:identification}, $\Delta^+(\YY_B, \ZZ)$ is exactly the set of all integral points in $\Delta^+(\YY_B, \RR)$. 
\hfill $\diamond$
}
\ere

\subsection{Admissible functions on $\Aprin$}\label{ss:admissible-Aprin}
For an arbitrary full rank extended mutation matrix $(B \bb P)$ as in $\S$\ref{ss:Aprin-trop}, one can define 
admissible functions on $\AA_{(B\bb P)}$ using tropical points of
$\YY_{(B^T \; P^T)}$, and we refer to \cite{cao:F-invariant} for this general setting. For our paper, we only 
need to consider admissible functions
on positive $\bfA$-spaces with principal coefficients.

Let $B$ be any $r \times r$ mutation matrix, and consider the Fock-Goncharov dual pair  
\[
(\Aprin = \AA_{(B \bb I_r)}, \;\; \Yprin = \YY_{(B^T\; I_r)}).
\]
For $t \in \TT_r$, let $\wt{B}_t = (B_t \bb C_t)$ be given in 
\eqref{eq:wt-Bt}, and denote the clusters of $\Aprin$ as $\{\wt{{\bf u}}_t = ({\bf u}_t, {\bf p})\}_{t \in \TT_r}$,
where ${\bf p} = (p_1, \ldots, p_r)$ are the frozen variables.
Set $\widehat{\bf y}_t= \widetilde{\bf u}_t^{\widetilde B_t} 
={\bf u}_t^{B_t}{\bf p}^{C_t}$ for $t \in \TT_r$.

\bde{:admi-prin}
1) Given $\widetilde{\rho} = \{\widetilde{\rho}_t \in \ZZ^{2r}_{\rm row}\}_{t \in \TT_r} \in \YY^{\rm prin}_{B^T}(\Zmax)$, an 
element $\widetilde{u} \in \calU(\AA^{\rm prin}_B)$ is said to be {\it $\widetilde{\rho}$-admissible} if for each
$t \in \TT_r$, one has
$\widetilde{u} = \widetilde{{\bf u}}_t^{-\widetilde{\rho}_t^T} F_t(\widehat{\bf y}_t)$
for some $F_t \in \ZZ_{\geq 0}[y_1, \ldots, y_r]$ with constant term $1$. The set of all $\widetilde{\rho}$-admissible elements in
$\calU(\Aprin)$ is denoted by $\calU^{\rm admi}_{\widetilde{\rho}}(\Aprin)$;

2) An element $\widetilde{u} \in \calU(\AA^{\rm prin}_B)$ that is $\widetilde{\rho}$-admissible for some 
$\widetilde{\rho} \in \YY_{B^T}^{\rm prin}(\Zmax)$ is called an {\it admissible function on $\Aprin$}.
The set of all admissible functions on $\Aprin$ is denoted by $\calU^{\rm admi}(\AA^{\rm prin}_B)$.
\ede

\bre{:prin-unique}
{\rm 
As $(B \bb I_r)$ has full rank, the same arguments in \reref{:admissible-1} show that any 
$\widetilde{u} \in \calU^{\rm admi}(\AA^{\rm prin}_B)$ is $\widetilde{\rho}$-admissible for a unique 
$\widetilde{\rho}  \in \YY_{B^T}^{\rm prin}(\Zmax)$. 
\hfill $\diamond$
}
\ere

An {\it extended cluster monomial} of $\Aprin$ is any $\wt{u} \in \FF_{>0}(\Aprin)$ such that 
$\wt{u} = {\bf u}_t^{{\bf m}} {\bf p}^{{\bf n}}$ for some $t \in \TT_r$, ${\bf m} \in \ZZ_{\geq 0}^r$, and ${\bf n} 
\in \ZZ^r$. Denote by $\calU^{\rm mono}(\Aprin)$ the set of all extended cluster monomials on  $\Aprin$. 
Let 
\begin{equation}\label{eq:Delat-pulus-prin}
\Delta^+(\Yprin, \ZZ) = \{\wt{\rho}  \in \Yprin(\Zmax):  \mbox{first}\; r\; \mbox{coordinates of}\; \wt{\rho}\; \mbox{non-positive
 at some}\; t \in \TT_r\}.
\end{equation}
The following statement and definition can again be found in \cite[Definition 5.10 and Lemma 7.10]{GHKK}.
We give a direct proof using \ldref{:globals-A-are-admissible}. 

\bld{:admi-prin}
Let $\wt{u} \in \calU^{\rm mono}(\Aprin)$. Let $t \in \TT_r$, ${\bf m} \in \ZZ_{\geq 0}^r$, and ${\bf n} 
\in \ZZ^r$ be such that $\wt{u} = {\bf u}_t^{{\bf m}} {\bf p}^{{\bf n}}$, let $\wt{\rho}(\wt{u})$ 
be the unique element in $\Delta^+(\Yprin, \ZZ)\subset \Yprin(\Zmax)$ with coordinate vector $-({\bf m}^T, {\bf n}^T)$ at $t$.  Then
the map
\begin{equation}\label{eq:rho-Aprin}
\calU^{\rm mono}(\Aprin) \longrightarrow \Delta^+(\Yprin, \ZZ), \;\; \wt{u} \longmapsto \wt{\rho}(\wt{u}),
\end{equation}
is a well-defined bijection and $\wt{u}$ is $\wt{\rho}(\wt{u})$-admissible. We call $\wt{\rho}(\wt{u}) \in \Yprin(\Zmax)$ 
the {\it extended $g$-vector} of $\wt{u} \in \calU^{\rm mono}(\Aprin)$.
\eld

\begin{proof}
Suppose that $\wt{u} = {\bf u}_t^{{\bf m}} {\bf p}^{{\bf n}}$ for some $t \in \TT_r$, ${\bf m} \in \ZZ_{\geq 0}^r$, and ${\bf n} 
\in \ZZ^r$. Then
$\wt{u}^\prime := {\bf u}_t^{{\bf m}}$ is a cluster monomial on  $\AA_{\overline{B}} = (\FF(\AA_{\overline{B}}), \{({\bf u}_w, {\bf p})\}_{w \in \TT_{2r}})$, where
$\overline{B}$ is given in \eqref{eq:o-B}. By \ldref{:globals-A-are-admissible}, $\wt{u}^\prime$ is
$\overline{\rho}$-admissible, where $\overline{\rho}$ is the tropical point in $\YY_{\overline{B}^T}(\Zmax)$ with coordinate vector  $(-{\bf m}^T,0)\in\mathbb Z_{\rm row}^{2r}$ at the vertex $t$.
 For any $w \in \TT_r \subset \TT_{2r}$, one has
\[
    \wt{u}^\prime ={\bf u}_t^{{\bf m}}= ({\bf u}_w, {\bf p})^{-\overline{\rho}_w^T} \, ({\bf F}_t^{w} (({\bf u}_w, {\bf p})^{\overline{B}_w}))^{{\bf m}'}, \quad {\bf m}' = \begin{pmatrix}
        {\bf m}\\ 0
    \end{pmatrix}\in \mathbb Z^{2r},
\]
where $\{\overline{B}_w\}_{w\in \TT_{2r}}$ is the matrix pattern of $\AA_{\overline{B}}$, and 
${\bf F}^w_t = (F^{\overline{B}_w;w}_{t;1}, \ldots , F^{\overline{B}_w;w}_{t;2r})$
are $F$-polynomials (see \cite[\S 5]{FZ:ClusterIV}). Since $w \in \TT_r \subset \TT_{2r}$ and by a recursive argument using \cite[(5.2), (5.3) and (5.9)]{FZ:ClusterIV}, one can show
that $F^{\overline{B}_w;w}_{t;i}$ depends only on $\widehat{{\bf y}}_w={\bf u}_w^{B_w}{\bf p}^{C_w}$ and that  $F^{\overline{B}_w;w}_{t;r+i} = 1$ for each $i \in [1,r]$. Moreover, it follows from the mutation rule of
the coordinate vectors of tropical points of $\Yprin$ that 
$\overline{\rho}_w -\wt{\rho}(\wt{u})_w = (0, {\bf n}^T)$ for every $w \in \TT_r$. Applying 
\ldref{:globals-A-are-admissible} to $\wt{u}^\prime$ and by  $\wt{u} = {\bf u}_t^{{\bf m}} {\bf p}^{{\bf n}}=\wt{u}^\prime{\bf p}^{{\bf n}}$, one has that
$\wt{u}$ is $\wt{\rho}(\wt{u})$-admissible. By the same arguments as in the proof of \ldref{:globals-A-are-admissible}, one can show that $\wt{\rho}(\wt{u})$ is independent of the choice of $t$ and thus it is well-defined. By the same arguments as in the proof of \leref{:Delta-plus}, one can show that the map
in \eqref{eq:rho-Aprin} is bijective.
\end{proof}

The following lemma will be used in \S \ref{ss:thmE}.

\ble{:gvec-mono-coeff-t0}
Let $\AA_B^{\rm prin}$ be the positive $\bfA$-space with principal coefficients defined by $B$ at  $t_0\in\mathbb T_r$, and let $\widetilde{u}$ be an extended cluster monomial of $\AA_B^{\rm prin}$ with extended 
    $g$-vector $\widetilde{\rho} = \{(\rho_t, \mu_t) \in \ZZ^r_{\rm row} \times \ZZ^r_{\rm row} : t \in \TT_r\}$. If 
    $t \in \TT_r$ is such that $\widetilde{u}$ is a Laurent monomial in $\widetilde{{\bf u}}_t$, then 
    $\mu_{t} = \mu_{t_0}$.
\ele
\begin{proof}
     Let ${\rm Trop}({\bf p})$ be the tropical semi-field generated by ${\bf p} = (p_1, \ldots, p_r)$, and let
\begin{equation}\label{eq:T}
T: \;\; \FF_{>0}(\AA_B^{\rm prin}) \longrightarrow {\rm Trop}({\bf p})
\end{equation}
be the unique  semi-field homomorphism which sends 
$p_i \in \FF_{>0}(\AA_{B}^{\rm prin})$ to $p_i \in {\rm Trop}({\bf p})$ for every $i \in [1, r]$ and sends every cluster 
variable of $\Aprin$ to $1 \in {\rm Trop}({\bf p})$ (see \cite[Lemma 2.2]{GHL:friezes}). Suppose that $t \in \TT_r$ is such that 
$\widetilde{u} = {\bf u}_t^{-\rho_t^T}
     {\bf p}^{-\mu_t^T}$. Then $T(\widetilde{u}) = {\bf p}^{-\mu_t^T}$.
On the other hand, since $\widetilde{u}$ is $\widetilde{\rho}$-admissible, one has
\begin{equation}\label{eq:mono-at-t-t0-exp}
        \widetilde{u} = {\bf u}_{t_0}^{-\rho_{t_0}^T} {\bf p}^{-\mu_{t_0}^T} F_{t_0}(\widehat{\bf y}_{t_0}),
     \end{equation}
for some polynomial $F_{t_0}$ with non-negative integer coefficients and constant term $1$. Since $T(\widehat{\bf y}_{t_0})={\bf p}$, 
\[
T(\widetilde{u}) = T\left({\bf u}_{t_0}^{-\rho_{t_0}^T} {\bf p}^{-\mu_{t_0}^T} F_{t_0}(\widehat{\bf y}_{t_0})\right) 
={\bf p}^{-\mu_{t_0}^T} T(F_{t_0}(\widehat{\bf y}_{t_0})) = {\bf p}^{-\mu_{t_0}^T}.
\]
It follows that $\mu_t = \mu_{t_0}$.
\end{proof}

\subsection{Admissible functions on $\YY_B$}\label{ss:admissible-YB}
We continue with the notation  in \eqref{eq:AveeB-YB}.

\bde{:admissible-Y}
1) Given $\delta^\sv = \{\delta_t^\sv \in \ZZ^r_{\rm row}\}_{t \in \TT_r} \in \AA_{B^T}(\Zmax)$, 
an element  $y \in \mathcal{U}(\YY_{B})$ is said to be {\it $\delta^\sv$-admissible} 
if for each $t \in \TT_r$ one has 
\begin{equation}\label{eq:a-g-y}
y = ({\bf y}_t)^{-(\delta_t^\sv)^T} \left(1 + \sum_{0 \neq {\bf v} \in V_t} \mu_{\bf v} {\bf y}_t^{\bf v}\right) \quad 
\end{equation}
for some finite subset $V_t$ of $\ZZ_{\geq 0}^r$ and $\mu_{\bf v} \in \ZZ_{\geq 0}$ for every ${\bf v} \in V_t$. 
 The set of all $\delta^\sv$-admissible elements of $\calU(\YY_B)$ is denoted by
$\calU^{\rm admi}_{\delta^\sv}(\YY_B)$.

2) An element in $\calU(\YY_B)$ that is $\delta^\sv$-admissible for some 
$\delta^\sv \in \AA_{B^T}(\Zmax)$ is called an {\it admissible function on $\YY_B$}. The set of all admissible functions on $\YB$ 
is denoted by $\calU^{\rm admi}(\YY_B)$. 
\hfill $\diamond$
\ede

\bre{:delta-unique}
{\rm
Simpler than the case of positive $\bfA$-spaces, every non-zero $y \in \calU^+(\YY_B)$ gives rise to a well-defined assignment 
$\delta^\sv = \{\delta_t^\sv \in \ZZ^r_{\rm row}\}_{t \in \TT_r}$ via its Laurent expansion
\[
y = {\bf y}_t^{-(\delta_t^\sv)^T} F_t({\bf y}_t)
\]
at every $t \in \TT_r$, where $F_t({\bf y}_t) \in \ZZ_{\geq 0}[{\bf y}_t]$ is not divisible by $y_{t; i}$ for any 
$i \in [1, r]$, and $(\delta_t^\sv)^T  \in \ZZ^r$ is the {\it denominator vector} of $y$ at $t$.
Then $y \in \calU^+(\YY_B)$ is admissible if and only if $\delta^\sv = \{\delta_t^\sv \in 
\ZZ^r_{\rm row}\}_{t \in \TT_r}$ is a tropical point of $\AA_{B^T}$ and $F_t$ has constant term $1$ 
for every $t \in \TT_r$. In particular, every $y \in \calU^{\rm admi}(\YY_B)$ is
$\delta^\sv$-admissible for a unique $\delta^\sv \in \AA_{B^T}(\Zmax)$.
\hfill $\diamond$
}
\ere

Following the idea used in \cite{GHKK}
to construct theta functions on $\YY_B$ through those on $\Aprin$ and in \cite{kimura-qin-wei} on bases of $\calU(\YY_B)$,
we now show that $\calU^{\rm admi}(\YY_B)$ can be embedded in 
$\calU^{\rm admi}(\Aprin)$.

Let the notation be as in $\S$\ref{ss:admissible-Aprin}. 
Since $\widetilde{B}$ has full rank, the variables in $\widehat{{\bf y}}_{t_0}$ are algebraically independent in
$\FF(\AA^{\rm prin}_B)$. One thus has the injective 
field homomorphism 
\[
\widetilde{p}_{\!{}_B}^*: \;\; \FF(\YY_B) \longrightarrow \FF(\AA_B^{\rm prin}), \;\; 
{\bf y}_{t_0} \longmapsto \widehat{{\bf y}}_{t_0}.
\]
By \cite[Proposition 3.9]{FZ:ClusterIV}, one has $\widetilde{p}_{\!{}_B}^*({\bf y}_t) = \widehat{{\bf y}}_{t}$
for every $t \in \TT_r$. Consequently, 
for any $y \in \calU(\YY_B)\backslash\{0\}$, if $y$ has denominator vector $(\delta_t^\sv)^T \in \ZZ^r$ at $t \in \TT_r$ 
and Laurent expression 
$y = {\bf y}_t^{-(\delta_t^\sv)^{T}} F_t({\bf y}_t)$, then 
\begin{equation}\label{eq:p-star}
\widetilde{p}_{\!{}_B}^*(y) = \wh{{\bf y}}_t^{-(\delta_t^\sv)^{T}} F_t(\widehat{\bf y}_t)
=\widetilde{\bf u}_t^{-\widetilde{\rho}_t^T} F_t(\widehat{\bf y}_t),
\hs \mbox{where}\hs\widetilde{\rho}_t = (\delta_t^\sv B_t^T, \; \delta_t^\sv C_t^T) \in \ZZ^{2r}_{\rm row}.
\end{equation}
In particular $\widetilde{p}_{\!{}_B}^*$ restricts to an injective ring homomorphism 
\begin{equation}\label{eq:p-YB-Aprin}
\widetilde{p}_{\!{}_B}^*:\;\; \calU(\YY_B) \longrightarrow \calU(\Aprin).
\end{equation}
Consider the torus $(\QQ^\times)^r$ with its character lattice identified with $\ZZ^r$, and consider the
$(\QQ^\times)^r$ action on $\FF(\Aprin)$ such that the variables 
in $({\bf u}_{t_0}, {\bf p})$ have the corresponding 
columns of  $(I_r, -B)$ as their weights.  Then (see \cite[\S 7]{GHKK})  $\widetilde{p}_{\!{}_B}^*$
identifies $\calU(\YY_B)$ with 
the subring of all  $(\QQ^\times)^r$-invariant elements of $\calU(\Aprin)$, or, equivalently, all
$\wt{u} \in \calU(\Aprin)$ that can be written as Laurent polynomials in $\widehat{{\bf y}}_t$
for every $t \in \TT_r$. 

Let $\calU^{\rm mono}(\YB)$ be the set of all global monomials on  $\YB$.
Part 2) of the following \prref{:YB-to-Aprin-admi} is \cite[(3) of Lemma 7.10]{GHKK}). 

\bpr{:YB-to-Aprin-admi} Under the injective ring homomorphism $\widetilde{p}_{\!{}_B}^*$,

1) $\calU^{\rm admi}(\YY_B)$ is mapped  bijectively to the set of all $(\QQ^\times)^r$-invariant elements of 
$\calU^{\rm admi}(\Aprin)$;

2) $\calU^{\rm mono}(\YB)$ is mapped bijectively to the set of all $(\QQ^\times)^r$-invariant elements of 
$\calU^{\rm mono}(\Aprin)$.
\epr

\begin{proof}
Recall from \reref{:YB-to-Aprin-trop} that each assignment
$\delta^\sv = \{\delta_t^\sv \in \ZZ^r_{\rm row}\}_{t \in \TT_r}$ gives the assignment
$\beta_{B^T}(\delta^\sv) = \{(\delta^\sv_t B_t^T, \, \delta_t^\sv C_t^T) \in \ZZ^{2r}_{\rm row}\}_{t \in \TT_r}$,
and that $\delta^\sv \in \AA_{B^T}(\Zmax)$ if and only if $\beta_{B^T}(\delta^\sv) \in \Yprin(\Zmax)$.
Using \eqref{eq:p-star}, one checks directly that for $\delta^\sv \in \ABTZ$ and $y \in \calU(\YB)$, one has 
$y \in \calU^{\rm admi}_{\delta^\sv}(\YB)$ 
if and only if 
$\widetilde{p}_{\!{}_B}^*(y) \in \calU^{\rm admi}_{\beta_{B^T}(\delta^\sv)}(\Aprin)$. It is also clear that 
$y \in \calU^{\rm mono}(\YB)$ if and only if $\widetilde{p}_{\!{}_B}^*(y) \in \calU^{\rm mono}(\Aprin)$.
\end{proof}

Let $\Delta^+(\AA_{B^T}, \ZZ)$ be the set\footnote{Our
$\Delta^+(\AA_{B^T}, \ZZ)$ is $\Delta^+_{\YY_B}(\ZZ)$ in the notation of \cite[Definition 7.9]{GHKK}.} of 
all $\{\delta^\sv_t \in \ZZ^r_{\rm row}\}_{t \in \TT_r} \in \AA_{B^T}(\Zmax)$ such that 
$-\delta^\sv_t B_t^T \geq 0$ for some $t \in \TT_r$. Alternatively, with 
$\Delta^+({\YY_{B^T}}, \ZZ) \subset \YY_{B^T}(\Zmax)$ given in \eqref{eq:Delta-plus-YBZ}, one has 
\begin{equation}\label{eq:Delta-Delta}
\Delta^+(\AA_{B^T}, \ZZ) = p_{\!{}_{B^T}}^{-1}(\Delta^+({\YY_{B^T}}, \ZZ)) \subset \ABTZ,
\end{equation}
where
$p_{\!{}_{B^T}}:  \ABTZ \rightarrow \YY_{B^T}(\Zmax),  \{\delta^\sv_t\}_{t \in \TT_r}
\mapsto \{\delta^\sv_t B_t^T\}_{t \in \TT_r}$ is the ensemble map.

\bld{:Y-g-vector}
Let $y \in \calU^{\rm mono}(\YB)$. Let
 $t \in \TT_r$ and ${\bf m} \in \ZZ^r$ be such that $y  = {\bf y}_t^{\bf m}$, and let 
$\delta^\sv(y)$ be the unique tropical point of $\ABT$ whose coordinate vector at $t$ is 
$-{\bf m}^T$.  Then the map
\begin{equation}\label{eq:mono-g-vec-YB}
\calU^{\rm mono}(\YY_B) \longrightarrow \Delta^+(\AA_{B^T}, \ZZ), \;\; y \longmapsto \delta^\sv(y),
\end{equation}
is a well-defined bijection and $y$ is $\delta^\sv(y)$-admissible. We call $\delta^\sv(y) \in \AA_{B^T}(\Zmax)$ the $g$-vector of $y \in \calU^{\rm mono}(\YB)$
\cite[Definition 5.10]{GHKK}.
\eld

\begin{proof}
The statements follow from \leref{:Y-global-mono}, \ldref{:admi-prin}, \prref{:YB-to-Aprin-admi}, and the definitions.
\end{proof}

\subsection{Global variables of $\bfY$-spaces}\label{ss:d-g}
Consider again  the Fock-Goncharov dual pair $(\AA_{B^T}, \YY_{B})$.

\bde{:global-Y-variables}
{\rm  
Elements in $\fX^{\rm global}(\YB) \stackrel{\rm def}{=}\fX(\YB) \cap \calU(\YB)$ are called {\it global variables of $\YB$}.
Global variables of $\YB$ are thus 
global monomials on  $\YB$ which are also $\bfY$-variables.
}
\ede

As we have already seen in \reref{:not-in}, not every $y \in \fX(\YB)$ is global. 
Denote again the variables of $\ABT$ and of $\YB$ by $x_{t; i}^\sv \in \fX(\ABT)$ and 
$y_{t; i} \in \fX(\YB)$, for $t \in \TT_r$ and $i \in [1, r]$. By \cite[Theorem 7.6]{ckq_2022} one has the well-defined map 
\[
\mathfrak{X}(\YY_{B}) \longrightarrow  \mathfrak{X}(\AA_{B}),\;\; y_{t;i}\longmapsto x_{t;i}, \hs t \in \TT_r, \, i \in [1, r],
\]
and by \cite[Lemma 4.13]{cao-gyoda:bongartz} one has the 
well-defined map 
\[
\mathfrak{X}(\AA_{B}) \longrightarrow \mathfrak{X}(\AA_{B^T}), \;\;x_{t; i}\longmapsto x_{t;i}^\sv, \hs t \in \TT_r, \,i \in [1, r].
\]
Combining the two, one thus has the well-defined map 
\begin{equation}\label{eq:lambda-y-x}
\lambda: \;\; \mathfrak{X}(\YY_{B}) \longrightarrow  \mathfrak{X}(\AA_{B^T}), \;\; y_{t; i} \longmapsto x^\sv_{t; i},
\hs t \in \TT_r, \,i \in [1, r].
\end{equation}

Recall from \ldref{:x-degree} and \ldref{:y-degree} that each $x^\sv \in \fX(\ABT)$ has a
$d$-tropical point $d_{x^\sv} \in \ABTZ$ and each $y \in \fX(\YB)$ has a $d$-tropical point $d_y \in \YBZ$. 
On the other hand, by \ldref{:globals-A-are-admissible} every $x^\sv \in \fX(\ABT)$ has a $g$-vector $\rho(x^\sv) \in \YBZ$, and
by \ldref{:Y-g-vector}
every $y \in 
\mathfrak{X}^{\rm global}(\YY_{B})$ has a $g$-vector $\delta^\sv(y) \in \ABTZ$. 

\ble{:d-g}
For every $y \in \fX(\YB)$ one has $d_y = \rho(\lambda(y)) \in \YBZ$, and for every $y \in
\mathfrak{X}^{\rm global}(\YY_{B})$ one has $\delta^\sv(y) = d_{\lambda(y)} \in \ABTZ$.
\ele

\begin{proof}
Let $y \in \fX(\YB)$ and suppose that  $y = y_{t; i}$ for $t \in \TT_r$ and $i \in [1, r]$. Then $\lambda(y) = 
x^\sv_{t; i}$, and the two tropical points $d_y$ and $\rho(\lambda(y))$ of $\YB$ have the same coordinate vector 
$-e^i$ at $t$. Thus
$d_y = \rho(\lambda(y))$. Similarly, $\delta^\sv(y) = d_{\lambda(y)}$ if 
$y \in \mathfrak{X}^{\rm global}(\YY_{B})$.
\end{proof}

\bco{:injective}
The restriction of $\lambda: \mathfrak{X}(\YY_{B}) \to \mathfrak{X}(\AA_{B^T})$ to 
$\mathfrak{X}^{\rm global}(\YY_{B}) \subset \fX(\YB)$   is injective.
\eco

\begin{proof} 
Suppose that $y, y' \in \fX^{\rm global}(\YB)$ are such that $\lambda(y) = \lambda(y')$. 
Then $d_{\lambda(y)} = d_{\lambda(y')}$, and thus
$\delta^\sv(y) = \delta^\sv(y')$ by \leref{:d-g}. By \ldref{:Y-g-vector}, one has $y = y'$. 
\end{proof}

We show in \leref{:Y-admi-finite} that $\lambda: \fX^{\rm global}(\YB) \to \fX(\ABT)$ is a bijection if $B$  is of finite type.

\subsection{Theta functions are admissible}\label{ss:theta-admi}
In \cite{GHKK}, for every $r \times r$ mutation matrix $B$ 
the authors defined subsets  
\[
\Theta(\Aprin) \subset \Yprin(\Zmax), \hs
\Theta(\AA_{B^T}) \subset \YY_B(\ZZ^{\max}), \hs \mbox{and} \hs \Theta(\YY_B) \subset \AA_{B^T}(\ZZ^{\max}),
\]
along with maps, which we will all denote as 
$\theta:  \Theta(V) \rightarrow \calU(V),  q \mapsto \theta_q \stackrel{\rm def}{=} \nu(\vartheta_q)$,
where $V = \AA_{B^T}, \YY_B$ or $\Aprin$, and 
$\nu$ and $\vartheta_q$ are as in  \cite[Theorem 0.3]{GHKK}. 
The elements $\theta_q$ are called {\it theta functions} on their respective spaces.
It follows by their construction in \cite[\S 7]{GHKK} (see also \cite[Proposition 3.6]{GHKK} and \cite[Theorem A.1.3 and
Theorem A.1.4]{qin-2022}) that the theta functions on $\Aprin$ are admissible. By \cite[Theorem 7.16]{GHKK}, 
$\theta_\rho$ is $\rho$-admissible for every $\rho \in \Theta(\ABT)$. 
By the construction in  \cite[Construction 7.11]{GHKK}
and \prref{:YB-to-Aprin-admi}, 
$\theta_{\delta^\sv}$ is $\delta^\sv$-admissible for $\delta^\sv \in \Theta(\YB)$.

\bco{:admissible-exist}
If $B$ is mutation equivalent to an acyclic mutation matrix, then 

1) $\rho$-admissible functions on $\AA_{B^T}$ exist for every $\rho\in \YY_{B}(\ZZ^{\max})$;

2) $\delta^\sv$-admissible funtions on $\YY_{B}$ exist for every $\delta^\sv\in \AA_{B^T}(\ZZ^{\max})$.
\eco

\begin{proof}
We may assume that $B$ is acyclic. By \cite[Proposition 8.24 and Proposition 8.25]{GHKK}, one has $\Theta(\Aprin) = 
\Yprin$, and by \cite[Definition 7.12 and Definition-Lemma 7.14]{GHKK}, 
\[
\Theta(\AA_{B^T}) = \YY_B(\Zmax) \hs \mbox{and} \hs 
\Theta(\YY_B) = \AA_{B^T}(\Zmax).
\]
As theta functions are admissible, the statements are proved.
\end{proof}

\bre{:large-complex}
{\rm
The same proof shows that \coref{:admissible-exist} holds if $\AA_B$ has {\it large cluster complex} 
\cite[Definition 8.23 and Proposition 8.25]{GHKK}.
\hfill $\diamond$
}
\ere

\section{Fock-Goncharov pairing from the perspective of admissible functions}\label{s:FG-pairing}

\subsection{A question}\label{ss:a-question} Motivated by \cite[Conjecture 4.3]{FG:ensembles}, we ask the following question
on a mutation matrix $B$ in terms of the Fock-Goncharov dual pair $(\AA_{B^T}, \YY_B)$.

\bqu{:FG-pairing-conj}
What mutation matrices $B$ have {\it Property $\mathfrak{P}$}, in the sense that 
$\calU^{\rm admi}_\rho(\AA_{B^T})\neq \emptyset$ and $\calU^{\rm admi}_{\delta^\sv}(\YY_{B}) \neq \emptyset$
for every $\rho \in \YY_B(\Zmax)$ and $\delta^\sv \in \AA_{B^T}(\Zmax)$ and that
\begin{equation}\label{eq:val-val-4}
\val_{\delta^\sv}(x^\sv) = \val_\rho(y)
\end{equation}
for any  $x^\sv \in \calU^{\rm admi}_\rho(\AA_{B^T})$ and $y \in \calU^{\rm admi}_{\delta^\sv}(\YY_{B})$? 
\equ

If the mutation matrix $B$ has Property $\mathfrak{P}$,  one then has a 
well-defined {\it Fock-Goncharov pairing}
\[
        \langle \cdot , \cdot \rangle : \;\; \AA_{B^T}(\ZZ^{\max}) \times \YY_B(\ZZ^{\max}) \longrightarrow \ZZ, \;\; 
        \langle \delta^\sv, \rho\rangle \stackrel{\text{def}}{=} {\rm val}_{\delta^\sv}(x^\sv) = {\rm val}_{\rho}(y)
\]
by taking any $x^\sv \in \calU^{\rm admi}_{\rho}(\AA_{B^T})$ and 
$y \in \calU^{\rm admi}_{\delta^\sv}(\YY_B)$ for $\delta^\sv \in \AA_{B^T}(\Zmax)$ and $\rho \in \YY_B(\Zmax)$.

\bre{:theta-pairing}
{\rm 
Recall from $\S$\ref{ss:theta-admi}  the theta functions given by the maps
\[
\Theta(\ABT) \longrightarrow \calU(\AA_{B^T}), \;\; \rho \longmapsto \theta_\rho \hs \mbox{and} \hs
\Theta(\YB) \longrightarrow \calU(\YY_B), \;\; \delta^\sv \longmapsto \theta_{\delta^\sv}.
\]
The {\it theta pairing conjecture} in \cite[Remark 9.11]{GHKK} says that
${\rm val}_{\delta^\sv}(\theta_\rho) = {\rm val}_{\rho}(\theta_{\delta^\sv})$ for all
$\delta^\sv \in \Theta(\YB)\subset \AA_{B^T}(\ZZ^{\max})$ and 
$\rho \in \Theta(\ABT) \subset \YY_B(\ZZ^{\max})$.
The theta pairing conjecture for skew-symmetric mutation matrices has recently been proved in \cite{Muller_talk}.
\hfill $\diamond$
}
\ere

We now discuss the values of the two expressions in \eqref{eq:val-val-4}. Following \cite[Definition 9.1 and Lemma 9.3]{GHKK}, for
\begin{equation}\label{eq:rho-delta-sv}
\rho = \{\rho_t \in \ZZ^r_{\rm row}\}_{t \in \TT_r} \in\YY_B(\mathbb Z^{\max}) \hs \mbox{and} \hs
\delta^\sv = \{\delta_t^\sv \in \ZZ^r_{\rm row}\}_{t \in \TT_r} \in \AA_{B^T}(\mathbb Z^{\max}),
\end{equation}
we say that 
$\rho$ is {\it optimized at $w \in \TT_r$} if $-\rho_w \geq 0$, and that $\delta^\sv$ is 
{\it optimized at $w \in \TT_r$} if $-\delta^\sv_w B_w^T \geq 0$.

\ble{:val-rho-delta}
Let $\rho$ and $\delta^\sv$ be as in \eqref{eq:rho-delta-sv}, and let  $x^\sv \in \calU^{\rm admi}_{\rho}(\AA_{B^T})$
and $y \in \calU^{\rm admi}_{\delta^\sv}(\YY_B)$. Then

1) ${\rm val}_{\delta^\sv}(x^\sv) \geq -\delta_t^\sv \rho_t^T$ for every $t \in \TT_r$; if $\delta^\sv$ is 
optimized at $w \in \TT_r$, then ${\rm val}_{\delta^\sv}(x^\sv) = -\delta_w^\sv \rho_w^T$;

2) ${\rm val}_{\rho}(y) \geq -\delta_t^\sv \rho_t^T$ for every $t \in \TT_r$; if $\rho$ is optimized at $w \in \TT_r$, then 
${\rm val}_{\rho}(y) = -\delta_w^\sv \rho_w^T$.
\ele

\begin{proof}
1) Let $t \in \TT_r$ be arbitrary and write 
$x^\sv = ({\bf x}_t^\sv)^{-\rho_t^T} (1+L_t)$ as in \eqref{eq:a-g-x}, so $L_t$ is a 
sum of Laurent monomials of the form  $({\bf x}_t^\sv)^{B_t^T{\bf u}}$ with ${\bf u}\in\mathbb Z_{\geq 0}^r$.
 Then
\[
{\rm val}_{\delta^\sv}(x^\sv) = -\delta_t^\sv \rho_t^T + \val_{\delta^\sv}(1+L_t) 
= -\delta_t^\sv \rho_t^T + {\rm max}\{0, \val_{\delta^\sv}(L_t)\}  \geq -\delta_t^\sv \rho_t^T.
\]
Suppose that $\delta^\sv$ is optimized at $w \in \TT_r$. Then $-{\delta_w^\sv}B_w^T{\bf u}\geq 0$ for every 
monomial term $({\bf x}_w^\sv)^{B_w^T{\bf u}}$ of $L_w$. Thus 
${\rm val}_{\delta^\sv}(x^\sv)=-\delta_w^\sv \rho_w^T + \val_{\delta^\sv}(1+L_w)  = -\delta_w^\sv \rho_w^T.$
One proves 2) similarly.
\end{proof}

The following  special case of \eqref{eq:val-val-4}, which has been noted in \cite[Remark 9.11]{GHKK} for theta functions
(see also \cite[4 of Conjecture 4.3 and Theorem 5.2]{FG:ensembles}),
is key to our study in 
$\S$\ref{s:ca-and-trop-pts} 
of tropical friezes and cluster-additive functions. 

\bpr{:duality-pairing} Let $x^\sv \in \calU^{\rm admi}_{\rho}(\AA_{B^T})$
and $y \in \calU^{\rm admi}_{\delta^\sv}(\YY_B)$, where $\rho  \in\YY_B(\mathbb Z^{\max})$
and $\delta^\sv \in \AA_{B^T}(\mathbb Z^{\max})$.
If $x^\sv$ is a global monomial on $\ABT$, or if $y$ is a global monomial on $\YY_B$, then 
\begin{equation}\label{eq:val-val}
{\rm val}_{\delta^\sv}(x^\sv)={\rm val}_{\rho}(y) = {\rm max}\{-\delta_t^\sv\rho_t^T: t \in \TT_r\},
\end{equation}
and the maximum is achieved at any $w \in \TT_r$ such that either $x^\sv$ or $y$ is a Laurent monomial at $w$.
\epr

\begin{proof} 
Suppose that $x^\sv$ is a global monomial on $\AA_{B^T}$ and $y \in \calU^{\rm admi}_{\delta^\sv}(\YY_B)$. 
By \leref{:A-global-mono}, $x^\sv=({\bf x}^\sv_w)^{-\rho_w^T}$ for some
$w\in\mathbb T_r$ with $-\rho_w\geq 0$. Then
${\rm val}_{\delta^\sv}(x^\sv)={\rm val}_{\delta^\sv}({\bf x}_w^{-\rho_w^T})=-\delta_w^\sv\rho_w^T$.
Since $\rho$ is optimized at $w$, by \leref{:val-rho-delta} one has $\val_\rho(y) = -\delta_w^\sv\rho_w^T=
{\rm val}_{\delta^\sv}(x^\sv)$. 

Suppose that  $y$ is a global monomial on $\YY_B$ and
$x^\sv \in \calU^{\rm admi}_\rho(\AA_{B^T})$. By \leref{:Y-global-mono}, $y={\bf y}_w^{-(\delta_w^\sv)^T}$ for some
 $w \in \TT_r$ and $\delta^\sv$ is optimized at $w$. By \leref{:val-rho-delta} again, 
${\rm val}_{\delta^\sv}(x^\sv)=-\delta_w^\sv\rho_w^T = \val_\rho(y)$.
\end{proof}

Recall from \reref{:admissible-1} that when $B$ has full rank, 
any $x^\sv \in \calU^{\rm admi}(\ABT)$ is $\rho$-admissible for a unique $\rho \in \YBZ$. 
We do not know whether or not such a uniqueness statement holds in general, 
but we have the following consequence of \prref{:duality-pairing}.

\bco{:rho-unique}
For a mutation matrix $B$, if the semi-field $\FF_{>0}(\YB)$ has a set of generators consisting of global monomials on $\YB$, then
every $x^\sv \in \calU^{\rm admi}(\ABT)$ is $\rho$-admissible for a unique $\rho\in \YBZ$. 
\eco

\begin{proof}
Let $\{z_1, \ldots, z_l\}$ be a set of generators of 
$\FF_{>0}(\YY_B)$ consisting of global monomias on $\YB$, and for 
$i \in [1, l]$, let $\delta_i^\sv \in \ABTZ$ be the $g$-vector of $z_i$.
Suppose that
$x^\sv \in \calU^{\rm admi}_{\rho}(\AA_{B^T}) \cap \calU^{\rm admi}_{\rho'}(\AA_{B^T})$ for $\rho, \rho' \in \YY_B(\Zmax)= {\rm Hom}_{\rm sf}(\FF_{>0}(\YY_B),\mathbb Z^{\rm max})$.
Then by \prref{:duality-pairing},  one has 
\[
\rho(z_i)=\val_{\rho}(z_i) = \val_{\delta_i^\sv}(x^\sv) =\val_{\rho'}(z_i)=\rho'(z_i)
\]
for every $i \in [1, l]$. As $\{z_1, \ldots, z_l\}$ generates 
$\FF_{>0}(\YY_B)$, one has $\rho = \rho'$.
\end{proof}

\bre{:admi-unique}
{\rm
Under the condition on the mutation matrix $B$ in \coref{:rho-unique},  
one has a well-defined  map
\[
\calU^{\rm admi}(\ABT) \longrightarrow \YBZ, \;\;x^\sv \longmapsto \rho(x^\sv),
\]
such that $x^\sv$ is $\rho(x^\sv)$-admissible
for $x^\sv \in \calU^{\rm admi}(\ABT)$, generalizing the assignment of the $g$-vector to
a cluster monomial on $\ABT$. In contrast, \reref{:delta-unique} says that one has the well-defined map
\[
\calU^{\rm admi}(\YB) \longrightarrow \ABTZ, \;\; y \longmapsto \delta^\sv(y),
\]
for every mutation matrix $B$, where  $y$ is $\delta^\sv(y)$-admissible
for $y \in \calU^{\rm admi}(\YB)$. 
\hfill $\diamond$
}
\ere

We will see in \reref{:xim-yim} that when $B$ is acyclic, there exists a free generating set of 
$\FF_{>0}(\YY_B)$ consisting of global variables of $\YY_B$. We thus have the following consequence.

\bco{:acyclic-rho-unique}
Assume that $B$ is mutation equivalent to an acyclic mutation matrix. Then every $x^\sv \in \calU^{\rm admi}(\ABT)$ is $\rho$-admissible for a unique $\rho\in \YBZ$. 
\eco



\subsection{Fock-Goncharov pairing in finite type}\label{ss:finite-case}
In this subsection, assuming that $B$ is of finite type, i.e., $\AA_B$, equivalently, 
$\ABT$, has finitely many cluster variables, we prove in 
\thref{:finite-admi-mono} that admissible functions on $\ABT$ and on
$\YB$ are precisely all the global monomials, and we recover in \coref{:duality-pairing-finite-type} 
the well-defined Fock-Goncharov pairing established in \cite[$\S$5.1]{FG:ensembles}. 

For any mutation matrix $B$, recall
that $\calU^{\rm mono}(\ABT)$ denotes the set of all cluster monomials on $\ABT$ and $\calU^{\rm mono}(\YB)$
that of all global monomials on $\YB$.
Recall also that the $g$-vector of $x^\sv \in \calU^{\rm mono}(\ABT)$ is denoted as $\rho(x^\sv) \in \YBZ$,
and  the $g$-vector of $y \in \calU^{\rm mono}(\YB)$ is denoted as
$\delta^\sv(y) \in \ABTZ$. 
It is well-known \cite{HPS_2018} that  the $g$-vector fan 
of $\ABT$ (see \reref{:g-fan}) is complete if $B$ is of finite type.

\ble{:finite-bijection-1}
When $B$ is of finite type, one has bijections
\begin{align*}
&\calU^{\rm mono}(\AA_{B^T}) \longrightarrow  \YBZ,\;\; x^\sv \longmapsto \rho(x^\sv),\\
&\calU^{\rm mono}(\YY_B) \longrightarrow \ABTZ, \;\; y \longmapsto \delta^\sv(y).
\end{align*}
\ele

\begin{proof}
Recall  from $\S$\ref{ss:admissible-ABT} 
that $\Delta^+(\YB, \ZZ)$ is 
the set of the $g$-vectors of all cluster monomials on  $\ABT$.
As the $g$-vector fan of $\ABT$ is complete in this case, 
$\Delta^+(\YB, \ZZ) = \YBZ$.
By \eqref{eq:Delta-Delta}, one also has
$\Delta^+(\ABT, \ZZ) = \ABTZ$, where $\Delta^+(\ABT, \ZZ)$ is the set of the $g$-vectors of all global
monomials on  $\YB$. \leref{:finite-bijection-1} now follows from \leref{:Delta-plus} and \ldref{:Y-g-vector}.
\end{proof}

\bth{:finite-admi-mono}
When the mutation matrix $B$ is of finite type, one has
\[
\calU^{\rm admi}(\AA_{B^T}) = \calU^{\rm mono}(\AA_{B^T}) \hs \mbox{and} \hs 
\calU^{\rm admi}(\YY_B) = \calU^{\rm mono}(\YY_B).
\]
\eth

\begin{proof}
Let $\rho \in \YBZ$ and assume that $z^\sv \in \calU(\ABT)$ is $\rho$-admissible. Since $B$ is of finite type and by \leref{:finite-bijection-1},
there exists a cluster monomial
$x^\sv$ of $\ABT$ whose $g$-vector is $\rho$. It suffices to show  $z^\sv = x^\sv$. 

Let $a^\sv \in \fX(\ABT)$ and let $d_{a^\sv}\in \ABTZ$ be the $d$-tropical point of $a^\sv$ defined in
\ldref{:x-degree}. One has $\val_{d_{a^\sv}}(u^\sv)=(a^\sv \, \|\, u^\sv)_d$ for any $u^\sv\in\mathcal U^+(\AA_{B^T})\setminus\{0\}$. 
 By \leref{:finite-bijection-1},  there exists a global monomial $y \in \calU^{\rm mono}(\YB)$ whose $g$-vector is  $d_{a^\sv}$.
Then by \prref{:duality-pairing}, one has 
$\val_{d_{a^\sv}}(z^\sv) = \val_\rho(y) = \val_{d_{a^\sv}}(x^\sv)$.
Hence,
$(a^\sv \, \|\, z^\sv)_d = (a^\sv \, \|\, x^\sv)_d$ for any $a^\sv \in \fX(\ABT)$.

Since $B$ is of finite type, cluster monomials on  $\ABT$ form an atomic $\ZZ$-basis of $\calU(\ABT)$ (see \cite[Theorem 10.2]{FT:orbifold-bases}). Thus the element $z^\sv$ admits a unique expression $z^\sv = x^\sv_1  + \ldots + x^\sv_q$, where each
$x_j^\sv$ is a cluster monomial on  $\ABT$. It follows that 
\begin{equation}
(a^\sv \, \| \, x^\sv)_d = (a^\sv \, \| \, z^\sv)_d = 
{\rm max}\{(a^\sv \, \| \, x_1^\sv)_d, \, \ldots, \, (a^\sv \, \| \, x_q^\sv)_d\}, \hs \forall\; a^\sv \in \fX(\ABT).
\end{equation}
Thus $(a^\sv \, \| \, x^\sv)_d \geq (a^\sv \, \| \, x_j^\sv)_d$ for every $a^\sv \in \fX(\ABT)$ and every $j \in [1, q]$.
By \leref{:u-v}, $x^\sv = x_j^\sv$ for each $j \in [1, q]$. It follows that $q = 1$ and $z^\sv = x^\sv$. Hence, $\calU^{\rm admi}(\AA_{B^T}) = \calU^{\rm mono}(\AA_{B^T})$.

Assume now that $y \in \calU^{\rm admi}(\YB)$. By  \prref{:YB-to-Aprin-admi},
$\widetilde{p}_{\!{}_B}^*(y) \in \calU^{\rm admi}(\Aprin)$, where the map $\widetilde{p}_{\!{}_B}^*$ is given in
\eqref{eq:p-YB-Aprin}. Let $\phi: \calU(\Aprin) \to \calU(\AA_B)$ be the coefficient specialization
which sends every cluster variable $u_{t; i}$ of $\Aprin$ to the corresponding cluster variable $x_{t; i}$ of $\AA_B$
and every frozen variable $p_j$ to $1$. Then $\phi(\widetilde{p}_{\!{}_B}^*(y)) \in \calU^{\rm admi}(\AA_B)$, so
$\phi(\widetilde{p}_{\!{}_B}^*(y))$ is a cluster monomial on  $\AA_B$, which implies that  $\widetilde{p}_{\!{}_B}^*(y)$
is an extended cluster monomial on $\Aprin$. By \prref{:YB-to-Aprin-admi} again, $y$ must be a global monomial on  $\YB$. Thus $\calU^{\rm admi}(\YY_B) = \calU^{\rm mono}(\YY_B)$.
\end{proof}

\bco{:duality-pairing-finite-type} \cite[$\S$5.1]{FG:ensembles}
Every mutation matrix $B$ of finite type has Property ${\mathfrak{P}}$, and one has the
Fock-Goncharov pairing $\langle \cdot , \cdot \rangle:\,  \AA_{B^T}(\ZZ^{\max})\times 
\YY_{B}(\ZZ^{\max}) \rightarrow \ZZ$ given by
\[ 
\langle \delta^\sv, \rho \rangle \stackrel{\text{def}}{=}  {\rm val}_{\delta^\sv}(x^\sv) =\val_\rho(y)
= {\rm max}\{-\delta_t^\sv\rho_t^T: t \in \TT_r\},
\]
where  for $\delta^\sv= \{\delta_t^\sv \in \ZZ^r_{\rm row}\}_{t \in \TT_r} \in  \AA_{B^T}(\ZZ^{\max})$,  $y$ 
is the unique global monomial on $\YB$ with $g$-vector $\delta^\sv$, and
for $\rho = \{\rho_t \in \ZZ^r_{\rm row}\}_{t \in \TT_r} \in \YY_{B}(\ZZ^{\max})$,
$x^\sv$ is the unique cluster monomial on $\ABT$ with $g$-vector $\rho$. 
\eco
\begin{proof}
    This follows from \thref{:finite-admi-mono} and \prref{:duality-pairing}.
\end{proof}

\section{Tropical friezes and cluster-additive functions}\label{s:ca-and-trop-pts}
Fix an $r \times r$ symmetrizable generalized Cartan matrix $A$ and let $B = B_A$
as in \eqref{eq:BA-intro}. We then have the Fock-Goncharov dual pair 
$\OO_A = (\AA_{B^T}, \YY_B)$. 
In this section, we give realizations of elements in 
$\CfriezeAT$ and $\CaddA$ using tropical points of and admissible functions on $\ABT$ and $\YB$. 
We also generalize Ringel's definition of cluster-hammock functions which was defined in \cite{Ringel:additive}
only for $A$  of finite and simply-laced type, and
we show that  cluster-hammock functions are also tropical friezes. 

\subsection{Generic frieze patterns}\label{ss:patterns-and-trop} With $B = B_A$ as in \eqref{eq:BA-intro}, we write
\[
    \ABT = (\FF(\ABT), \{{\bf x}^\sv_t = 
    (x^\sv_{t; 1}, \ldots, x^\sv_{t; r})\}_{t \in \TT_r}) \;\;\; \mbox{and} \;\;\;
    \YY_B = (\FF(\YY_B), \{{\bf y}_t = (y_{t; 1}, \ldots, y_{t, r})\}_{t \in \TT_r}).
\]
Set $t(1, 0) = t_0$ and consider  the 2-regular sub-tree $\TT_r^\flat$
of $\TT_r$ given by (recall \eqref{eq:right} and \eqref{eq:left})
\[
\begin{xy}
(0,1)*+{\cdots}="A1",(15,1)*+{t(r, -2)}="A2",(33.5,1)*+{t(1, -1)}="A3",(48.5,1)*+{\cdots}="A4",(64,1)*+{t(r, -1)}="A5",
(85,1)*+{t(1,0) = t_0}="A6",(102.5,1)*+{\cdots}="A7",(116,1)*+{t(r,0)}="A8",(132,1)*+{t(1, 1)}="A9",(145.5,1)*+{\cdots}="A10",
\ar@{-}^{r-1\;\;\;\;}"A1";"A2",\ar@{-}^{r}"A2";"A3",\ar@{-}^{\;\;\;\;1}"A3";"A4", 
\ar@{-}^{r-1\;\;\;\;}"A4";"A5",\ar@{-}^{r\;\;\;}"A5";"A6",\ar@{-}^{\;\;\;\;\;\;\;1}"A6";"A7",
\ar@{-}^{r-1\;\;\;}"A7";"A8",\ar@{-}^{r}"A8";"A9",\ar@{-}^{\;\;\;1}"A9";"A10"
\end{xy}
\]
As we have done in $\S$\ref{ss:thmA-intro}, introduce, for $(i, m) \in [1, r] \times \ZZ$, 
\begin{equation}\label{eq:xim-yim-1}
x^\sv(i,m) = x^\sv_{t(i,m); \,i} \in \fX(\ABT) \hs \mbox{and} \hs 
y(i,m) = y_{t(i,m); \,i} \in \fX(\YY_B).
\end{equation}

\ble{:B-im-column}
For every $(i, m) \in [1, r] \times \ZZ$, one has $B_{t(i, m+1)} = B_{t(i, m)}$, and
 the $i^{\rm th}$ column of  
$B_{t(i, m)}$ 
is $(-a_{1,i}, \ldots, -a_{i-1,i}, 0, -a_{i+1,i}, \ldots, -a_{r,i})^T\geq 0$. Consequently, $y(i, m)$ is a global variable
of $\YY_B$ for every $(i, m) \in [1, r] \times \ZZ$.
\ele

\begin{proof}
The first statement is proved  by a direct calculation (see also \cite{GHL:friezes}). The second statement follows from 
the first statement and 
\leref{:Y-global-mono}.
\end{proof}

\bpr{:xim-yim} 
For every $(i, m) \in [1, r] \times \ZZ$, one has
\begin{align}\label{eq:xsv-im}
    &x^\sv(i,m) \, x^\sv(i,m+1) = 1 + \prod_{j=i+1}^r x(j,m)^{-a_{i, j}} \; \prod_{j=1}^{i-1} x(j,m+1)^{-a_{i, j}},\\
\label{eq:y(m,i)-relations} 
&y(i,m)\, y(i,m+1) = \prod_{j=i+1}^r \left(1 + y(j,m)\right)^{-a_{j,i}} \, \,\prod_{j=1}^{i-1} (1 + y(j,m+1))^{-a_{j,i}}.
\end{align}
\epr

\begin{proof}
The equation \eqref{eq:xsv-im} is proved in 
\cite[Lemma 5.4]{GHL:friezes}.
To prove \eqref{eq:y(m,i)-relations}, for $i \in [1, r]$ and $m \in \ZZ$,
let $\Sigma_{t(i,m)} = ({\bf y}_{t(i,m)}, B_{t(i,m)})$ and set 
${\bf y}_{t(1,m)} = (y_1^{(m)}, \dots, y_r^{(m)})$.
Since $\Sigma_{t(i,m)} =  \mu_{i-1} \cdots \mu_2\mu_1(\Sigma_{t(1,m)})$, by \eqref{eq:Y-mut} and \leref{:B-im-column}, we have
\begin{equation}\label{eq:Y-cluster-to-slice-1}
y_i^{(m)}=  y(i,m) \prod_{j=1}^{i-1} (1 + y(j,m))^{a_{j,i}}, \quad i \in [1, r].
\end{equation}
On the other hand, we have $\Sigma_{t(i,m)} = \mu_i \mu_{i+1} \cdots \mu_r(\Sigma_{t(1,m+1)})$. Applying \eqref{eq:Y-mut} and by the fact that the j$^{th}$ variable of $\mu_j(\Sigma_{t(j,m)})$ is $y(j,m)^{-1}$, we get
\begin{equation}\label{eq:Y-cluster-to-slice-2}
y_i^{(m+1)} =y(i,m)^{-1} \prod_{j=i+1}^{r} (1 + y(j,m))^{-a_{j,i}}, \quad i \in [1, r].
\end{equation}
Replacing $m$ by $m+1$ in  \eqref{eq:Y-cluster-to-slice-1} and comparing with \eqref{eq:Y-cluster-to-slice-2}, one arrives at  
\eqref{eq:y(m,i)-relations}.
\end{proof}

\bde{:yim-xim}
1) The {\it generic frieze pattern of $\ABT$} is the map
\[
[1, r] \times \ZZ \longrightarrow \fX(\ABT), \;\; (i, m) \longmapsto x^\sv(i, m).
\]

2) The {\it generic $\bfY$-frieze pattern of $\YY_B$} is the map
\[
[1, r] \times \ZZ \longrightarrow \fX(\YY_B), \;\; (i, m) \longmapsto y(i, m).
\]
We will also denote the two generic frieze patterns  as the respective
assignments
\begin{equation}\label{eq:frieze-patrn-nota}
\mathfrak{X}^\flat(\ABT) = \{ x^\sv(i,m)\}_{(i,m) \in [1,r] \times \ZZ} \hs \mbox{and} \hs
\mathfrak{X}^\flat(\YY_B) = \{ y(i,m)\}_{(i,m) \in [1,r] \times \ZZ}.
\end{equation}
Abusing notation, we also regard $\mathfrak{X}^\flat(\ABT)$ and $\mathfrak{X}^\flat(\YY_B)$ as subsets of 
$\fX(\ABT)$ and $\fX(\YY_B)$ respectively. By \leref{:B-im-column}, one has $\mathfrak{X}^\flat(\YY_B)\subset \fX^{\rm global}(\YB)$.
\ede

\bre{:xim-yim}
For any $(i, m) \in [1, r] \times \ZZ$, it is proved in \cite[Lemma 5.4]{GHL:friezes} that
\[
{\bf x}^\sv_{t(i,m)}=(x^\sv(1,m+1),\,\ldots,\,x^\sv(i-1,m+1),\,x^\sv(i,m),\, x^\sv(i+1,m),\, \ldots,\,x^\sv(r,m)).
\]
In particular, we have ${\bf x}^\sv_{t(1, m)} = (x^\sv(1, m), \ldots, x^\sv(r, m))$ for every $m \in \ZZ$.
The $\bfY$-cluster ${\bf y}_{t(i,m)}$, however, is in general {\it not} formed by the variables from $\mathfrak{X}^\flat(\YY_B)$, as is already seen from the $A_2$ example in \reref{:not-in}. On the other hand,
it follows from \eqref{eq:Y-cluster-to-slice-1} that for each $m \in \ZZ$, the $r$-tuple
$(y(1, m), \ldots, y(r, m))$ of global variables is a set of free generators of $\FF(\YY_B)$ and of 
$\FF_{>0}(\YY_B)$, a fact that 
has been used in the proof of  \coref{:acyclic-rho-unique}. It then follows from \eqref{eq:y(m,i)-relations} that
for every $(i, m) \in [1, r] \times \ZZ$, the $r$-tuple
\[
(y(1,m+1),\,\ldots,\,y(i-1,m+1),\, y(i, m), \,  y(i+1,m),\, \ldots,\,y(r,m))
\]
is also a set of free generators of $\FF(\YY_B)$ and of 
$\FF_{>0}(\YY_B)$. 

\ere

\subsection{Two realizations of tropical friezes}\label{ss:dual-trop-frieze}
Recall that $C_{\rm frieze}(A^T)$ denotes 
the set of all tropical friezes associated to $A^T$. 
 For 
$\delta^\sv = \{\delta^\sv_t=(\delta^\sv_{t; 1}, \ldots, \delta^\sv_{t; r})\}_{t \in \TT_r} \in \ABT(\Zmax)$, define 
\begin{equation}\label{eq:f-delta-i-m}
f_{\sdeltasv}:\;\; [1, r] \times \ZZ \longrightarrow \ZZ, \;\; f_{\sdeltasv}(i, m) = 
\val_{\sdeltasv} (x^\sv(i, m)) = \delta^\sv_{t(i, m); i},
\end{equation}
where the last identity is due to the definition of $x^\sv(i, m)$ as the $i^{\rm th}$ variable of the cluster ${\bf x}^\sv_{t(i, m)}$.

\bpr{:trop-frieze-dual-1}
For any symmetrizable generalized Cartan matrix $A$ and with $B = B_A$ as in \eqref{eq:BA-intro}, 
$f_{\sdeltasv} \in C_{\rm frieze}(A^T)$ for every $\delta^\sv \in \ABT(\Zmax)$, and one has a bijection 
\begin{equation}\label{eq:f-deltasv}
\ABT(\Zmax) \longrightarrow C_{\rm frieze}(A^T), \;\;\delta^\sv \longmapsto f_{\sdeltasv}.
\end{equation}
\epr

\begin{proof}   By \prref{:xim-yim}, the variables $\{x^\sv(i, m)\}$ in $\fX^\flat(\ABT)$ satisfy \eqref{eq:xsv-im}.
Applying $\val_{\sdeltasv}$ to both sides of \eqref{eq:xsv-im}, one sees that
$f_{\sdeltasv} \in C_{\rm frieze}(A^T)$ for every $\delta^\sv \in \ABT(\Zmax)$. As 
$x^\sv_{t_0} = (x^\sv(1, 0), \ldots, x^\sv(r, 0))$ is a set of free generators of $\FF_{>0}(\ABT)$, the bijectivity of the map 
in \eqref{eq:f-deltasv} follows from \leref{:trop-intro}.
\end{proof}

Recall from \deref{:slice} that for $m \in \ZZ$, the 
$m^{\rm th}$ slice of any $f \in C_{\rm frieze}(A^T)$ is defined as 
\[
f|_{\{m\}} = (f(1, m), \, \ldots, \, f(r, m))^T \in \ZZ^r.
\] 
We have the  following description of the inverse of the
bijection in \eqref{eq:f-deltasv}.

\ble{:inverse-f-deltasv}
For any $\delta^\sv \in \ABTZ$ and any $m \in \ZZ$, one has  
$(\delta^\sv_{t(1, m)})^T=f_{\sdeltasv}|_{\{m\}}$.
\ele
\begin{proof}
This follows from \eqref{eq:f-delta-i-m} and the fact that  ${\bf x}^\sv_{t(1, m)} = (x^\sv(1, m), \ldots, x^\sv(r, m))$ (see \reref{:xim-yim}).
\end{proof}

For the second realization of elements in $C_{\rm frieze}(A^T)$,
for $y \in \calU^{\rm admi}(\YB)$ define
\begin{equation}\label{eq:f-ysv}
f_{y}: \;\; [1, r] \times \ZZ \longrightarrow \ZZ, \;\; f_{y}(i, m) = \val_{\rho(i, m)}(y),
\end{equation}
where for $(i, m) \in [1, r] \times \ZZ$, 
 $\rho(i, m) \in \YB(\Zmax)$ is the $g$-vector of $x^\sv(i, m) \in \fX(\ABT)$.

\bpr{:trop-frieze-dual-2}
For any symmetrizable generalized Cartan matrix and with $B = B_A$ as in \eqref{eq:BA-intro}, 
$f_{y} \in C_{\rm frieze}(A^T)$ for every $y \in \calU^{\rm admi}(\YB)$,  and one has a surjective map 
\begin{equation}\label{eq:f-ysv-1}
\calU^{\rm admi}(\YB) \longrightarrow C_{\rm frieze}(A^T), \;\; y \longmapsto f_{y}.
\end{equation}
Moreover, if $y \in \calU^{\rm admi}(\YB)$ and $\delta^\sv \in \ABT(\Zmax)$, then $f_{y} = f_{\sdeltasv}$
 if and only if $y$ is $\delta^\sv$-admissible. 
\epr

\begin{proof}
Suppose that $y \in \calU^{\rm admi}(\YB)$ is $\delta^\sv$-admissible for 
$\delta^\sv \in \ABT(\Zmax)$. By \prref{:duality-pairing},
\[
f_{y} (i, m) = \val_{\rho(i, m)}(y) = \val_{\sdeltasv}(x^\sv(i, m)) = f_{\sdeltasv}(i, m), 
\hs  (i, m) \in [1, r] \times\ZZ.
\]
Thus $f_{y} = f_{\sdeltasv}\in C_{\rm frieze}(A^T)$.
Since $B$ is acyclic, we know by \coref{:admissible-exist} that $\delta^\sv$-admissible functions on $\YB$ exist for every 
$\delta^\sv \in \ABT(\Zmax)$.
By \prref{:trop-frieze-dual-1}, the map in \eqref{eq:f-ysv-1} is surjective.

Suppose that $y \in \calU^{\rm admi}(\YB)$ and $\delta^\sv \in \YBZ$ are such that $f_{y} = f_{\sdeltasv}$. If 
$y$ is $\delta_1^\sv$-admissible for
$\delta_1^\sv \in \ABT(\Zmax)$, then
$f_{{\scriptstyle \delta_1^\sv}} = f_{y} = f_{\sdeltasv}$. So $\delta_1^\sv = \delta^\sv$ by \prref{:trop-frieze-dual-1}. 
Thus $y$ is $\delta^\sv$-admissible.
\end{proof}

We have the following alternative description of $f_y \in C_{\rm frieze}(A^T)$ for  $y \in 
\calU^{\rm admi}(\YB)$.

\ble{:f-y-alter}
For any $y \in 
\calU^{\rm admi}(\YB)$ and any $m \in \ZZ$, $f_y|_{\{m\}}$ is the denominator vector of $y$ at $t(1, m)$. 
\ele

\begin{proof}  
Let $y\in \calU^{\rm admi}_{\delta^\sv}(\YB)$ for $\delta^\sv \in \ABTZ$.
By definition (see \reref{:delta-unique}), the denominator vector of $y$ at $t$ is 
$(\delta^\sv_t)^T$ for every $t \in \TT_r$. \leref{:f-y-alter} now follows from \leref{:inverse-f-deltasv}.
\end{proof}

Recall from \deref{:y-degree} that for every $(i, m) \in [1, r] \times \ZZ$ and any $y \in \FF_{>0}(\YB)$ one has the 
integer $(y(i, m) \, \|\, y)_d$, 
called the {\it $d$-compatibility degree} of $y \in \FF_{>0}(\YB)$ with $y(i, m)$.

\bpr{:trop-frieze-dual-3}
For  $y \in \calU^{\rm admi}(\YB)$ and  $f_y \in \CfriezeAT$ defined in \eqref{eq:f-ysv}, one also has
\begin{equation}\label{eq:f-ysv-3}
f_{y}(i, m) = (y(i, m) \, \|\, y)_d, \hs (i, m) \in [1, r] \times \ZZ.
\end{equation}
\epr

\begin{proof}
By \deref{:y-degree}, one has 
$(y(i, m) \, \|\, y)_d = \val_{d_{y(i, m)}}(y)$, where $d_{y(i, m)} \in \YB(\Zmax)$ is the $d$-tropical point of $y(i, m)$ 
defined in \ldref{:y-degree}. Recall that $\rho(i, m)\in \YB(\Zmax)$ denotes $g$-vector of $x^\sv(i, m) \in \fX(\ABT)$.  Since the two tropical points $\rho(i, m)$ and $d_{y(i, m)}$ in $\YB(\Zmax)$ have the same coordinate vector $-e^i\in\mathbb Z_{\rm row}^r$ at $t(i,m)\in\mathbb T_r$, one has 
$\rho(i, m) = d_{y(i, m)}$. Thus 
\eqref{eq:f-ysv-3} follows from 
\eqref{eq:f-ysv}. 
\end{proof}

The following decomposition property of tropical friezes follows from the definition in \eqref{eq:f-ysv}.

\ble{:trop-frieze-decomposition}
If $y_1, y_2 \in \calU^{\rm admi}(\YB)$ are such that $y_1 y_2 \in \calU^{\rm admi}(\YB)$, then 
$f_{y_1 y_2} = f_{y_1} + f_{y_2}$.
\ele

\subsection{Two realizations of cluster-additive functions}\label{ss:dual-additive}
Recall that $C_{\rm add}(A)$ denotes 
the set of all cluster-additive functions associated to $A$. Parallel to what we did in $\S$\ref{ss:dual-trop-frieze}, we now 
 give two realizations of elements in $C_{\rm add}(A)$ using 
the (same) Fock-Goncharov dual pair
$\OO_{A} = (\AA_{B^T}, \YY_B)$. 
As in \eqref{eq:k-rho-intro}, for 
$\rho = \{(\rho_{t; 1}, \ldots, \rho_{t; r})\}_{t \in \TT_r} \in \YBZ$, define 
\begin{equation}\label{eq:k-rho-i-m}
k_\rho:\;\; [1, r] \times \ZZ \longrightarrow \ZZ, \;\; k_\rho(i, m) = \val_\rho (y(i, m)) = \rho_{t(i, m); i},
\end{equation}
where the last identity is due to the definition of $y(i, m)$ as the $i^{\rm th}$ variable of the cluster ${\bf y}_{t(i, m)}$.

\bpr{:additive-dual-1}
For any symmetrizable generalized Cartan matrix $A$ and with $B = B_A$ as in \eqref{eq:BA-intro}, 
$k_\rho \in C_{\rm add}(A)$ for every $\rho \in \YBZ$, and one has a bijection
\begin{equation}\label{eq:rho-k}
\YY_{B}(\Zmax) \longrightarrow C_{\rm add}(A), \;\;\rho \longmapsto k_\rho.
\end{equation}
\epr

\begin{proof}
By \prref{:xim-yim}, the variables $\{y(i, m)\}$ in $\fX(\YB)$ satisfy \eqref{eq:y(m,i)-relations}.
Applying $\val_\rho$ to both sides of \eqref{eq:y(m,i)-relations}, one sees that
$k_\rho \in C_{\rm add}(A)$ for every $\rho \in \YBZ$. By  \eqref{eq:Y-cluster-to-slice-1}, $(y(1, 0), \ldots, y(r, 0))$ is a set of free generators of $\FF_{>0}(\YY_{B})$. By \leref{:trop-intro}, 
the map in \eqref{eq:rho-k} is bijective. 
\end{proof}

\bre{:additive-first}
{\rm
While the $\bfY$-variables $\{y(i, m)\}_{(i, m)\in [1, r] \times \ZZ}$ were not considered in 
\cite{GHL:friezes}, the facts that every $\rho \in \YBZ$ gives rise to an element in $\CaddA$ via 
$(i, m) \mapsto \rho_{t(i, m); i}$ and that every element in $C_{\rm add}(A)$ is
of this form for a unique $\rho \in \YBZ$ were first proved in  \cite[$\S$5]{GHL:friezes}.
\hfill $\diamond$
}
\ere

For every $\rho \in \YBZ$ and every $m \in \ZZ$, applying $\val_\rho$ to \eqref{eq:Y-cluster-to-slice-1} and 
\eqref{eq:Y-cluster-to-slice-2} respectively, one has
\begin{align}\label{eq:ca-to-rho-coord}
&\rho_{t(1,m);i} = k_\rho(i, m)  + \sum_{j=1}^{i-1} a_{j,i} [k_\rho(j,m)]_+, \qquad i \in [1,r], \\
\label{eq:ca-to-rho-coord-1}
&\rho_{t(1, m+1); i} = -k_\rho(i, m) -\sum_{j=i+1}^r a_{j,i} [k_\rho(j,m)]_+, \qquad i \in [1,r].
\end{align}
In view of \eqref{eq:ca-to-rho-coord} and \eqref{eq:ca-to-rho-coord-1}, we introduce
two bijective  piece-wise linear maps 
\[
E_A^+, \; E_A^-:\;\;  \ZZ^r \longrightarrow \ZZ_{\rm row}^r,
\]
where $E_A^+({\bf d}) = ((E_A^+({\bf d}))_1, \ldots, (E_A^+({\bf d}))_r)$
and $E_A^-({\bf d}) = ((E_A^-({\bf d}))_1, \ldots, (E_A^-({\bf d}))_r)$ for 
${\bf d} = (d_1, \ldots, d_r)^T \in \ZZ^r$ are respectively given by
\[
(E_A^+({\bf d}))_i = d_i + \sum_{j=1}^{i-1} a_{j, i} [d_j]_+ \hs \mbox{and} \hs 
(E_A^-({\bf d}))_i = d_i + \sum_{j=i+1}^{r} a_{j, i} [d_j]_+, \hs i \in [1, r],
\]
or, in a more compact form, 
\begin{equation}\label{eq:EA-plus}
 E_A^+({\bf d}) = {\bf d}^T + [{\bf d}]_+^T U_A \hs \mbox{and} \hs
E_A^-({\bf d}) = {\bf d}^T + [{\bf d}]_+^T L_A.
\end{equation}
where 
\begin{equation}\label{eq:UA}
U_A = \begin{pmatrix}
        0 & a_{1,2} & \cdots &  a_{1,r}\\
        0 & 0 & \cdots & a_{2,r} \\
        \vdots & \vdots & \ddots & \vdots \\
        0 & 0 & \cdots & 0
    \end{pmatrix} \hs \mbox{and} \hs 
L_A= \begin{pmatrix}
        0 & 0 & \cdots & 0 \\
        a_{21} & 0 & \cdots & 0\\
        \vdots & \vdots & \ddots & \vdots \\
        a_{r,1} & a_{r,2} & \cdots & 0
    \end{pmatrix}.
\end{equation}

\bre{:k-additive-E}
{\rm
Note that the condition \eqref{eq:cl-additive} for a map $k: [1, r]\times \ZZ \to \ZZ$ to be in $\CaddA$ can be rewritten as
$E_A^- (k|_{\{m\}}) + E_A^+ (k|_{\{m+1\}}) = 0$ for every $m \in \ZZ$, or, equivalently,
\begin{equation}\label{eq:k-additive-m}
k|_{\{m+1\}} = (E_A^+)^{-1} (-E_A^-(k|_{\{m\}})), \hs \forall \; m \in \ZZ,
\end{equation}
where $k|_{\{m\}} \in \ZZ^r$ is the $m^{\rm th}$ slice of $k$ (see \deref{:slice}).
}
\ere

Rewriting \eqref{eq:ca-to-rho-coord} and \eqref{eq:ca-to-rho-coord-1}, 
we obtain the following descriptions of the inverse of the bijection 
in \eqref{eq:rho-k} (compare with the corresponding statements for tropical friezes in \leref{:inverse-f-deltasv}).

\bpr{:rho-chi-A}
For any symmetrizable generalized Cartan matrix $A$ and with $B = B_A$ as in \eqref{eq:BA-intro}, and  for any $\rho \in \YBZ$ and any
$m \in \ZZ$, one has
\begin{equation}\label{eq:rho-EA}
\rho_{t(1, m)} = E_A^+(k_\rho|_{\{m\}}) \hs \mbox{and} \hs
\rho_{t(1, m+1)} = -E_A^-(k_\rho|_{\{m\}}).
\end{equation}
Consequently,  one has
\begin{equation}\label{eq:rho-chi-A}
\rho_{t(1, m+1)} = -E_A^- (E_A^+)^{-1} (\rho_{t(1, m)}).
\end{equation}
\epr

For the second realization of elements in $\CaddA$, for each $(i, m) \in [1, r] \times \ZZ$,
recall from \leref{:B-im-column} 
that $y(i, m)$ is a global variable of $\YB$, and let $\delta^{\sv}(i, m) \in \ABTZ$ be the $g$-vector of $y(i, m)$. 
For $x^\sv \in \calU^{\rm admi}(\ABT)$, define
\begin{equation}\label{eq:k-xsv}
k_{x^\sv}: \;\; [1, r] \times \ZZ \longrightarrow \ZZ, \;\; k_{x^\sv}(i, m) = \val_{\delta^\sv(i, m)}(x^\sv).
\end{equation}

\bpr{:additive-dual-2}
For any symmetrizable generalized Cartan matrix $A$ and with $B = B_A$ as in \eqref{eq:BA-intro}, 
$k_{x^\sv} \in C_{\rm add}(A)$ for every $x^\sv \in \calU^{\rm admi}(\ABT)$,  and one has a surjective map
\begin{equation}\label{eq:k-xsv-2}
\calU^{\rm admi}(\ABT) \longrightarrow C_{\rm add}(A), \;\; x^\sv \longmapsto k_{x^\sv}.
\end{equation}
Moreover, for $x^\sv \in \calU^{\rm admi}(\ABT)$ and $\rho \in \YBZ$, one has $k_{x^\sv} = k_\rho$
 if and only if $x^\sv$ is $\rho$-admissible. 
\epr

\begin{proof} The proof is similar to that of \prref{:trop-frieze-dual-2}: 
suppose that $x^\sv \in \calU^{\rm admi}(\ABT)$ is $\rho$-admissible for 
$\rho \in \YBZ$. By \prref{:duality-pairing}, for every $(i, m) \in [1, r] \times\ZZ$, one has
\[
k_{x^\sv} (i, m) = \val_{\delta^\sv(i, m)}(x^\sv) = \val_\rho(y(i, m)) = k_\rho(i, m).
\]
Thus $k_{x^\sv} = k_\rho\in C_{\rm add}(A)$.
Since $B$ is acyclic, by \coref{:admissible-exist} we know that $\rho$-admissible functions on $\ABT$ exist for every $\rho \in \YBZ$.
By \prref{:additive-dual-1}, the map in \eqref{eq:k-xsv-2} is surjective.

Suppose that $x^\sv \in \calU^{\rm admi}(\ABT)$ and $\rho \in \YBZ$ are such that $k_{x^\sv} = k_\rho$. If 
$x^\sv$ is $\rho_1$-admissible for
$\rho_1\in \YBZ$, then
$k_{\rho_1} = k_{x^\sv} = k_{\rho}$. So $\rho_1 = \rho$ by \prref{:additive-dual-1}. 
Thus $x^\sv$ is $\rho$-admissible.
\end{proof}

By \deref{:x-degree}, the $d$-compatibility degree of any $x^\sv \in \FF_{>0}(\ABT)$ with $x^\sv(i, m)$ is defined as
\begin{equation}\label{eq:d-xsvim}
(x^\sv(i, m) \, \|\, x^\sv)_d = \val_{d_{x^\sv(i, m)}}(x^\sv),
\end{equation}
where $d_{x^\sv(i, m)} \in \ABT(\Zmax)$ is the $d$-tropical point  of $x^\sv(i, m)$ defined in \ldref{:x-degree}.

\bpr{:additive-dual-3}
For any $x^\sv \in \calU^{\rm admi}(\ABT)$ and $k_{x^\sv} \in \CaddA$ defined as in \eqref{eq:k-xsv}, one has
\begin{equation}\label{eq:k-xsv-3}
k_{x^\sv}(i, m) = (x^\sv(i, m) \, \|\, x^\sv)_d, \hs (i, m) \in [1, r] \times \ZZ.
\end{equation}
\epr

\begin{proof}
By \leref{:d-g}, one has $\delta^{\sv}(i, m) = d_{x^\sv(i, m)} \in \ABT(\Zmax)$, so
\eqref{eq:k-xsv-3} follows \eqref{eq:k-xsv} and  \eqref{eq:d-xsvim}. 
\end{proof}

For $x^\sv \in \calU^{\rm admi}(\ABT)$ and $t \in \TT_r$,  let 
${\bf d}^{x^\sv}_t \in \ZZ^r$
be the denominator vector of $x^\sv$ at $t$.
We have the following alternative description of $k_{x^\sv} \in \CaddA$. Compare with \leref{:f-y-alter}.

\bco{:kxsv-slice}
For any $x^\sv \in \calU^{\rm admi}(\ABT)$ and any $m\in\mathbb Z$, one has $k_{x^\sv}|_{\{m\}} = {\bf d}^{x^\sv}_{t(1, m)}$. 
\eco

\begin{proof}
The statement follows from \eqref{eq:k-xsv-3} and the fact that 
${\bf x}^\sv_{t(1, m)} = (x^\sv(1, m), \ldots, x^\sv(r, m))$.
\end{proof}
Consider the well-defined map (see \reref{:admi-unique}) 
\[
\calU^{\rm admi}(\ABT) \longrightarrow \YBZ, \;\; x^\sv \longmapsto \rho(x^\sv),
\]
where for $x^\sv \in \calU^{\rm admi}(\ABT)$, $\rho(x^\sv)$ is the unique tropical point of $\YB$ such that
$x^\sv$ is $\rho(x^\sv)$-admissible (see \coref{:acyclic-rho-unique}).
By \ldref{:globals-A-are-admissible}, $\rho(x^\sv)$ is the $g$-vector of $x^\sv$ if $x^\sv$ is a cluster monomial.

\bco{:rho-xsv-d} For any symmetrizable generalized Cartan matrix $A$ and with $B = B_A$ as in \eqref{eq:BA-intro}, and
for any $x^\sv \in \calU^{\rm admi}(\ABT)$ and $m \in \ZZ$, one has
\begin{equation}\label{eq:rho-EE}
\rho(x^\sv)_{t(1, m)} =E_A^+ \left({\bf d}^{x^\sv}_{t(1, m)}\right) = 
-E_A^{-}\left({\bf d}^{x^\sv}_{t(1, m-1)}\right).
\end{equation}
\eco

\begin{proof}
\eqref{eq:rho-EE} follows directly from \eqref{eq:rho-EA} and \coref{:kxsv-slice}.
\end{proof}

\bre{:g-d}
{\rm
When $x^\sv$ is a cluster monomial, the relation in \eqref{eq:rho-EE} between the $g$-vector and denominator vector of $x^\sv$ at 
an acyclic seed can be found in \cite[Proposition 8]{rupel-stella:consequences}.
}
\ere

The following decomposition property of cluster-additive functions follows from the definition in \eqref{eq:k-xsv}.

\ble{:ca-decomposition}
If $x_1^\sv, x_2^\sv \in \calU^{\rm admi}(\AA_{B^T})$ are such that $x_1^\sv x_2^\sv \in \calU^{\rm admi}(\AA_{B^T})$, then 
$k_{x_1^\sv x_2^\sv} = k_{x_1^\sv} + k_{x_2^\sv}$.
\ele
\subsection{The shift on cluster-additive functions}\label{ss:shift}
Fix an $r \times r$ symmetrizable generalized Cartan matrix $A$.
For any $k: [1, r] \times \ZZ \to \ZZ$, define 
\[
k_{[1]}: \;\; [1, r] \times \ZZ \longrightarrow  \ZZ,\;\; k_{[1]}(i, m) = k(i, m-1).
\]
It is then clear from 
\eqref{eq:cl-additive} that $k \in \CaddA$ if and only if $k_{[1]} \in \CaddA$. 

 We now determine the maps
on $\YBZ$ and $\calU^{\rm admi}(\ABT)$ that correspond to the shift
\[
\CaddA \longrightarrow \CaddA, \;\; k \longmapsto k_{[1]}
\]
under the two maps 
respectively given in \prref{:additive-dual-1} and 
\prref{:additive-dual-2}.

Let $B = B_A$ as in \eqref{eq:BA-intro} and 
$\{B_t\}_{t \in \TT_r}$ the matrix pattern with $B_{t_0} = B$. Denote by 
\[
t \longmapsto t[1], \hs t \in \TT_r,
\]
the unique tree automorphism of $\TT_r$ such that $t_0[1] = t(1, 1)$, where $t(1,1) \in \TT_r$ is given in \eqref{eq:right}. 
Note that $B_{t_0[1]}=B_{t_0} =B$ by \leref{:B-im-column}, and one  
checks directly that the extended mutation matrix of $\Aprin$ at $t_0[1]$ is 
$(B \bb -I_r)$.
As $B_t = B_{t[1]}$ for all $t \in \TT_r$, 
for each $\rho = \{\rho_t \in \ZZ^r_{\rm row}\}_{t \in \TT_r} \in \YBZ$, one then has 
\[
\rho[1] 
=\{(\rho[1])_t \in \ZZ^r_{\rm row}\}_{t \in \TT_r} \in \YY_{B}(\Zmax)
\]
such that
$(\rho[1])_{t[1]} = \rho_t$ for every $t \in \TT_r$. We define the {\it shift map on $\YBZ$} to be the bijection
\begin{equation}\label{eq:shift-1}
\YY_{B}(\Zmax) \longrightarrow \YY_{B}(\Zmax), \;\; \rho \longmapsto \rho[1].
\end{equation}
Note that under the tree automorphism $t\mapsto t[1]$, one has $t(i,m)[1]=t(i,m+1)$ for any $(i,m)\in[1,r]\times \mathbb Z$.
Thus for every $\rho \in \YY_{B}(\Zmax)$, one has
\[ (\rho[1])_{t(1, m)}= (\rho[1])_{t(1, m-1)[1]}=\rho_{t(1,m-1)}.\]
On the other hand, by \prref{:rho-chi-A}, one has $\rho_{t(1, m)} = -E_A^- (E_A^+)^{-1} (\rho_{t(1, m-1)})$. Hence, one has
\begin{equation}\label{eq:rhosv-shift-1}
(\rho[1])_{t(1, m)}=\rho_{t(1,m-1)}= E_{A}^+ (E_{A}^-)^{-1} (-\rho_{t(1, m)}), \hs \forall \; m \in \ZZ.
\end{equation}

\ble{:k-rho-shift-1}
For any $\rho \in \YBZ$, one has $k_{\rho[1]} = (k_\rho)_{[1]} \in \CaddA$.
\ele

\begin{proof}
It suffices to show $  k_{\rho[1]}|_{\{m\}}= (k_\rho)_{[1]}|_{\{m\}}$, which can be proved as follows:
\begin{eqnarray}
   k_{\rho[1]}|_{\{m\}}  &\xlongequal{{\rm  by}\;\;  \eqref{eq:rho-EA}}& (E_A^+)^{-1} \left((\rho[1])_{t(1, m)}\right)
   \xlongequal{{\rm  by}\;\;  \eqref{eq:rhosv-shift-1}}  (E_{A}^-)^{-1} (-\rho_{t(1, m)})\nonumber\\
   &\xlongequal{{\rm  by}\;\;  \eqref{eq:rho-EA}}&k_\rho|_{\{m-1\}}
    \xlongequal{{\rm  by \;\;def.}}
    (k_\rho)_{[1]}|_{\{m\}}.\nonumber
\end{eqnarray}
\end{proof}

Recall that the clusters of $\AA_{B^T}$ are denoted as $\{\bfx_t^\sv\}_{t \in \TT_r}$. Again due to the fact that
$B_{t} = B_{t[1]}$ 
for every $t \in \TT_r$, we have the well-defined algebra automorphism $\tau$ of $\calU(\ABT)$ satisfying
\[
\tau:\;\; \calU(\ABT) \longrightarrow \calU(\ABT), \;\; {\bf x}^\sv_t \longmapsto {\bf x}^\sv_{t[1]}, \hs \forall \; t \in \TT_r.
\]
For $x^\sv \in \calU(\ABT)$, we also write $\tau(x^\sv)$ as $x^\sv[1]$. Recall again that for $x^\sv \in \calU^{\rm admi}(\ABT)$, $\rho(x^\sv)$ denotes the unique point in $\YBZ$  such that $x^\sv$ is 
$\rho(x^\sv)$-admissible. 

\ble{:tau}
An element $x^\sv \in \calU(\ABT)$ is admissible if and only if 
$x^\sv[1]$ is admissible, and in this case,
\begin{equation}\label{eq:rho-kxsv-shift}
\rho(x^\sv[1]) = \rho(x^\sv)[1] \hs \mbox{and} \hs k_{x^\sv[1]} = (k_{x^\sv})_{[1]} \in \CaddA.
\end{equation}
\ele

\begin{proof}
Let $x^\sv \in \calU^{\rm admi}(\ABT)$ and let $\rho = \rho(x^\sv) = \{\rho_t\}_{t \in \TT_r}$.
For every $t \in \TT_r$, one has
\[
x^\sv = ({\bf x}_t^\sv)^{-\rho_t^T} F\left(({\bf x}^\sv_t)^{B_t^T}\right)
\]
for some $F \in \ZZ_{\geq 0}[y_1, \ldots, y_r]$ with constant term $1$.
Applying $\tau$, one has
\[
\tau(x^\sv) = ({\bf x}^\sv_{t[1]})^{-\rho_t^T} F\left(({\bf x}^\sv_{t[1]})^{B^T_t}\right) =
({\bf x}^\sv_{t[1]})^{-(\rho[1]_{t[1]})^T} F\left(({\bf x}^\sv_{t[1]})^{B^T_{t[1]}}\right).
\]
It follows that $\tau(x^\sv) \in \calU^{\rm admi}_{\rho[1]}(\ABT)$. The second identity in \eqref{eq:rho-kxsv-shift} follows from 
\prref{:additive-dual-2}. Applying the same arguments to $\tau^{-1}$, one sees that 
$x^\sv \in \calU^{\rm admi}(\ABT)$ if and only if
$x^\sv[1] \in \calU^{\rm admi}(\ABT)$.
\end{proof}

\subsection{The ensemble map from $C_{\rm frieze}(A)$ to $C_{\rm add}(A)$}
Let again $B = B_A$ as in \eqref{eq:BA-intro}. Since 
$(B_{A^T})^T = -B$, replacing $A$ by $A^T$ in \prref{:trop-frieze-dual-1}, we have the bijection 
\[
\AA_{-B}(\Zmax) \longrightarrow C_{\rm frieze}(A), \;\; \delta \longmapsto f_\delta.
\]
Denote the generic frieze pattern of $\AA_{-B}$ by
\begin{equation}\label{eq:xim-new}
\fX^\flat(\AA_{-B}) = \{x(i, m)\}_{(i, m) \in [1, r] \times \ZZ}.
\end{equation}
Then by \prref{:xim-yim}, it satisfies \eqref{eq:x(i,m)s}, and one has 
\begin{equation}\label{eq:f-delta}
f_\delta(i, m) = \val_\delta(x(i, m)), \hs (i, m) \in [1, r] \times \ZZ,
\end{equation}
for $\delta \in 
\AA_{-B}(\Zmax)$.  
As the mutation rules for the coordinate vectors
of elements in $\AA_B(\Zmax)$ are the same as that for elements in $\AA_{-B}(\Zmax)$, by identifying $\AA_B(\Zmax)$ with
$\AA_{-B}(\Zmax)$, one gets a map
\[
p_\sA: C_{\rm frieze}(A) \longrightarrow C_{\rm add}(A), \quad f_\delta \longmapsto k_{\pB(\delta)},
\]
where $\pB: \AA_B(\Zmax) \to \YY_B(\Zmax)$ is defined in \eqref{eq:p-map}. 
One checks directly that 
\[
    k_{\pB(\delta)}(i,m) = \sum_{j=i+1}^r (-a_{j,i}) f_\delta(j,m) + \sum_{j=1}^{i-1} (-a_{j,i}) f_\delta(j,m+1),
    \hs (i, m) \in [1, r] \times \ZZ.
\]
The map $p_\sA$ is studied in greater generality by the second author in \cite[\S 2]{Antoine:Y-frieze-patterns}.

\subsection{Cluster-hammock functions}\label{ss:hammock} 
Let $(i,m) \in [1,r]\times \ZZ$. For a Cartan matrix $A$ of finite and simply-laced type, Ringel defined in 
\cite{Ringel:additive} the
{\it cluster-hammock function} $h_{(i,m)}$ as the unique cluster-additive function associated to $A$ such that
for $j \in [1, r]$,
\begin{equation}\label{eq:him-00}
h_{(i,m)}(j,m) = \begin{cases} -1, & \hs j = i, \\ 0, & \hs j \neq i\end{cases}.
\end{equation}
We define $h_{(i, m)} \in \CaddA$ via \eqref{eq:him-00} for an arbitrary symmetrizable generalized Cartan matrix $A$. In light of \prref{:additive-dual-1} and \prref{:additive-dual-2}, we have the following description of $h_{(i, m)}$.

\bpr{:hammock-1}
For any $(i,m) \in [1,r] \times \ZZ$, one has 
$h_{(i,m)}   = k_{\rho(i,m)} = k_{x^\sv(i, m)}$. Consequently, for 
$(j, n) \in [1,r] \times \ZZ$, one has 
\begin{align}\label{eq:him-0}
h_{(i,m)}(j,n) = (x^{\sv}(j, n)\, \| \, x^{\sv}(i, m))_d = (y(i,m) \, \| \, y(j, n))_d 
= \begin{cases}
        -1, & \quad \text{ if }\; x^\sv(j,n) = x^\sv(i,m), \\
        \geq 0, & \quad \text{ otherwise}. 
    \end{cases}
\end{align}
\epr

\begin{proof}
By \leref{:d}, the two cluster-additive functions  $h_{(i, m)}$ and $k_{x^\sv(i, m)}$ agree on $[1, r] \times \{m\}$, so they are equal.
It follows from \eqref{eq:k-xsv-3} and \leref{:d} that for every $(j,n) \in [1, r] \times \ZZ$
\[
h_{(i,m)}(j,n) = (x^{\sv}(j, n)\, \| \, x^{\sv}(i, m))_d = 
\begin{cases}
        -1, & \quad \text{ if }\; x^\sv(j,n) = x^\sv(i,m), \\
        \geq 0, & \quad \text{ otherwise}.
\end{cases}
\]
By \prref{:additive-dual-2}, $h_{(i, m)} = k_{x^\sv(i, m)} = k_{\rho(i, m)}$. By
\leref{:d-g}, $\rho(i, m) = d_{y(i, m)} \in \YBZ$,
where $d_{y(i, m)}$ is the $d$-tropical point of $y(i, m)\in \fX(\YB)$. Thus  for every $(j, n) \in [1, r] \times \ZZ$,
\[
h_{(i, m)}(j, n) = \val_{\rho(i, m)}(y(j, n)) = \val_{d_{y(i, m)}}(y(j, n)) = (y(i, m) \, \|\, y(j, n))_d.
\]
\end{proof}
We now show that cluster-hammock functions are also tropical friezes.

\bpr{:hammock-2}
One has $h_{(i, m)} \in C_{\rm frieze}(A)$ for every $(i, m) \in [1, r] \times \ZZ$.
\epr

\begin{proof}
Let $(i, m) \in [1, r] \times \ZZ$ and write $h = h_{(i,m)}$ for simplicity.
Fix $(j, n) \in [1, r] \times \ZZ$.
Since $h \in C_{\rm add}(A)$, 
\begin{equation}\label{eq:caf-equation-proof-ham-d-comp}
    h(j,n) + h(j,n+1) = \sum_{l=j+1}^r (-a_{l,j})[h(l,n)]_+ + \sum_{l=1}^{j-1} (-a_{l,j})[h(l,n+1)]_+.
\end{equation}
Recall from \reref{:xim-yim} that 
\[
{\bf x}_{t(j,n)}^\sv=(x^\sv(1,n+1),\,\ldots,\,x^\sv(j-1,n+1),\,x^\sv(j,n),\, x^\sv(j+1,n),\, \ldots,\,x^\sv(r,n)).
\]

{\bf Case 1:} $x^\sv(i,m)$ belongs to the cluster ${\bf x}^\sv_{t(j,n)}$ of $\AA_{B^T}$. 
By  \eqref{eq:him-0} and \leref{:d},  one has $h(l,n) \in \{0, -1\}$ for $l \in [j+1, r]$ and 
$h(l,n+1) \in \{0, -1\}$ for $l \in [1, j-1]$. In particular, the right-hand side of \eqref{eq:caf-equation-proof-ham-d-comp} is 
equal to $0$ which is also equal to 
\[
    \left[\sum_{l=j+1}^r (-a_{l,j})h(l,n) + \sum_{l =1}^{j-1} (-a_{l,j})h(l,n+1)\right]_+.
\]

{\bf Case 2:} $x^\sv(i,m)$ does not belong to the cluster ${\bf x}_{t(j,n)}^\sv$.
Again by \eqref{eq:him-0} and \leref{:d}, one has $h(l,n) \geq 0$ for $l \in [j+1, r]$ and $h(l,n+1) \geq 0$ for $l \in [1, j-1]$. In particular, the right-hand side of \eqref{eq:caf-equation-proof-ham-d-comp} is again equal to
\[
    \left[\sum_{l=j+1}^r (-a_{l,j})h(l,n) + \sum_{l =1}^{j-1} (-a_{l,j})h(l,n+1)\right]_+.
\]
Thus $h = h_{(i, m)} \in C_{\rm frieze}(A)$.   
\end{proof}

Recall from \eqref{eq:xim-new} the variables $x(i,m)$ are in 
the generic frieze pattern of $\AA_{-B}$. 
Let 
\[
\fX^\flat(\YY_{-B^T}) = \{y^\sv(i,m)\}_{(i, m) \in [1, r] \times \ZZ}
\]
be the
generic frieze pattern of $\YY_{-B^T}$, and let
$\delta(i, m) \in \AA_{-B}(\Zmax)$ be the $g$-vector of $y^\sv(i, m)$.  By 
\prref{:trop-frieze-dual-3}, one  has $f_{y^\sv(i, m)} \in C_{\rm frieze}(A)$ and 
\[
f_{y^\sv(i, m)}(j, n) = (y^{\sv}(j, n)\, \| \, y^{\sv}(i, m))_d, \hs(j, n) \in [1, r]\times \ZZ.
\]

\bco{:hammock-3}
For each $(i,m) \in [1,r] \times \ZZ$, one has 
$h_{(i,m)} = f_{\delta(i,m)} = f_{y^\sv(i, m)} \in  C_{\rm frieze}(A)$. Consequently,  for $(j, n) \in [1,r] \times \ZZ$, one  has
\begin{equation}\label{eq:him-2}
h_{(i,m)}(j,n)  = (x(i,m) \, \| \, x(j, n))_d = (y^{\sv}(j, n)\, \| \, y^{\sv}(i, m))_d.
\end{equation}
\eco

\begin{proof}
By \prref{:trop-frieze-dual-2}, $f_{\delta(i,m)} = f_{y^\sv(i, m)}$. 
By \leref{:d-g}, $\delta(i,m) = d_{x(i,m)}$,
the $d$-tropical point of $x(i, m) \in \fX(\AA_{-B})$ defined in \ldref{:x-degree}. Thus
for $(j, n) \in [1, r] \times \ZZ$, one has 
\[
f_{\delta(i, m)}(j, n) = \val_{\delta(i, m)}(x(j, n)) = (x(i, m) \, \| \, x(j, n))_d = 
(y^{\sv}(j, n)\, \| \, y^{\sv}(i, m))_d,
\]
where in the last identity we used the definition of $f_{y^\sv(i, m)}$.  
Now the two tropical friezes $h_{(i,m)}$ and $f_{\delta(i,m)}$ agree on $[1, r] \times \{m\}$ by \leref{:d}, so they are equal,
and \eqref{eq:him-2} holds. 
\end{proof}

\bre{:guo}
When $A$ is of finite and simply-laced type, the fact that every cluster-hammock function is also a tropical frieze
has been proved by Guo \cite[\S 3.4]{Guo:tropical-frieze} using categorical arguments.
\hfill $\diamond$
\ere

\subsection{A duality of $d$-compatibility degrees}
By \cite[Lemma 4.13]{cao-gyoda:bongartz}, for any mutation matrix $B$ the assignment
$x_{t;i} \mapsto x^\sv_{t;i}$, for  $(t , i ) \in  \TT_r \times [1,r]$, gives a 
well-defined bijection 
\begin{equation}\label{eq:x-xvee}
\mathfrak{X}(\AA_{-B}) \longrightarrow \mathfrak{X}(\AA_{B^T}), \;\; x \longmapsto x^\sv.
\end{equation}
Following \cite{Reading-Stella_2018}, a mutation matrix $B$ is said to have {\it Property D} if
\begin{equation}\label{eq:xzxz}
(x \, \|\,   z)_d = (z^\sv \, \|\,   x^\sv)_d, \hs \forall \; x, z \in \mathfrak{X}(\AA_{-B}).
\end{equation}
It is shown in \cite[Theorem 2.2 and Theorem 2.3]{Reading-Stella_2018} that $B$ has Property $D$ if it is of rank $2$ or of finite
type. An example where \eqref{eq:xzxz} does not hold is, however, given in \cite[Example 4.21]{Fu-Gyoda}. 
As $\mathfrak{X}^\flat(\AA_{-B}) = \mathfrak{X}(\AA_{-B})$ and $\mathfrak{X}^\flat(\AA_{B^T}) = \mathfrak{X}(\AA_{B})$ (see \eqref{eq:frieze-patrn-nota} for the notation) when $B$ has rank $2$, or when $B$ is of finite type
\cite[Theorem 3.1]{keller-knitting-algo},   we have the following generalization of a result of N. Reading and S. Stella in \cite[Theorem 2.2 and Theorem 2.3]{Reading-Stella_2018}.

\bpr{:d-dual}
For any $B = B_A$ as in \eqref{eq:BA-intro} and for any  $x,z \in \mathfrak{X}^\flat(\AA_{B})$, one has
\[
(x \, \|\,    z)_d=(z^\sv \, \|\,  x^\sv)_d.
\]
\epr

\begin{proof}
The statement follows by combining \eqref{eq:him-0} and 
\eqref{eq:him-2}.
\end{proof}

\subsection{An auxiliary lemma}\label{ss:auxi}  In this section we prove \leref{:z-x-g-vector} which will be used in 
$\S$\ref{ss:thmE}.

Let again $A = (a_{i, j})$ be any
$r \times r$ symmetrizable generalized Cartan matrix, and let $B = B_A$ as in \eqref{eq:BA-intro}.
Consider the two positive spaces 
\[
\AA_{B} = (\FF(\AA_B), \{{\bf z}_t\}_{t \in \TT_r}) \quad\text{and}\quad \AA_{-B} = (\FF(\AA_{-B}), \{{\bf x}_t\}_{t \in \TT_r}).
\]
As the cluster variables of $\AA_B$ and of $\AA_{-B}$ obey the same mutation rules, one has the well-defined
algebra isomorphism $\varphi: \calU(\AA_B) \rightarrow \calU(\AA_{-B})$ satisfying
\[
\varphi(z_{t; i}) = x_{t; i}, \hs t \in \TT_r, \; i \in [1, r],
\]
and $\varphi$ restricts to a bijection $\varphi: \calU^{\rm mono}(\AA_B) \rightarrow \calU^{\rm mono}(\AA_{-B})$.
For $z \in \calU^{\rm mono}(\AA_B)$, we now determine the $g$-vector of $\varphi(z)$ in terms of the $g$-vector of $z$.
Note first the well-defined bijection
\[
\YY_{B^T}(\Zmax) \longrightarrow \YY_{-B^T}(\Zmax), \;\; \rho^\sw =\{(\rho^\sw)_t
\in \ZZ^r_{\rm row}\}_{t \in \TT_r}\longmapsto -\rho^\sw
:=\{-(\rho^\sw)_t\in \ZZ^r_{\rm row}\}_{t \in \TT_r}.
\]
Consequently, every $\rho^\sw \in \YY_{B^T}(\Zmax)$ gives rise to $k_{-\rho^\sw} \in C_{\rm add}(A^T)$.

For $z \in \calU^{\rm admi}(\AA_B)$, let $\rho^\sw(z)$ be the unique point 
in $\YY_{B^T}(\Zmax)$ such that $z$ is $\rho^\sw(z)$-admissible.

\ble{:z-denom}
For any $z \in \calU^{\rm admi}(\AA_B)$ and any $m \in \ZZ$, one has 
\[
k_{-\rho^\sw(z)}|_{\{m-1\}} ={\bf d}^z_{t(1, m)},
\]
where ${\bf d}^z_{t(1, m)} \in \ZZ^r$ is the denominator vector of 
$z$ at $t(1, m) \in \TT_r$.
\ele

\begin{proof}
Fix $m \in \ZZ$ and consider the path
\begin{equation}\label{eq:tt-mm}
\begin{xy}
(0,1)*+{t(1, m-1)}="A2",(25,1)*+{t(2, m-1)}="A3", (46,1)*+{\;\;\cdots\;\;}="A4", (67,1)*+{t(r, m-1)}="A5",(89,1)*+{t(1, m)}="A6",
\ar@{-}^{\;\;1\;\;}"A2";"A3", \ar@{-}^{\;\;\;\;\;\;2}"A3";"A4", 
\ar@{-}^{r-1\;\;\;\;\;\;}"A4";"A5",\ar@{-}^{\;\;\;r}"A5";"A6" 
\end{xy}
\end{equation}
in $\TT_r$. For $i \in [1, r]$, let $t'(i, m-1) \in \TT_r$ be such that
\begin{xy}(0,1)*+{t(i, m-1)}="A",(25,1)*+{t'(i, m-1)}="B",\ar@{-}^i"A";"B" \end{xy}, i.e., 
\[
t'(i, m-1) = \begin{cases} t(i+1, \, m-1), & i \in [1, r-1], \\
t(1, m), & i = r.
\end{cases}
\]
Let $\{B_t\}_{t \in \TT_r}$ be the matrix pattern determined by $B_{t_0} = B$, and 
let $\bfy^\sw_t = (y^\sw_{t; 1}, \,\ldots, \,y^\sw_{t; r})$ be the cluster of $\YY_{B^T}$.
Following the sequence of mutations of the seeds of $\YY_{B^T}$ from $t(1, m)$ to $t(1, m-1)$ along the path in \eqref{eq:tt-mm},
one checks directly that for every $i \in [1, r]$, the $i^{\rm th}$ column of $B_{t'(i, m-1)}^T$ is non-negative, and
\[
(y^\sw_{t(i, m-1); i})^{-1} = y^\sw_{t'(i, m-1); i} \hs \mbox{and} \hs z_{t(1, m); i} = z_{t'(i, m-1); i}.
\]
By \leref{:Y-global-mono}, 
$y^\sw_{t'(i, m-1); i}$
is a global variable of $\YY_{B^T}$  for every $i \in [1, r]$.
Note that the $g$-vector 
of $(y^\sw_{t(i, m-1); i})^{-1} = y^\sw_{t'(i, m-1); i}$
coincides with the
$d$-tropical point $d_{z_{t(1, m); i}}$ of $z_{t(1, m); i} = z_{t'(i, m-1);i}$ (see \ldref{:x-degree}), as they are both in $\AA_{B}(\Zmax)$ and
have the same coordinate vector $-e^i$ at the vertex $t'(i,m-1)$.
Writing 
${\bf d}^z_{t(1, m)} = ({\bf d}^z_{t(1, m); 1}, \ldots, {\bf d}^z_{t(1, m); r})^T$ and by
\prref{:duality-pairing},  one has
\[
{\bf d}^z_{t(1, m); i} = \val_{d_{z_{t(1, m); i}}} (z) =\val_{\rho^\sw(z)}( (y^\sw_{t(i, m-1); i})^{-1})
= -(\rho^\sw(z))_{t(i, m-1); i}
=k_{-\rho^\sw(z)}(i, m-1),
\]
for  $i \in [1, r]$.
Thus $k_{-\rho^\sw(z)}|_{\{m-1\}} ={\bf d}^z_{t(1, m)}$.
\end{proof}

Recall from $\S$\ref{ss:shift} the  shift map $\YY_{-B^T}(\Zmax) \to \YY_{-B^T}(\Zmax),
\rho^\sv \mapsto \rho^\sv[1]$.

\ble{:z-x-g-vector}
If $z \in \calU^{\rm mono}(\AA_B)$ has $g$-vector $\rho^\sw(z) \in \YY_{B^T}(\Zmax)$,  then 
$\varphi(z) \in \calU^{\rm mono}(\AA_{-B})$ has $g$-vector
\[
\rho^\sv(\varphi(z)) = (-\rho^\sw(z))[1]  \in \YY_{-B^T}(\Zmax).
\]
\ele

\begin{proof}
Let $z \in \calU^{\rm mono}(\AA_B)$ and denote by $\rho^\sv(x)\in \YY_{-B^T}(\Zmax)$ the $g$-vector of $x:=\varphi(z)\in \calU^{\rm mono}(\AA_{-B})$. By the definition of $\varphi$, one can see that ${\bf d}^x_t = {\bf d}^z_t$ for every $t \in \TT_r$, where
${\bf d}^x_t$ and  ${\bf d}^z_t$ are the respective  denominator vectors of $z$ and $x$ at $t$.

Fix an integer $m\in\mathbb Z$.
By 
\coref{:rho-xsv-d} and \leref{:z-denom},
\[
\rho^\sv(x)_{t(1, m)} = E_{A^T}^+({\bf d}^x_{t(1, m)}) = E_{A^T}^+({\bf d}^z_{t(1, m)}) = E_{A^T}^+(k_{-\rho^\sw(z)}|_{\{m-1\}}).
\]
By  \prref{:rho-chi-A}, one has $(-\rho^\sw(z))_{t(1, m-1)}=E_{A^T}^+(k_{-\rho^\sw(z)}|_{\{m-1\}})$. Thus
\[
\rho^\sv(x)_{t(1, m)} =E_{A^T}^+(k_{-\rho^\sw(z)}|_{\{m-1\}})=(-\rho^\sw(z))_{t(1, m-1)}= ((-\rho^\sw(z))[1])_{t(1, m)}.
\]
As the two tropical points $\rho^\sv(x)$ and $(-\rho^\sw(z))[1]$ in $\YY_{-B^T}(\Zmax)$ have the same coordinate vector at the vertex $t(1,m)$, one has $\rho^\sv(\varphi(z)) =\rho^\sv(x) = (-\rho^\sw(z))[1]$. 
\end{proof}


\section{Finite type: Ringel's conjecture and applications}\label{s:finite} 
We now prove the
results in Theorem D and Theorem E in $\S$\ref{ss:finite-intro} for Cartan matrices of finite type.

\subsection{Global $\bfY$-variables in finite type}\label{ss:finite-patterns} 
Recall that for any symmetrizable generalized Cartan matrix $A$ and with $B = B_A$ as in 
\eqref{eq:BA-intro} we have set 
\[
\fX^{\rm global}(\YB) = \mathfrak{X}(\YY_{B})\cap \calU(\YY_{B}),
\]
and that we have an injective map 
$\lambda: \mathfrak{X}^{\rm global}(\YY_{B}) \to \fX(\ABT)$ by \coref{:injective}.
Recall also 
\[
\fX^\flat(\ABT) = \{x^\sv(i, m)\}_{(i, m) \in [1, r]\times \ZZ} \hs \mbox{and} \hs
\fX^\flat(\YB) = \{y(i, m)\}_{(i, m) \in [1, r]\times \ZZ}.
\]
For $A$ of finite type, one has $\mathfrak{X}^\flat(\ABT) = \mathfrak{X}(\ABT)$ by \cite[Theorem 3.1]{keller-knitting-algo}.
See also \cite[$\S$7.5]{GHL:friezes}.


\ble{:Y-admi-finite} 
When $A$ is of finite type, $\mathfrak{X}^\flat(\YY_{B})= \mathfrak{X}^{\rm global}(\YY_{B})$,  and one has a  bijection
\[
\lambda: \;\; \mathfrak{X}^{\rm global}(\YY_{B})\longmapsto \fX(\ABT), \;\;  y(i,m) \longmapsto x^\sv(i, m),
\hs (i, m) \in [1, r] \times \ZZ,
\]
\ele

\begin{proof} By \leref{:B-im-column}, $\mathfrak{X}^\flat(\YY_{B})\subset \mathfrak{X}^{\rm global}(\YY_{B})$.
As $\lambda$ is injective and $\lambda|_{\mathfrak{X}^\flat(\YY_{B})}: \mathfrak{X}^\flat(\YY_{B}) \to \mathfrak{X}(\AA_{B^T})$ is surjective,
the statements follow.
\end{proof}

\subsection{Lie theory preliminaries}\label{ss:Lie}
Fix a realization of an Cartan matrix $A$ of finite type, i.e., an $r$-dimensional real vector space $\h$ together with a basis
$\{\alpha_1, \ldots, \alpha_r\}$ of $\h^*$, called the {\it simple roots}, and a basis
$\{\alpha_1^\sv, \ldots, \alpha_r^\sv\} \subset \h$, called the {\it simple co-roots}, such that
$a_{i, j} = (\alpha_i^\sv, \, \alpha_j)$ for $i, j \in [1, r]$.
Let $\{\omega_1, \ldots, \omega_r\}\subset \h^*$ be the fundamental weights, i.e., the
basis of $\h^*$ dual to $\{\alpha_1^\sv, \ldots, \alpha_r^\sv\} \subset \h$,
and let $\calP = {\rm Span}_\ZZ\{\omega_1, \ldots, \omega_r\}\subset \h^*$ be the corresponding weight lattice.

For $i \in [1, r]$, let $s_i = s_{\alpha_i}: \h^* \to \h^*, s_i(\mu) = \mu - ( \alpha_i^\sv, \, \mu)\alpha_i$ for $\mu \in \h^*$.
Let $W$ be the Weyl group generated by $s_1, \ldots, s_r$ and let
$w_0 \in W$ be the longest element.
For $i \in [1, r]$, let $i^* \in [1, r]$ such that $w_0 \omega_i = -\omega_{i^*}$.
Consider the Coxeter element $c = s_1s_2 \cdots s_r \in W$. 
By \cite[Proposition 1.3]{YZ:Lcc}, for each $i \in [1, r]$ there exists a positive integer $h(i; c)$ such that
\begin{equation}\label{eq:omega-list}
\omega_i > c\omega_i > \cdots > c^{h(i; c)-1}\omega_i > c^{h(i; c)} \omega_i = -\omega_{i^*},
\end{equation}
where for $\lambda, \mu \in \calP$, we write $\lambda > \mu$ if $\lambda \neq \mu$ and
$\lambda - \mu \in \sum_{i=1}^r \ZZ_{\geq 0} \alpha_i$.

Recall that $A$ is said to be {\it decomposable} if there exists a partition $[1,r]=I\sqcup J$ of the set $[1,r]$ such that $I,J$ are non-empty and  $a_{i,j} = 0$ whenever $i \in I, j \in J$. In this case, we write 
$A = A_I \times A_J$, where $A_I = (a_{i_1,i_2})_{i_1,i_2 \in I}$ and $A_{J} = (a_{j_1,j_2})_{j_1,j_2 \in J}$. We say that $A$ is 
{\it indecomposable} if $A$ is not decomposable. Note that if $A$ is of finite type and $A = A_I \times A_J$, then both $A_I$ and $A_J$ are 
Cartan matrices of finite type. For $K = I$ or $J$, the triple 
$({\rm Span}_{\mathbb R}(\alpha_k^\sv : k \in K), \{\alpha_k : k \in K\}, \{\alpha_k^\sv : k \in K\} )$ is a realization of 
$A_K$, and the involution $i \mapsto i^*$ on $[1,r]$ restricts to an involution on $K$. Setting $c_K = s_{k_1} s_{k_2} \cdots s_{k_n}$ where 
$k_1 < k_2 < \ldots < k_n$ with $\{k_1, \ldots, k_n\} = K$, it is clear that $c = c_I c_J$ and that 
\begin{equation}\label{eq:hic-decomp}
    h(k;c) = h(k;c_K), \qquad k \in K.
\end{equation} 
It is thus enough to understand the integers $h(i;c)$ when $A$ is {\it indecomposable}. 

Suppose that $A$ is indecomposable, and let ${\bf h}$ be the Coxeter number of $A$, i.e., the order of the element $c \in W$. 
When the Cartan matrix $A$ is of type $A_1, B,C,D_{2n}, n \geq 2,  E_7,E_8, F_4$ or $G_2$, 
the longest element $w_0$ in $W$ acts as the negative of the identity operator on $\h^*$, so $i = i^*$ and 
$h(i; c) = \frac{1}{2} {\bf h}$ for every $i \in [1, r]$. 
In the remaining types of  $A_n, D_{2n+1}$, $n \geq 2$, or $E_6$, for each $i$ the number $h(i; c)$ can be read off from the 
Dynkin quiver defined by $B_A$, and we refer to \cite[Remark 7.23]{GHL:friezes} for the detail.

\subsection{The fundamental domain and periodicity}\label{ss:Ringel-period} Let $A$ be a Cartan matrix of finite type. Set
\begin{equation}\label{eq:F-intro}
{\mathfrak{F}}_A: \;\; [1, r] \times \ZZ \longrightarrow [1, r] \times \ZZ , \;\; {\mathfrak{F}}_A(i,m) = (i^*, \; m+1+h(i^*; c)).
\end{equation}
Note that ${\mathfrak{F}}_A$ is bijective and has fundamental domain
\begin{equation}\label{eq:calD}
\calD_A = \{(i, m): i \in [1, r], \, m \in [0, h(i; c)]\}.
\end{equation}

\bpr{:F-invariant-frieze} All the four generic frieze patterns 
\begin{align*}
&(i, m) \longmapsto x^\sv(i, m) \in \fX(\ABT), \hs (i, m) \longmapsto y(i, m) \in \fX(\YB), \hs (i, m) \in [1, r] \times \ZZ,\\
&(i, m) \longmapsto x(i, m) \in \fX(\AA_{-B}), \hs 
(i, m) \longmapsto y^\sv(i, m) \in \fX(\YY_{-B^T}), \hs (i, m) \in [1, r] \times \ZZ,
\end{align*}
are ${\mathfrak{F}}_A$-invariant and have $\calD_A$ as a fundamental domain.
\epr

\begin{proof}
If $A$ is indecomposable, the statements for $\AA_{-B}$ is proved in \cite[Corollary 7.22]{GHL:friezes} and that 
for $\ABT$ follows from the bijection
in \eqref{eq:x-xvee}. The statements for $\YB$ and $\YY_{-B^T}$ follow from \leref{:Y-admi-finite}. We now describe how the case of arbitrary 
$A$ reduces to the indecomposable one. We do this for $\YB$, the other three statements being entirely parallel. Suppose that $A$ 
is decomposable and $A = A_I \times A_J$. For $K = I$ or $J$, set $B_K = B_{A_K}$ and denote 
the generic $\bfY$-frieze pattern of $\YY_{B_K}$ by $\{y_K(k,m)\}_{(k,m) \in K \times \ZZ}$. Let $\phi_K: \FF(\YY_{B_K}) \rightarrow \FF(\YB)$ be the field embedding defined by
$\phi_K(y_K(k,0)) = y(k,0)$ for $k \in K$.
Using \eqref{eq:y(m,i)-relations} and the fact that $a_{i,j} = 0$ for $i \in I, j \in J$, we have 
$\phi_K(y_K(k,m)) = y(k,m)$ for  $k \in K$ and $m \in \ZZ$. 
By \eqref{eq:hic-decomp}, we have $\calD_A = \calD_I \cup \calD_J$ and 
\[
    \phi_K(y_K(\mathfrak{F}_{A_K}(k,m))) = y(\mathfrak{F}_A(k,m)), \quad k \in K, m \in \ZZ,
\]
where $\mathfrak{F}_{A_K}(k,m) = (k^*, m + 1 + h(k^*;c_K))$ and $\calD_K = \{(k,m) : k \in K, m \in [0, h(k;c_K)]\}$ (see \S \ref{ss:Lie} for the definition of $c_K$). Thus, the generic \bfY-frieze pattern 
of $\YY_B$ is $\mathfrak{F}_A$-periodic with fundamental domain $\calD_A$ if and only if the generic \bfY-frieze pattern 
of $\YY_{B_K}$ is $\mathfrak{F}_{A_K}$-periodic with fundamental domain $\calD_{A_K}$ for $K = I, J$. Since every Cartan matrix of finite type is a 
finite product of indecomposable Cartan matrices of finite type, our claim follows.
\end{proof}

\bre{:FA}
{\rm For the standard Cartan matrix of type 
$A_r$, i.e., $A = (a_{i, j})$ with $a_{i, i} = 2$ for all $i \in [1, r]$, $a_{i, i+1} = a_{i+1, i} = -1$ for all $i \in [1, r-1]$,
and $a_{i, j} = 0$ otherwise, the $\mathfrak{F}_A$-invariance stated in \prref{:F-invariant-frieze} is called
the {\it gliding symmetry} of the generic frieze patterns. }
\hfill $\diamond$
\ere

As a direct consequence of  \prref{:trop-frieze-dual-1}, \prref{:additive-dual-1} and 
\prref{:F-invariant-frieze}, we have the following periodicity of tropical friezes and cluster-additive functions, generalizing 
Ringel's conjecture \cite[Page 477]{Ringel:additive} when $A$ is simply-laced.

\bth{:ringel-conj-period} (Ringel's Conjecture: periodicity) For any Cartan matrix $A$ of finite type, 
all elements in $C_{\rm frieze}(A)$, $\CfriezeAT$, $\CaddA$ and $C_{\rm add}(A^T)$ are 
$\mathfrak{F}_A$-invariant and have $\calD_A$ as a fundamental domain.
\eth

\subsection{Decomposition of cluster-additive functions}\label{ss:Ringel-decomp} Recall from 
$\S$\ref{ss:dual-additive} that each cluster monomial $x^\sv \in \calU^{\rm mono}(\ABT)$ on $\ABT$ 
gives rise to $k_{x^\sv} \in C_{\rm add}(A)$ via
\[
k_{x^\sv}(i, m) =\val_{\delta^\sv(i,m)}(x^\sv) =  (x^\sv(i, m) \, \| \, x^\sv)_d, \hs (i, m) \in [1, r] \times \ZZ.
\]
Since $A$ is of finite type, cluster monomials on $\ABT$ are in bijection with tropical points of $\YB$ via their $g$-vectors, and
by \prref{:additive-dual-2},  we have the bijection
\begin{equation}\label{eq:U-k}
\calU^{\rm mono}(\ABT) \longrightarrow C_{\rm add}(A), \;\; x^\sv \longmapsto k_{x^\sv}.
\end{equation}

We now have the following generalization of 
Ringel's conjecture  
\cite{Ringel:additive} on 
cluster-additive functions for Cartan matrices of finite and simply-laced type.

\bth{:ringel-conj-decomp} (Ringel's Conjecture: decomposition)
Every cluster-additive function associated to any Cartan matrix of finite type 
is a finite sum of cluster-hammock functions.
\eth

\begin{proof}
Let $k \in \CaddA$ and let $k = k_{x^\sv}$ for a unique $x^\sv \in \calU^{\rm mono}(\ABT)$. 
Write
\begin{equation}\label{eq:xsv-lam}
x^\sv =\prod_{j=1}^n (x^\sv_j)^{\lambda_j},
\end{equation}
where $x_j^\sv\in \fX(\AA_{B^T})$ and $\lambda_j \in \ZZ$
is positive for  $j \in [1,n]$. By \leref{:ca-decomposition},
$k = k_{x^\sv} = \sum_{j=1}^n \lambda_j k_{x^\sv_j}$.
As $A$ is of finite type, $x^\sv_j \in \fX^\flat(\ABT)$ for each $j$. By \prref{:hammock-1},
each $k_{x^\sv_j}$ is a cluster-hammock function.
\end{proof}

To get an explicit decomposition of elements in $C_{\rm add}(A)$, 
for $k \in C_{\rm add}(A)$, let 
\begin{equation}\label{eq:DAk}
\calD_A(k) = \{(i,m) \in \calD_A: k(i, m) < 0\}.
\end{equation}

\bpr{:k-decomp}
For any $k \in C_{\rm add}(A)$, one has 
\begin{equation}\label{eq:k-k}
k = \sum_{(i,m) \in \calD_A(k)} (-k(i,m)) h_{(i,m)} = \sum_{(i, m) \in \calD_A} [-k(i, m)]_+ h_{(i, m)}.
\end{equation}
\epr

\begin{proof}
Let $k = k_{x^\sv}$ for a unique $x^\sv \in \calU^{\rm mono}(\ABT)$ as in \eqref{eq:xsv-lam}. Write $x^\sv =\prod_{j=1}^n (x^\sv_j)^{\lambda_j}
$, where $\{x_1^\sv, \ldots, x_n^\sv\}$ is a subset 
of some cluster of $\ABT$ and each $\lambda_j$ is a  positive integer. 
By \leref{:d}, 
$\{x_1^\sv, \ldots, x_n^\sv\} = \{x^\sv(i,m): (i,m) \in \calD_A(k)\}$, 
and $k(i,m) = -\lambda_j$ if
$(i, m) \in \calD_A(k)$ and $x^\sv(i,m) = x_j^\sv$, 
leading to \eqref{eq:k-k}.
\end{proof}

\bre{:trasjective-hammock}
{\rm
When $A$ is not necessarily of finite type, 
following \cite{Assem-Dupont:euclidean}, a cluster variable $x^\sv$ of $\ABT$ is said to be {\it transjective} if 
$x^\sv = x^\sv(i,m)$ for some $(i, m) \in [1, r] \times \ZZ$, and  a cluster monomial is said to be transjective if
it is a monomial of transjective cluster variables from some cluster. The same arguments as for \thref{:ringel-conj-decomp} show that
if $x^\sv$ is a transjective cluster monomial on  $\ABT$, then the cluster-additive function $k_{x^\sv} \in C_{\rm add}(A)$
is a (finite) sum of cluster-hammock functions.
\hfill $\diamond$
}
\ere

\subsection{Fock-Goncharov duality via tropical friezes and cluster-additive functions}\label{ss:thmE}
Let again $A$ be any Cartan matrix  of finite type, let $B = B_A$ as in \eqref{eq:BA-intro}, and consider the Fock-Goncharov dual pair 
$\OO_A = (\ABT, \, \YB)$. As $B$ is of finite type, by \leref{:finite-bijection-1}
we have bijections 
\begin{equation}\label{eq:AY}
\ABTZ \longrightarrow \calU^{\rm mono}(\YB), \;\; \delta^\sv \longmapsto y_{\sdeltasv}, \hs \mbox{and} \hs 
\YBZ \longrightarrow \calU^{\rm mono}(\ABT), \;\; \rho \longrightarrow x^\sv_\rho,
\end{equation}
and by \coref{:duality-pairing-finite-type}  we have the Fock-Goncharov pairing 
\[
\langle \delta^\sv, \; \rho\rangle = \val_{\delta^\sv} (x^\sv_\rho) = \val_\rho(y_{\sdeltasv}), 
\]
where for $\delta^\sv \in \ABTZ$ and $\rho \in \YBZ$, $y_{\sdeltasv} \in \calU^{\rm mono}(\YB)$ and 
$x^\sv_\rho \in \calU^{\rm mono}(\ABT)$ are the respective global monomials 
with $\delta^\sv$ and $\rho$ as their  $g$-vectors. 
We now prove Theorem E in $\S$\ref{ss:finite-intro} which expresses the two maps in \eqref{eq:AY} and the 
Fock-Goncharov pairing via tropical friezes and cluster-additive functions.

\begin{proof}[Proof of 1) and 2) in Theorem E]
The same arguments as in the proof of \prref{:k-decomp} prove 1) of Theorem E, namely for every $\rho \in \YBZ$ one has
\begin{equation}\label{eq:xsv-11}
x^\sv_\rho= \prod_{(i,m) \in \calD_A} x^\sv (i,m)^{[-k_\rho(i,m)]_+}.
\end{equation}

For $\delta^\sv \in \ABTZ$ and $\rho \in \YBZ$, writing $x^\sv_\rho$ as in \eqref{eq:xsv-11}, one has 
\[
\langle \delta^\sv , \rho \rangle = 
{\rm val}_{\delta^\sv}(x^\sv_\rho) = 
\sum_{(i,m) \in \calD_A} [-k_\rho(i,m)]_+ {\rm val}_{\delta^\sv}(x^\sv(i,m)) 
= \sum_{(i,m) \in \calD_A} f_{\sdeltasv}(i,m) [-k_\rho(i,m)]_+,
\]
which proves 2) of Theorem E. 
\end{proof}

 To prove 3) in Theorem E, we first prove a generalization of 2) in Theorem E to the case of principal coefficients.
Recall from $\S$\ref{ss:admissible-Aprin}  the Fock-Goncharov 
dual pair $(\Aprin,\, \; \Yprin)$ with principal coefficients.
As in $\S$\ref{ss:admissible-Aprin}, denote the  clusters of $\Aprin$  as 
$\{\wt{{\bf u}}_t = ({\bf u}_t, {\bf p})\}_{t \in \TT_r}$,
where ${\bf p} = (p_1, \ldots, p_r)$ are the frozen variables. 
Define the generic frieze pattern of $\AA_B^{\rm prin}$ by 
\[
    (i,m) \longmapsto u(i,m) \stackrel{\rm def}{=} u_{t(i,m);i}, \qquad (i,m) \in [1,r] \times \ZZ.
\]

\bpr{:2E-prin}
Let $A$ be a Cartan matrix of finite type, and let $B = B_A$ as in \eqref{eq:BA-intro}. 
If $\widetilde{u}$ is an extended cluster monomial on $\AA_B^{\rm prin}$ with $g$-vector 
$\widetilde{\rho} = \{(\rho^\wedge_t, \mu_t) \in \ZZ^r_{\rm row} \times \ZZ^r_{\rm row}\}_{t \in \TT_r} \in 
\YY_{B^T}^{\rm prin}(\Zmax)$, 
then $\rho^\wedge = \{\rho^\wedge_t\in \ZZ^r_{\rm row}\}_{t \in \TT_r}  \in \YY_{B^T}(\Zmax)$, and one has
\[
    \widetilde{u} = {\bf p}^{-(\mu_{t_0})^T} \prod_{(i,m) \in \calD_A} u(i,m)^{[-k_{-\rho^\sw}(i,m-1)]_+}.
\]
\epr

\begin{proof}
Since $\widetilde{\rho}\in \YY_{B^T}^{\rm prin}(\Zmax)$ and by \leref{:breve-B-p-map}, one has $\rho^\wedge  \in \YY_{B^T}(\Zmax)$.
So $-\rho^\wedge\in \YY_{-B^T}(\Zmax)$ and it defines a cluster-additive function $k_{-\rho^\wedge}$ associated to the Cartan matrix $A^T$.

Let $w \in \TT_r$ be such that 
$\widetilde{u} = {\bf p}^{-(\mu_w)^T} \, {\bf u}_w^{-(\rho_w^\wedge)^T}$, where $-(\rho^\wedge_w)^T\in\mathbb Z_{\geq 0}^r$.
As in $\S$\ref{ss:auxi}, denote the clusters of $\AA_B$ and of $\AA_{-B}$ respectively as $\{{\bf z}_t\}_{t \in \TT_r}$ 
and $\{{\bf x}_t\}_{t \in \TT_r}$ and consider the algebra isomorphism
\[
\varphi: \;\;\calU(\AA_B) \longrightarrow \calU(\AA_{-B}): \;\; {\bf z}_t \longmapsto {\bf x}_t, \hs \forall \; t \in \TT_r.
\]
Let $z = {\bf z}_w^{-(\rho_w^\wedge)^T} \in \calU^{\rm mono}(\AA_B)$.
Then the $g$-vector of $z$ is $\rho^\wedge \in \YY_{B^T}(\Zmax)$. By \leref{:z-x-g-vector}, the $g$-vector of
$x = \varphi(z)={\bf x}_w^{-(\rho_w^\wedge)^T} \in \calU^{\rm mono}(\AA_{-B})$ is $(-\rho^\sw)[1] \in \YY_{-B^T}(\Zmax)$. 
Applying  1) of Theorem E to the Cartan matrix $A^T$ in which case $B_{A^T}=-B^T$, one has
\[
{\bf x}_w^{-(\rho_w^\wedge)^T}=x = \prod_{(i,m) \in \calD_A} x(i,m)^{[-k_{(-\rho^\sw)[1]}(i,m)]_+} = \prod_{(i,m) \in \calD_A} x(i,m)^{[-k_{-\rho^\sw}(i,m-1)]_+}.
\]
It follows that 
\[
       {\bf u}_w^{-(\rho_w^\wedge)^T} = \prod_{(i,m) \in \calD_A} u(i,m)^{[-k_{-\rho^\sw}(i,m-1)]_+}. 
\]
By \leref{:gvec-mono-coeff-t0}, $\mu_{w} = \mu_{t_0}$. Thus $\widetilde{u} = {\bf p}^{-(\mu_w)^T} \, {\bf u}_w^{-(\rho_w^\wedge)^T}={\bf p}^{-(\mu_{t_0})^T}\prod_{(i,m) \in \calD_A} u(i,m)^{[-k_{-\rho^\sw}(i,m-1)]_+}$.
\end{proof}

\begin{proof}[Proof of 3) in Theorem E] Let $\delta^\sv \in \ABTZ$, and let 
$y_{\sdeltasv}$ be the global monomial on $\YY_B$ with $g$-vector $\delta^\sv$. By \prref{:YB-to-Aprin-admi} and \eqref{eq:p-star},
$\widetilde{p}^*(y_{\sdeltasv})$ given in \eqref{eq:p-YB-Aprin} is an extended cluster monomial of $\AA_{B}^{\rm prin}$ whose
extended $g$-vector has coordinates $(\delta^\sv_{t_0}B^T, \delta^\sv_{t_0})$ at $t_0$.  By \prref{:2E-prin}, 
\[
    \widetilde{p}^*(y_{\sdeltasv}) = {\bf p}^{-(\delta^\sv_{t_0})^T} \, \prod_{(i,m) \in \calD_A} \, 
    u(i,m)^{[-k_{-\rho^\sw}(i,m-1)]_+}. 
\]
Let $Z$ be the ring homomorphism $\ZZ[{\bf u}_{t_0}^{\pm 1}, {\bf p}^{\pm 1}] \rightarrow \ZZ[{\bf p}^{\pm 1}]$ by setting each variable in ${\bf u}_{t_0}$ to
$1$. 
Then 
\[
    Z(\widetilde{p}^*(y_{\sdeltasv})) = {\bf p}^{-(\delta^\sv_{t_0})^T} \, \prod_{(i,m) \in \calD_A} \, 
    (F(i,m)({\bf p}))^{[-k_{-\rho^\sw}(i,m-1)]_+}. 
\]
where $F(i,m) \in \ZZ[{\bf p}]$ is the $F$-polynomial of $u(i,m)$. Our claim follows from the fact that 
$Z \circ \widetilde{p}^*: \ZZ[{\bf y}_{t_0}^{\pm 1}] \to \ZZ[{\bf p}^{\pm 1}]$ is a ring isomorphism.  This concludes the proof of Theorem E.
\end{proof}

Recall from \eqref{eq:wt-Bt} that the matrix pattern of $\AA_B^{\rm prin}$ is $\{(B_t \bb C_t) : t \in \TT_r\}$. 
For $(i, m) \in [1, r] \times \ZZ$, let $c(i,m)$ be the $i^{\rm th}$ column of $C_{t(i,m)}$. The following 
recursive relations on the $F$-polynomials is obtained by specializing the exchange relations of the generic frieze pattern 
of $\AA_B^{\rm prin}$ at ${\bf u}_{t_0} = 1$.  

\ble{:F(i,m)}
    The polynomials $\{F(i,m) : i \in [1,r], m \in \ZZ \} \subset \ZZ[{\bf p}]$ are uniquely determined by the initial conditions 
$F(i,0) = 1$ for $i \in [1, r]$ and the recursive relations
\[
F(i,m) \, F(i,m+1) = {\bf p}^{[-c(i,m)]_+} + {\bf p}^{[c(i,m)]_+} \prod_{j=i+1}^r F(j,m)^{-a_{j,i}} \prod_{j=1}^{i-1}F(j,m+1)^{-a_{j,i}}.
\]
\ele

\bre{:Guo}
{\rm
We give some remarks comparing our results with Guo's in \cite{Guo:tropical-frieze} on tropical friezes.
In \cite{Guo:tropical-frieze}, the author defines tropical friezes on any $2$-Calabi-Yau Hom-finite triangulated category 
$\calC$ with cluster-tilting objects. More precisely, for a cluster-tilting object $T$ of $\calC$ and
$m \in K_0({\rm End}_{\calC}(T))$, the author defines a map $f_{T, m}: {\rm Ob}(\calC) \to \ZZ$ and gives a condition in
\cite[Theorem 3.1]{Guo:tropical-frieze} for $f_{T, m}$ to be a tropical frieze on $\calC$. 
As explained in \cite[$\S$3.3]{Guo:tropical-frieze}, under the condition in \cite[Theorem 3.1]{Guo:tropical-frieze}
$f_{T, m}$ is the composition of a cluster character 
${\rm Ob}(\calC) \to \QQ({\bf x})_{>0}$ defined by $T$ and a semi-field homomorphism 
$\QQ({\bf x})_{>0} \to \Zmax$ defined by $m$. In the terminology of our paper, 
$f_{T, m}$ corresponds to the tropical evaluation  map
$x \mapsto \val_{\delta}(x)$ for $\delta \in \AA_{-B}(\Zmax)$ but computed at $w \in \TT_r$ where 
$\delta$ is optimized (see \leref{:val-rho-delta}).
As every $\delta \in \AA_{-B}(\Zmax)$ is optimized at some $w \in \TT_r$ when $A$ is of finite type,
the main Theorem 5.1 of \cite{Guo:tropical-frieze} on tropical friezes, which says that 
every tropical frieze of $\calC$ is of the form $f_{T, m}$ when $\calC$ is the cluster category of a simply-laced Dynkin quiver, 
corresponds to our statement that tropical friezes associated to $A$ are in bijection with tropical points of $\AA_{-B}$. 
\hfill $\diamond$
}
\ere

\bibliographystyle{alpha}
\bibliography{myref}
\end{document}